\documentclass[english]{smfart}

\usepackage[english]{babel}
 \usepackage{amsfonts}
\usepackage{epsfig,graphics}
\usepackage{amssymb}
\usepackage{amscd}
\usepackage{enumerate}
\usepackage{smfthm}

\author{Charles Frances}
\address{IRMA, 7 rue Ren\'e Descartes, 67000 Strasbourg.}
\email{cfrances@math.unistra.fr}
\urladdr{}
\title[Lorentz dynamics on closed $3$-manifolds]{Lorentz dynamics on closed $3$-manifolds}
\newtheorem{theoreme}{Theorem}[section]

\newtheorem{fact}[theoreme]{Fact}

\newtheorem{proposition}[theoreme]{Proposition}
       
\newtheorem{definition}[theoreme]{Definition}

\newtheorem{lemme}[theoreme]{Lemma}

\newtheorem{remarque}[theoreme]{Remark}

\newtheorem{corollaire}[theoreme]{Corollary}

\newtheorem{thmx}{Theorem}
\newtheorem{corox}[thmx]{Corollary}

\newcommand{\Ker}{\operatorname{Ker }}

\newcommand{\Hom}{\operatorname{Hom }}
\newcommand{\iso}{\operatorname{Iso }}

\newcommand{\ad}{\operatorname{ad}}

\newcommand{\SL}{\operatorname{SL}}
\newcommand{\Sol}{\operatorname{SOL}}
\newcommand{\Nor}{\operatorname{Nor}}
\newcommand{\GL}{\operatorname{GL}}
\newcommand{\PSL}{\operatorname{PSL}}
\newcommand{\PGL}{\operatorname{PGL}}
\newcommand{\SO}{\operatorname{SO}}
\newcommand{\OO}{\operatorname{O}}

\newcommand{\ein}{\operatorname{\bf Ein}_3}

\newcommand{\Conf}{\operatorname{Conf}}
\newcommand{\Iso}{\operatorname{Iso}}

\newcommand{\Ad}{\operatorname{Ad}}
\newcommand{\Aut}{\operatorname{Aut}}

\newcommand{\cald}{\mathcal{D}}
\newcommand{\calf}{\mathcal{F}}

\newcommand{\calm}{\mathcal{M}}
\newcommand{\calo}{\mathcal{O}}
\newcommand{\calw}{\mathcal{W}}
\newcommand{\tcalm}{\tilde{\mathcal{M}}}
\newcommand{\tm}{\tilde{{M}}}

\newcommand{\ddt}{\frac{\partial}{\partial t}}

\newcommand{\tX}{\tilde{{X}}}
\newcommand{\tY}{\tilde{{Y}}}
\newcommand{\tZ}{\tilde{{Z}}}
\newcommand{\tcalf}{\tilde{{\mathcal{F}}}}
\newcommand{\tFd}{\tilde{{F}}_{\Delta}}
\newcommand{\tF}{\tilde{{F}}}
\newcommand{\tD}{\tilde{{D}}}

\newcommand{\tth}{\tilde{h}}
\newcommand{\tcalfd}{\tilde{{\mathcal{F}_{\Delta}}}}

\newcommand{\calfd}{{{\mathcal{F}}}_{\Delta}}
\newcommand{\tcald}{\tilde{{\mathcal{D}}}}

\newcommand{\hcalm}{\hat{\mathcal{M}}}

\def\dk{{\mathcal{D}\kappa}}
\def\D{{\Delta}}

\def\NN{\mathbb{N}}
\def\TT{\mathbb{T}}

\def\RR{\mathbb{R}}

\def\ZZ{\mathbb{Z}}

\def\O{\mathcal{O}}
\def\OO{\operatorname{O}}

\newcommand{\mint}{{M^{\rm int}}}
\newcommand{\hmint}{{{\hat M}^{\rm int}}}

\newcommand{\R}{{\bf R}}

\newcommand{\BP}{{\bf P}}

\newcommand{\PO}{{\operatorname{PO}(2,3)}}

\newcommand{\hm}{{\hat{M}}}
\newcommand{\hx}{{\hat{x}}}

\newcommand{\hz}{{\hat{z}}}

\newcommand{\delx}{{\frac{\partial}{\partial x_1}}}
\newcommand{\dely}{{\frac{\partial}{\partial x_2}}}
\newcommand{\delz}{{\frac{\partial}{\partial x_3}}}

\newcommand{\ka}{{\kappa}}

\newcommand{\kil}{\operatorname{{\mathfrak{kill}}}}
\newcommand{\kiloc}{\operatorname{{\mathfrak{kill}^{loc}}}}
\newcommand{\heis}{{\operatorname{\mathfrak{heis}}}}
\newcommand{\Heis}{{\operatorname{{Heis}}}}
\newcommand{\SOL}{{\operatorname{{SOL}}}}

\newcommand{\liek}{{\mathfrak{k}}}

\newcommand{\lien}{{\mathfrak{n}}}

\newcommand{\lieg}{{\mathfrak{g}}}
\newcommand{\lieh}{{\mathfrak{h}}}
\newcommand{\liea}{{\mathfrak{a}}}

\newcommand{\tlien}{\tilde{\mathfrak{n}}}
\newcommand{\liep}{{\mathfrak{p}}}
\newcommand{\lies}{{\mathfrak{s}}}

\newcommand{\oo}{{\mathfrak{o}}}

\newcommand{\sld}{\operatorname{{\mathfrak{sl}}(2,\RR)}}

\newenvironment{preuve}{\medskip \noindent {\bf Proof: }}
   {$\diamondsuit$ }

\begin{document}
\frontmatter

\begin{abstract}
In this paper, we give a complete topological, as well as geometrical classification of 
closed $3$-dimensional Lorentz manifolds admitting a noncompact isometry group.
\end{abstract}
%

\maketitle
\section{Introduction}

A celebrated theorem of Myers and Steenrod \cite{myers-steenrod}, says that the isometry group of a closed Riemannian
 manifold is always a compact Lie  transformation group.  It is very well known that this compactness
  property is specific to the Riemannian world, and fails for general closed pseudo-Riemannian manifolds.  
  For instance the Lorentz torus 
    $  \RR^n /\ZZ^n$, endowed with the metric induced by $-dx_1^2+dx_2^2+\ldots+dx_n^2$,  has
    isometry group $\OO(1,n-1)_{\ZZ} \ltimes \TT^n$.  This group is noncompact, since $\OO(1,n-1)_{\ZZ}$ is a lattice in 
     $\OO(1,n-1)$.  
     
 The geometrical implications of those two antagonistic phenomena --noncompactness of the isometry group on the one hand, and
 compactness of the manifold on the other hand-- were much studied in the  Lorentzian case, on which we will focus here.
  A sample of significant results can be found in \cite{gromov}, \cite{zimmer.iso},  \cite{adams.stuck}, 
  \cite{zeghib.isom},  \cite{zeghibisom1}, \cite{zeghibisom2}, among a lot of other works.
  
  One remarkable point is that the noncompactness of the isometry group is also expected to have strong topological
   consequences.  This was first noticed by M. Gromov in \cite{gromov} when the isometry group is ``large'', for instance
   when it contains a noncompact simple Lie group (see also further developments in \cite{fisher.zimmer}).
     Without any extra asumption on the acting group, let us mention the following striking result:

  \begin{theoreme}\cite{dambra}
   \label{thm.dambra}
   Let $(M,g)$ be a closed Lorentz manifold.  We assume that $M$ and $g$ are real analytic, and $M$ is simply connected.
    Then $\Iso(M,g)$ is a compact  group.
  \end{theoreme}

  The analyticity condition is crucial in the proof of Theorem \ref{thm.dambra}, and 
  when the dimension of the manifold is $> 3$, we actually don't know if  the result holds in 
  the smooth category.

  The aim of this paper is to focus on the  $3$-dimensional situation, and to provide a thorough 
   study of all closed $3$-dimensional
   manifolds, which can be endowed with a Lorentz metric admitting a noncompact group of isometries.  
   
   \subsection{Statement of results}
   
   Let us recall a class of closed $3$-manifolds which will play a prominent role
   in the sequel, namely the torus bundles over the circle 
(torus bundles for short).  Let  $\TT^2$ be   a $2$-torus $ \RR^2/ \ZZ^2 $,
  and let us consider the product $[0,1] \times \TT^2$. We then make  the identification $(0,x) \simeq (1,Ax)$, where $A$ is a given element of 
 $ \SL(2,\ZZ)$.  The resulting $3$-manifold is denoted $\TT_A^3$.  When $A= id$, we just get the $3$-torus $\TT^3$.
   If $A \in \SL(2,\ZZ)$ is {\it hyperbolic}, namely is $\RR$-split with eigenvalues of modulus $\not =1$, we say that 
   $\TT_A^3$  is a 
   {\it hyperbolic torus bundle}.  
 If $A \in \SL(2,\ZZ)$ is {\it parabolic}, namely conjugated to a unipotent matrix ($A \not =id$), we say that $\TT_A^3$
    is a
    {\it parabolic torus bundle}.
     Finally, {\it elliptic  torus bundles} are those for which $A$ has finite order. 
     
  \subsubsection{A topological classification}   
  Our first result is a topological classification of closed Lorentz $3$-manifolds admitting 
  a noncompact isometry group.   
    
\begin{thmx}
 \label{thm.topological}
 Let $(M,g)$ be a smooth, closed $3$-dimensional Lorentz manifold.  Assume that $(M,g)$ is orientable 
 and time-orientable, and that  $\Iso(M,g)$ is noncompact.  Then $M$ is homeomorphic to
  one of the following spaces:
  \begin{enumerate}
   \item {A quotient $\Gamma \backslash \widetilde{\PSL}(2,\RR)$, where 
   $\Gamma \subset \widetilde{PSL}(2,\RR)$ is any uniform lattice.}
   \item{A $3$-torus ${\mathbb T}^3$, or a  torus bundle ${\mathbb T}_A^3$, where $A \in \SL(2,\ZZ)$ can be any 
   hyperbolic or  parabolic element. }
  \end{enumerate}
Conversely, any smooth compact $3$-manifold homeomorphic to one of 
the examples above can be endowed with a smooth Lorentz metric with a noncompact isometry group.
\end{thmx}

The assumption about orientability and time-orientability of the manifold is not really relevant, and one could drop it
 (adding a few allowed topological types) with extra  case-by-case arguments in our proofs.  Notice that any closed $3$-manifold will have
  a covering of order at most four satisfying the assumptions of Theorem \ref{thm.topological}.  We thus see that a lot
   of $3$-manifolds do not admit coverings appearing in the list of the theorem, where only four among the eight Thurston's geometries 
   are represented. Hyperbolic manifolds are notably missing,  and we can  state:
   
   \begin{corollaire}
    \label{coro.hyperbolicmanifold}
    Let $M$ be a smooth closed $3$-dimensional manifold, which is homeomorphic to a complete 
    hyperbolic manifold $\Gamma \backslash {\mathbb H}^3$.  Then for every smooth Lorentz metric $g$ on $M$, the group
     $\Iso(M,g)$  is compact.
   \end{corollaire}

\subsubsection{Continuous versus discrete isometries}
 
 It is interesting to compare the conclusions of Theorem \ref{thm.topological} to closely related results, and especially 
 to the work \cite{zeghibflot}, which was a great source of motivation for the present paper.  In \cite{zeghibflot},
   A. Zeghib studies $3$-dimensional closed manifolds admitting a non equicontinuous  isometric flow.  This hypothesis
    is actually equivalent to {\it the noncompactness of the identity component $\Iso^o(M,g)$}.  The classification can be
     briefly summarized as follows:

\begin{theoreme}\cite[Theorems 1 and 2]{zeghibflot}
 \label{thm.ghaniclassification}
 Let $(M,g)$ be a smooth, closed $3$-dimensional Lorentz manifold.  If the identity component $\Iso^o(M,g)$
  is not compact, then:
  \begin{enumerate}
   \item{Up to a finite cover, the manifold $M$ is homeomorphic either to a torus bundle $\TT_A^3$, with 
   $A \in \SL(2,\ZZ)$ hyperbolic, or to a quotient $\Gamma \backslash \widetilde{\PSL}(2,\RR)$, for a uniform lattice
   $\Gamma \subset \widetilde{\PSL}(2,\RR)$.}
  \item{The manifold $(M,g)$ is locally homogeneous.  It is flat when $M$ is a hyperbolic torus bundle, and locally 
  modelled on a Lorentzian, non-Riemannian, left-invariant metric on $\widetilde{\PSL}(2,\RR)$ otherwise.} 
  \end{enumerate}

\end{theoreme}
The definition of Lorentzian, non-Riemannian, left-invariant metrics on $\widetilde{\PSL}(2,\RR)$ will be made precise in 
   Section \ref{sec.exemplessl2}.
   
Does it make a big difference, putting the noncompactness assumption on  $\Iso^o(M,g)$ instead of $\Iso(M,g)$? 
At the topological level,  notice that $3$-tori and parabolic torus bundles do not show up in Theorem \ref{thm.ghaniclassification}.  
For Lorentz metrics on those manifolds, $\Iso^o(M,g)$ is always compact, but  we will see that there
exist suitable metrics $g$, 
 for which the full group  $\Iso(M,g)$ is noncompact.  It means that for those examples, the noncompactness comes from
  the {\it discrete part} $\Iso(M,g)/\Iso^o(M,g)$.  Actually, there are instances of $3$-manifolds 
  (see Section \ref{sec.panorama}), where the isometry group is discrete, isomorphic to $\ZZ$.
  
  To put more emphasis on  how the general case may differ from the conclusions of \cite{zeghibflot}, let us state
   the following existence result:

\begin{thmx}
 \label{thm.nonhomogeneous}
 Let $M$ be a closed $3$-dimensional manifold which is homeomorphic to a $3$-torus $\TT^3$, or a torus bundle $\TT_A^3$,
  with $A \in \SL(2,\ZZ)$ hyperbolic or parabolic. Then it is possible to endow $M$ with  time-orientable Lorentz metrics 
  $g$ having the following properties:
  \begin{enumerate}
   \item The isometry group $\Iso(M,g)$ is noncompact, but the identity component $\Iso^o(M,g)$ is compact.
   \item There is no open subset of $(M,g)$ which is locally homogeneous.
  \end{enumerate}

\end{thmx}
Observe that for any closed Lorentz manifold $(M,g)$ which is not locally homogeneous, $\Iso^o(M,g)$ is automatically compact by 
Theorem \ref{thm.ghaniclassification} above.

The constructions leading to theorem \ref{thm.nonhomogeneous} are rather flexible.  In particular, on $\TT^3$, or on
 any hyperbolic or pabolic torus bundle ${\TT_A^3}$, the moduli space of Lorentz metrics admitting a noncompact
  isometry group is by no mean finite dimensional.  This is again in sharp constrast  with the second point 
  of Theorem \ref{thm.ghaniclassification}.

\subsubsection{Geometrical results}
The topological classification given by Theorem \ref{thm.topological} comes as a byproduct of a finer, geometrical 
understanding of closed Lorentz $3$-manifolds with noncompact isometry group.  We actually get a quite complete geometrical
 description:

\begin{thmx}
 \label{thm.geometrical}
 Let $(M,g)$ be a smooth, closed $3$-dimensional Lorentz manifold.  Assume that $(M,g)$ is orientable 
 and time-orientable, and that  $\Iso(M,g)$ is noncompact.
 \begin{enumerate}
  \item{If $M$ is homeomorphic to $\Gamma \backslash \widetilde{\PSL}(2,\RR)$, then $(M,g)$ is locally homogeneous, modelled
   on a Lorentzian, non-Riemanniann, left-invariant metric on $\widetilde{\PSL}(2,\RR)$.}
  \item{If $M$  is homeomorphic to $\TT_A^3$, with $A \in \SL(2,\ZZ)$ hyperbolic, then there exists a  
  smooth, positive, periodic 
   function $a:\RR \to (0,\infty)$ such that the universal cover 
  $(\tilde{M}^3,{\tilde{g}})$ is isometric to $\RR^3$ endowed with the metric
  $$ {\tilde{g}}=dt^2+2a(t)dudv.$$
  If $g$ is locally homogeneous, it is flat.}
  \item{If $M$ is homeomorphic to $\TT_A^3$, with $A \in \SL(2,\ZZ)$ parabolic, then there exists a  
  smooth, positive, periodic 
   function $a:\RR \to (0,\infty)$ such that the universal cover 
  $(\tilde{M}^3,{\tilde{g}})$ is isometric to $\RR^3$ endowed with the metric
  $$ {\tilde{g}}=a(v)(dt^2+2dudv).$$
  If $g$ is locally homogeneous, it is either flat or modelled on the Lorentz-Heisenberg geometry.}
  \item{If $M$ is homeomorphic to a  $3$-torus $\TT^3$, then the universal cover $(\tilde{M}^3,{\tilde{g}})$ is
   isometric to $\RR^3$ with a metric of type $2)$ or $3)$ above.  If the metric $g$ is locally homogeneous, it is flat.}
 \end{enumerate}

\end{thmx}
  
 The Lorentz-Heisenberg geometry will be described in Section \ref{sec.lorentzheisenberg}.\\

We already emphasized that in some examples, the isometry group $\Iso(M,g)$  could be infinite discrete. 
However, it is worth mentioning that noncompactness of $\Iso(M,g)$  always produces somehow local continuous 
symmetries.  It is
 indeed easy to infer from Theorem \ref{thm.geometrical} the following result.
 
 \begin{corox}
  \label{coro.action.universel}
  Let $(M,g)$ be a closed $3$-dimensional Lorentz manifold.  If $\Iso(M,g)$ is noncompact, then $\Iso^o(\tm,\tilde{g})$
   is noncompact.  Actually $(\tm,\tilde{g})$  admits an isometric action of the group $\widetilde{\PSL}(2,\RR)$, $\Heis$ or $\Sol$.
 \end{corox}

 \subsection{General strategy of the proof, and organization of the paper}

 One aspect of the present work consists of existence results.  This is the topic of Section \ref{sec.panorama}, where
  we recollect well-known, and probably less known, examples of closed Lorentz $3$-manifolds having a noncompact 
  isometry group.  Examples are given, where $\Iso(M,g)$  is infinite discrete, or semi-discrete.
    This yields the existence part in Theorem \ref{thm.topological}, and a proof of Theorem \ref{thm.nonhomogeneous}.
    
 The remaining of the paper is then devoted to our classification results, namely Theorems \ref{thm.topological} and 
 \ref{thm.geometrical}.   The point of view we adopted, is  that of Gromov's theory of rigid geometric structures 
 \cite{gromov}.
   
   Section \ref{sec.gromov}  recall the main aspects of the theory, recast in the framework of Cartan geometry as in \cite{melnick}, 
   \cite{vincent}.  The key result is the existence of a dense open subset $\mint \subset M$, called the
   {\it integrability locus}, where Killing generators of finite order do integrate into genuine local Killing fields.
    Using the recurrence properties of the isometry group, this implies the crucial fact that the noncompactness of $\Iso(M,g)$ must produce a lot of local Killing fields
     (Proposition \ref{prop.dimension}).  Those continuous local symetries, arising from a potentially discrete 
     $\Iso(M,g)$, will be of great help to understand the geometry of the connected components of $\mint$, which
      can be roughly classified into three categories: constant curvature, hyperbolic, and parabolic 
      (see Section \ref{sec.components}).   To unravel the global structure of $M$, we must understand how all the
      components of
       $\mint$ are patched together (notice that there can be infinitely many such components).
     
     The first, and easiest case to study, is when all the components of $\mint$  are locally homogeneous.  
      Results of \cite{frances.open} show that $(M,g)$ itself is then locally homogeneous, allowing to  
      understand $(M,g)$ completely.  This is done in Section \ref{sec.homogeneous}.
      
     Section \ref{sec.hyperbolic} studies the case where one component of $\mint$ is 
     not locally homogeneous and hyperbolic. One then shows that $M$ is a $3$-torus or a hyperbolic torus bundle, and the 
     geometry is that of examples $2.$ and $4.$ of Theorem \ref{thm.geometrical}.  
     This is summarized in Theorem \ref{theo.hyperbolic}. 
      The key feature in this case is to show that $(M,g)$ contains a Lorentz $2$-torus, on 
      which an element  $h \in \Iso(M,g)$  acts as an Anosov diffeomorphism (Lemma \ref{lem.anosov.tore}).  We then
       show that it is possible to push this Anosov torus by a kind of normal flow, to recover the topological, as well 
       as geometrical structure of $(M,g)$.

     The most tedious case to study is when $(M,g)$ is not
     locally homogeneous, and there are no hyperbolic components at all.  This is the purpose of Sections 
     \ref{sec.local.geometry}, \ref{sec.geomeinstein} and \ref{sec.global}.  We show there that $M$ is a $3$-torus or a 
     parabolic torus bundle, and the geometry is the one described in cases $3.$ and $4.$ of Theorem 
     \ref{thm.geometrical}.  This
      is summarized in Theorem \ref{th.theoplat}.  The main observation here is  that the manifold $(M,g)$
       is conformally flat (Section \ref{sec.local.geometry}).  We then get a developing map 
       $\delta: \tilde{M}^3 \to \ein$, which is a 
       conformal immersion from the universal cover $\tilde{M}^3$ to a Lorentz model space $\ein$, 
       called Einstein's  universe.  After introducing relevant geometric aspects of $\ein$ in 
       Section \ref{sec.geomeinstein}, we are in position to study in details the map $\delta: \tilde{M}^3 \to \ein$ in Section 
       \ref{sec.global}.  We show that $\delta$ maps $\tilde{M}^3$ in a one-to-one way onto an open subset of $\ein$, which is 
        conformally equivalent to Minkowski space.  We are then reduced to the study of closed, flat, Lorentz $3$-manifolds with 
        noncompact isometry groups, which was already done in Section \ref{sec.homogeneous}.
        
        All those partial results are recollected in Section \ref{sec.fin}, where we see how they yield 
        Theorem \ref{thm.geometrical} and Corollary \ref{coro.action.universel}.

\section{A panorama of examples}
\label{sec.panorama}

The aim of this section is to construct a wide range  of closed $3$-dimensional Lorentz manifolds $(M,g)$, with noncompact
isometry group. Those examples will show  that all topologies appearing in Theorem \ref{thm.topological}  do really occur. 
 Moreover, Sections \ref{sec.exemplestorus}, \ref{sec.exempleshyperboliques} and \ref{sec.exemplesparaboliques}  prove
  our Theorem \ref{thm.nonhomogeneous}.  Part of the examples presented here are well known, others like those described in 
   Section \ref{sec.exemplestorus2}, \ref{sec.exempleshyperboliques}, \ref{sec.exemplesparaboliques} seem less classical,
   though elementary.
\subsection{Examples on quotients $\Gamma \backslash \widetilde{\PSL}(2,\RR)$}
\label{sec.exemplessl2}

The Lie group $\widetilde{\PSL}(2,\RR)$, universal cover of ${\PSL}(2,\RR)$, admits a lot of interesting left-invariant
Lorentzian
 metric.  The most symmetric one is the {\it anti-de Sitter} metric $g_{AdS}$.  It is obtained by left-translating 
 the Killing
  form of the Lie algebra $\mathfrak{sl}(2,\RR)$.  The space $(\widetilde{\PSL}(2,\RR),g_{AdS})$ is a complete Lorentz  manifold
   with constant sectional curvature $-1$, called anti-de Sitter space $\widetilde{\bf{AdS}}_3$.  Because the Killing form
  is $\Ad$-invariant, the metric $g_{AdS}$ is invariant by  left and right multiplications of $\widetilde{\PSL}(2,\RR)$
  on itself.  It follows that for any uniform lattice $\Gamma \subset \widetilde{\PSL}(2,\RR)$, the metric $g_{AdS}$
 induces a Lorentz metric $\overline{g}_{AdS}$ on the quotient manifold $\Gamma \backslash \widetilde{\PSL}(2,\RR)$, 
 with a noncompact isometry group coming from the right action of $\widetilde{\PSL}(2,\RR)$ 
 on $\Gamma \backslash \widetilde{\PSL}(2,\RR)$.
 
 There are other metrics than $g_{AdS}$ on $\widetilde{\PSL}(2,\RR)$, which allow the same kind of constructions.  They are
  obtained as follows. Exponentiating
   the linear space spanned by the matrix $\left( \begin{array}{cc} 
                                                0&1\\
                                                0&0\\
                                               \end{array}
 \right)$ (resp. $\left( \begin{array}{cc} 
                                                1&0\\
                                                0&-1\\
                                               \end{array}
 \right)$), one gets a unipotent (resp. $\RR$-split) $1$-parameter group  $\{\tilde{u}^t\}$ (resp. $\{ \tilde{h}^t  \}$) in 
 $\widetilde{\PSL}(2,\RR)$.  The adjoint action of each of those flows, admits invariant
 Lorentz scalar products on $\mathfrak{sl}(2,\RR)$, which are  not equal to  a multiple of the Killing form. One can  left-translate those scalar products 
 and
  get metrics $g_u$ and $g_h$ on $\widetilde{\PSL}(2,\RR)$ which are respectively $\widetilde{\PSL}(2,\RR)
  \times \{\tilde{u}^t\}$ and $\widetilde{\PSL}(2,\RR)
  \times \{\tilde{h}^t\}$-invariant.   Actually there are families of such metrics $g_u$ and $g_h$ which are not pairwise
  isometric.  Now, for each uniform lattice $\Gamma \subset \widetilde{\PSL}(2,\RR)$, the quotient
   $\Gamma \backslash \widetilde{\PSL}(2,\RR)$ can be endowed with induced metrics $\overline{g}_u$ or $\overline{g}_h$ 
   carrying an isometric, noncompact action of $\RR$, coming from the right actions of, respectively,
   $\{\tilde{u}^t\}$ and $\{\tilde{h}^t\}$ on $\widetilde{\PSL}(2,\RR)$.
  
 In the sequel, the metric $g_{AdS}$ and metrics of the form $g_u$ or $g_h$, will be refered to as 
 {\it Lorentzian, non-Riemannian, left-invariant metrics on} $\widetilde{\PSL}(2,\RR)$.  Those are the only left-invariant 
 metrics on $\widetilde{\PSL}(2,\RR)$, the isometry group of which does not preserve a Riemannian metric.

\subsection{Examples on hyperbolic torus bundles}
\label{sec.exempleshyperboliques}

Let us start with the space $\RR^3$ endowed with coordinates $(x_1,x_2,t)$ associated to a basis $(e_1,e_2,e_t)$. We 
 consider a hyperbolic matrix $A$ in $\SL(2,\ZZ)$.  Hyperbolic means that $A$ has two distinct real 
eigenvalues $\lambda$ and $\lambda^{-1}$ different from $\pm 1$.

 Let us  consider the  group $\Gamma$ generated by $\gamma_1=T_{e_1}$ (the translation
 of vector $e_1$), $\gamma_2=T_{e_2}$ and the affine
 transformation
 $\gamma_3=\left(  \begin{array}{cc} A & 0\\ 0 & 1  \end{array} \right)+ \left(  \begin{array}{c}  0\\0\\ 1  \end{array}\right)$.
  It is clear that $\Gamma$ is discrete, acts freely properly and  discontinuously on $\RR^3$, giving a quotient 
   manifold 
  $\Gamma \backslash \RR^3$ diffeomorphic to the hyperbolic torus bundle  ${\mathbb T}_A^3$.

  We see $A$ as a linear transformation of $\operatorname{Span}(e_1,e_2)$.  This transformation is of the form 
  $(u,v) \mapsto (\lambda u , \lambda^{-1}v)$ in suitable coordinates $(u,v)$.
    For any smooth function  $a: \RR \to (0,\infty)$,  which is 
   $1$-periodic, the group $\Gamma$ acts isometrically for the Lorentz metric 
   $$ g_a=dt^2+2a(t)dudv$$
  on $\RR^3$.
Hence the metric $g_a$ induces a Lorentz
  metric 
  $\overline{g}_a$  on $M={\mathbb T}_A^3$.  
  
  When $a$ is a constant, we get for $\overline{g}_a$ a flat metric on ${\mathbb T}_A^3$, and the flow of translations 
  $T_{e_3}^t$ acts on $({\mathbb T}_A^3,\overline{g}_a)$ as an Anosov flow. 
  Up to finite index, the isometry group coincides with this flow.  It is of course noncompact.
  
  More interesting examples arise if one takes a function $a: \RR \to (0,\infty)$ which is $1$-periodic, but is constant 
  on no sub-interval of $\RR$.  Then all Killing fields of $g_a$ must be tangent to the hyperplanes $t=t_0$.  
  There is no nonempty 
  open subset where the metric
   $g_a$ is locally homogeneous.  The same is true for  $\overline{g}_a$.  The linear transformation 
   $\left(  \begin{array}{cc} A & 0\\ 0 & 1  \end{array} \right)$ induces an isometry $f$ of 
   $({\mathbb T}_A^3,\overline{g}_a)$  which preserves individualy the Lorentz tori $t=t_0$ on ${\mathbb T}_A^3$, and
    acts on them by an Anosov diffeomorphism.  The isometry group $\Iso({\mathbb T}_A^3,\overline{g}_a)$ is thus noncompact. 
     This group  virtually coincides with  the subgroup $<f> \simeq \ZZ$ generated by $f$.  It is thus discrete.
     
  These examples prove Theorem \ref{thm.nonhomogeneous} for hyperbolic torus bundles.

\subsection{Examples on parabolic  torus bundles}
\label{sec.exemplesparaboliques0}

\subsubsection{Flat,  or non locally homogeneous examples}
\label{sec.exemplesparaboliques}
We consider now $\RR^3$ with coordinates $(u,t,v)$.     
Let us call $H$ the $3$-dimensional Lie group given by the affine transformations
$$ \left( \begin{array}{ccc} 1& z & -\frac{z^2}{2}\\ 0 & 1 & -z \\ 0& 0 & 1  \end{array}   \right) 
+ \left( \begin{array}{c} r\\s\\ z \end{array} \right)$$
where $r,s,z$ describe $\RR$.  Observe that $H$ is a subgroup isomorphic to the $3$-dimensional Heisenberg group
$\Heis$. The action of $H$ on $\RR^3$ is free and transitive.   Observe also that $H$ acts isometrically for the flat
  Lorentz metric $h_0=dt^2+2dudv$.  Let $a:\RR \to (0,\infty)$ be a smooth function, which is $1$-periodic, 
  and let us consider the metric
   $$ h_a=a(v)(dt^2+2dudv).$$
   When $a$ is not constant, it is no longer true that $h_a$ is $H$-invariant.  But it remains true that $h_a$
    is invariant under the action of the discrete subgroup $\Gamma \subset H$, comprising transformations of the form
  $$ \left( \begin{array}{ccc} 1& m & -\frac{m^2}{2}\\ 0 & 1 & -m \\ 0& 0 & 1  \end{array}  
  \right)  + \left( \begin{array}{c} \frac{n}{2}\\\frac{l}{2}\\ m \end{array} \right)$$
  where $m,n,l$ describe $\ZZ$.  The gluing map between planes $v=0$ and $v=1$ is made by the matrix 
  $A= \left(  \begin{array}{cc} 1& 1\\0 & 1 \end{array} \right)$.  
  Thus the quotient $\Gamma \backslash \RR^3 $ is diffeomorphic to $\TT_A^3$, with $A$ the 
  unipotent matrix above. All parabolic torus bundles can be obtained by considering finite index subgroups of $\Gamma$.
  
  The metric $h_a$ induces a Lorentz metric $\overline{h}_a$ on the parabolic torus bundle
   $\TT_A^3$, and the linear maps  $ B=\left( \begin{array}{ccc} 1& m & -\frac{m^2}{2}\\ 0 & 1 & -m \\ 0& 0 & 1  \end{array}  
  \right) $, $m \in \ZZ$, normalize $\Gamma$, hence  induce a group of isometries in $(\TT_A^3,\overline{h}_a)$.
  It is readily checked that 
   this group  does not have compact closure in $\Iso(\TT_A^3, \overline{h}_a)$.

  We now make the following observation.  Let $X=X_1 \frac{\partial}{\partial u}
  +X_2 \frac{\partial}{\partial t}+X_3 \frac{\partial}{\partial v}$ be a local conformal vector field 
  for the flat metric $h_0$, then 
   $L_Xh_0=\alpha_X h_0$ for a smooth function $\alpha_X$.  Assume that $X_3$ is nonzero on a small open set, the $X$
    will be a Killing field for $h_a$ if and only if 
    \begin{equation}
     \label{equa.kill}
   \frac{a'(v)}{a(v)}=-\frac{\alpha_X(x)}{X_3(x)}
   \end{equation}

   The set of local conformal Killing fields for $h_0$ is finite dimensional, 
   hence for a generic choice of smooth, $1$-periodic $a$,
    the relation (\ref{equa.kill}) won't be satisfied, whatever the conformal vector field $X$ we are considering. 
    It follows that for such a generic set of function, 
     there won't be any open subset of $\TT_A^3$  (resp. of $\TT^3$) where the metric $\overline{h}_a$ will be locally
     homogeneous.
    
     These examples prove Theorem \ref{thm.nonhomogeneous} for parabolic torus bundles.
 
 \subsubsection{Examples modelled on Lorentz-Heisenberg geometry}
 \label{sec.lorentzheisenberg}
 We call  $\heis$ the $3$-dimensional Heisenberg Lie algebra, and $\Heis$ the connected, simply connected, associated 
 Lie group.  Recall that $\heis$ admits a basis $X,Y,Z$, for which
  the only nontrivial bracket relation is $[X,Y]=Z$.  Let $A \in \SL(2,\ZZ)$ be a hyperbolic matrix, and consider
   the automorphism $\varphi$ of $\heis$, which in the basis $X,Y,Z$ writes 
   $\left(  \begin{array}{cc} A & 0\\ 0 & 1  \end{array} \right)$.  It defines an automorphism $\Phi$ of the Lie group
    $\Heis$. 
    
    The matrix $A$ is diagonal in some basis $X',Y'$ of $\operatorname{Span}(X,Y)$, with eigenvalues $\lambda,\lambda^{-1}$.
     The Lorentz scalar product defined by $<X',Y'>=1$, $<Z,Z>=1$, and all other products are zero, can be left-translated on 
      $\Heis$ to give an homogeneous Lorentz metric $g_{LH}$ called {\it the Lorentz-Heisenberg} metric on $\Heis$.
       By construction, $\Phi$ acts isometrically on $(\Heis,g_{LH})$.  We now consider the following lattice in $\Heis$:
       $$ H_{\ZZ}:=\{ \exp(aX+bY+cZ) \ | \ (a,b,c) \in \ZZ^3  \}.$$
       The quotient $H_{\ZZ} \backslash \Heis$ is homeomorphic to a parabolic torus bundle, on which the Lorentz-Heisenberg metric
        induces a metric $\overline{g}_{LH}$.  The automorphism $\Phi$ preserves $H_{\ZZ}$, hence induces an isometry
         $\overline{\Phi}$ on $(H_{\ZZ} \backslash \Heis,\overline{g}_{LH})$, and $\overline{\Phi}$ 
         generates a noncompact group.
 
   \subsection{Some examples on the $3$-torus $\TT^3$}
\label{sec.exemplestorus}

 \subsubsection{A flat example}
 \label{sec.exemplestorus1}
 The most classical example, already mentioned in the introduction, comes from the flat metric
 $$ g_{0}=-du^2+dv^2+dw^2.$$
 
 We call $\OO(1,2)$ the group of linear transformations preserving $g_0$, and we introduce $\Gamma$ the discrete subgroup generated by the translations 
 $T_u,T_v,T_w$ of vectors $u,v,w$.
  The quotient $\Gamma \backslash \RR^3$ inherits an induced (flat) metric $\overline{g}_0$ from $g_0$, and the isometry
   group of $(\TT^3,\overline{g}_0)$ is $\OO(1,2)_{\ZZ} \ltimes \ZZ^3$.  It is noncompact, because $\OO(1,2)_{\ZZ}$ is
    a lattice in $\OO(1,2)_{\ZZ}$.  The identity component is however compact in this case.
    
 \subsubsection{Non locally homogeneous examples}
 \label{sec.exemplestorus2}
 
 These examples are built in  the same way as those of  Sections \ref{sec.exempleshyperboliques} and 
 \ref{sec.exemplesparaboliques}, so that we will be rather sketchy in our description.  
 
 We consider the  metric $g_a$, introduced in Section  \ref{sec.exempleshyperboliques}, 
 for $a: \RR \to (0,\infty)$ a smooth $1$-periodic function.  
 
 The metric $g_a$ is invariant by the discrete group $\Gamma$ generated by the translations of vectors $e_1,e_2$ and $e_t$.
   Hence $g_a$ induces a metric $\overline{g}_a$ on $\TT^3$.  As in Section \ref{sec.exempleshyperboliques}, for generic
    choices of the function $a:\RR \to (0,\infty)$, there is no open set on which $\overline{g}_a$ is  locally homogeneous.  The isometry group
     is then $\ZZ \ltimes \TT^2$ (the $\ZZ$-factor comes from the transformation
     $\left(  \begin{array}{cc} A & 0\\ 0 & 1  \end{array} \right)$, as in \ref{sec.exempleshyperboliques}). 
  

We can also consider the metric $h_a$ introduced in Section \ref{sec.exemplesparaboliques}, and take for $\Gamma$
 the discrete subgroup generated by the translations of vectors $(e_u,e_t,e_v)$.  This yields a metric $\overline{h}_a$
  on $\TT^3=\Gamma \backslash \RR^3$.  For generic choices of the $1$-periodic function $a:\RR \to (0,\infty)$, there
   is no open set on which $\overline{h}_a$ is locally homogeneous, and the  isometry group is noncompact, isomorphic to
    $\ZZ \ltimes \TT^2$.
 
 These examples prove Theorem \ref{thm.nonhomogeneous} for $3$-dimensional tori.
\section{Curvature, recurrence, and local Killing fields}
\label{sec.gromov}
\subsection{Generalized curvature map and integrability locus}
 \label{sec.curvaturemap}

Let us consider $(M,g)$, a smooth Lorentz manifold of dimension $n \geq 2$.  All the material presented below holds actually in the much wider
framework of Cartan geometries, but  we won't need such a generality. 

\subsubsection{Cartan connection associated to the metric}
\label{sec.cartan.connection}
Let $\pi: \hm \to M$ denote the  bundle of orthonormal frames on $\hm$.  This 
 is a principal $\OO(1,n-1)$-bundle over $M$, and it is classical (see \cite{kobanomi}[Chap. IV.2 ])  that  the Levi-Civita connection associated to
  $g$ can be interpreted as an Ehresmann connection $\alpha$ on $\hm$, with values in the Lie algebra $\oo(1,n-1)$.  Let $\theta$ be the soldering 
  form on
   $\hm$, namely the $\RR^n$-valued $1$-form on $\hm$, which to every $\xi \in T_{\hx}\hm$  associates the coordinates of the vector 
   $\pi_*(\xi) \in T_xM$
    in the frame $\hx$.  The sum $\alpha + \theta$ is a $1$-form $\omega: T\hm \to \oo(1,n-1) \ltimes \RR^n$ called the {\it canonical Cartan 
    connection} associated 
    to $(M,g)$ (see \cite[Chap. 6]{sharpe} for a nice introduction to Cartan geometries).  
    
    In the following, we will denote by $\lieg$ the Lie algebra $\oo(1,n-1) \ltimes \RR^n$.
     The Cartan 
    connection $\omega$ satisfies the two crucial properties:
     
     - For every $\hx \in \hm$,
      $\omega_{\hx}: T_{\hx}\hm \to \lieg$ is an isomorphism of vector spaces.
      
      - The form $\omega$ is $\OO(1,n-1)$-equivariant (where $\OO(1,n-1)$ acts on 
      $\lieg$ via the adjoint action).

 \subsubsection{Generalized curvature map}
 \label{sec.generalized.curvature}
  The {\it curvature} of the Cartan connection $\omega$ is a $2$-form $K$ on $\hm$, with values in $\lieg$.  If $X$ and $Y$ are 
 two vector fields on $\hm$, it is given by the relation:
 $$ K(X,Y)=d \omega(X,Y) +[\omega(X),\omega(Y)].$$
 Because $\omega_{\hx}$ establishes an isomorphism between $T_{\hx}\hm$ and $\lieg$ at each point $\hx$ of $\hm$, 
   any $k$-differential form on $\hm$, with values in some vector space ${\mathcal W}$,  can be seen as a map from $\hm$ to 
  $ \Hom( \otimes^k  \lieg,{\mathcal W})$.  This remark applies in particular for the curvature form, yielding  a curvature map 
  $\kappa: \hm \to {\mathcal W}_0$, where the vector space ${\mathcal W}_0$ is $\Hom(\wedge^2(\lieg/\liep);\lieg)$ (the curvature
   is antisymmetric and vanishes when one argument is tangent to the fibers of $\hm$).  
   
  We can now differentiate $\kappa$, getting a map $D \kappa : T\hm \to {\mathcal W}_0$.  Our previous remark allows to  
    see $D\kappa$  as a map  $D \kappa : \hm \to {\mathcal W}_1$, with ${\mathcal W}_1=\Hom(\lieg,{\mathcal W}_0)$. 
    Applying this procedure $r$ times, we define 
    inductively the $r$-derivative of the curvature $D^r \kappa : \hm \to \Hom( \otimes^r  \lieg,{\mathcal W}_r)$ (with ${\mathcal W}_r$ defined inductively
     by ${\mathcal W}_r=\Hom(\lieg, {\mathcal W}_{r-1})$).  
    
    Let us now set $m = \dim G$.  The {\it generalized curvature map} of the Cartan geometry $(M, {\mathcal C})$
    is the map  $\dk =(D\kappa,\ldots,D^{m+1}\kappa)$.  The $\OO(1,n-1)$-module ${\mathcal W}_{m+1}$ will be simply denoted ${\mathcal W}$ 
    in the following.

%
%
%
 \subsubsection{Integrability locus}     
   \label{sec.integrability.locus}
   Let $x \in M$ and $\hx \in \pi^{-1}(x)$. Because $\dk: \hm \to {\mathcal W}$ is a $\OO(1,n-1)$-equivariant map, the rank of $\dk$ at $\hx$ 
     (namely the rank of the linear map $D_{\hx}(\dk):T_{ \hx} \hm  \to { {\mathcal W}}$) does not depend on $\hx$ in a same fiber $\pi^{-1}(x)$. 
   Hence, it makes sense to speak about the rank of $\dk$ at a point $x \in M$. 
    One defines {\it the integrability locus of $M$}, denoted $\mint$, as the set of points $x \in M$ at
     which the rank of $ \dk$ is locally constant. Notice that $ \mint$ is a dense open subset of $M$.
      We define in the same way $ \hmint  \subset  \hm$ as the inverse image $ \pi^{-1}( \mint)$ (which is thus the set of points where the rank of
       $\dk: \hm \to {\mathcal W}$ is locally constant).

\subsection{Integrability theorem and structure of $Is^{loc}$-orbits}

\label{sec.integration}

  For $x \in M$, the $Is^{loc}$-orbit of $x$ is the set of points $y \in M$ such that $y=f(x)$ for some local isometry 
  $f:U\subset M \to V \subset M$.  The $\kiloc$-orbit of $x$ is the set of points $y \in M$ that can be reached
   by flowing along (finitely many)  successive local Killing fields.  
   
   Notice that any local Killing field $X$ on $U \subset M$ lifts to a Killing field (still denoted $X$) on $\hm$, satisfying
    $L_X\omega=0$.  Indeed, local flows of isometries clearly induce local flows on the bundle of orthonormal frames.
     Conversely, local vector fields of $\hm$ such that $L_X\omega=0$ will project on local Killing fields on $M$.
      The same remark holds for local isometries.

   Observe finally that if $X$ is a Killing field on $\hm$  (namely $L_X \omega=0$), then the local flow of $X$ preserves 
   $\dk$, hence $X$ belongs to $\Ker(D_{\hx}\dk)$ at each point.  The integrability theorem below says that the converse is 
   true on $\mint$.

   \begin{theoreme}[Integrability theorem]
    \label{thm.integrabilite}
    Let $(M^n,g)$ be a Lorentz manifold, and $\hmint \subset \hm$ the integrability locus.
    
    \begin{enumerate}
     \item{For every $\hx \in \hmint$, and every $\xi \in Ker(D_{\hx}\dk)$, there exists a local Killing field $X$ around $\hx$ such that 
     $X(\hx)=\xi$.}
     \item{The $Is^{loc}$-orbits in $\mint$ are submanifolds of $\mint$, the connected components of which are  $\kiloc$-orbits.}
    \end{enumerate}

    \end{theoreme}

  The deepest, and most difficult part, of the theorem is the first point.  Such an integrability result as well as the structure of 
  $Is^{loc}$-orbits first appeared in \cite{gromov}.  The results were recast in the framework of Cartan geometry by K. Melnick 
  in the analytic case (see \cite{melnick}).  The reference \cite{vincent} gives an alterative approach for smooth Cartan geometries, 
  leading to the statement of Theorem \ref{thm.integrabilite}.  A proof that the integrability property 
   actually holds on the set where the rank of  $\dk$ is locally constant (first point of the theorem) 
   can be found in Annex A of \cite{frances.open}.

  Let us recall how the second point of Theorem \ref{thm.integrabilite} easily follows from the first one 
  (see also \cite[Sec. 4.3.2]{vincent}).  The generalized curvature map $\dk: \hm \to {\mathcal W}$ is invariant under all local isometries.  
  It follows that  $\hmint$ is invariant as well. Given 
  $\hx \in \hmint$, and $w=\dk(\hx)$, the $Is^{loc}$-orbit  $Is^{loc}(\hx)$ is contained
   in $\dk^{-1}(w) \cap \hmint$.  Now since $\dk$ has locally constant rank on $\hmint$, $\dk^{-1}(w) \cap \hmint$ is a submanifold of $\hmint$, and
    the first point of Theorem \ref{thm.integrabilite}  exactly means that the $\kiloc$-orbit of $\hx$ 
    coincides with the connected component of 
    $\dk^{-1}(w) \cap \hmint$ containing $\hx$, hence is a submanifold on $\hmint$. 
    The set $Is^{loc}(\hx)$ is a union of such connected components, 
    hence a submanifold too.  The point we have to check is that this property remains true when one projects everything on $M$.
     Observe first that the projection of $\dk^{-1}(w) \cap \hmint$ on $M$
     coincides with that of  $\dk^{-1}({\calo }.w) \cap \hmint$, where ${\calo}.w$ stands for the $\OO(1,n-1)$-orbit of $w$ in ${\mathcal W}$.
     Now, using the constancy of $\rm{rank}(\dk)$ on
    $ \dk^{-1}({\calo }.w) \cap \hmint$, the $\OO(1,n-1)$-equivariance of $\dk$, and the fact that $\OO(1,n-1)$-orbits in ${\mathcal W}$ are locally closed, 
    one  shows that 
    $\dk^{-1}({\calo }.w) \cap \hmint$  is  a submanifold of $\hmint$.  By $\OO(1,n-1)$-invariance of this set, its projection on $\mint$
     is a submanifold too.

  \subsection{Components of the integrability locus and $\kiloc$-algebra} 
  \label{sec.intlocus}
  
  For each $x \in M$, there exists a neighborhood $U$ of $x$, so that any local Killing field defined in 
  a neighborhood of $x$ is defined on $U$.  There is thus a good notion of Lie algebra of local Killing fields that
  we denote by  $\kiloc(x)$. 
  
  The integrability locaus $\mint$ splits into a union of connected components $\bigcup {\mathcal M}_i$.  The 
  ${\mathcal M}_i's$ will be just called {\it components} in the sequel. 
  The dimension of  $\kiloc(x)$ (which is finite) can not  decrease locally, while the 
   rank of $\dk$ can not decrease locally. It follows from Theorem \ref{thm.integrabilite} that on each component 
   $\calm$,
   the dimension of the 
   Lie algebra $\kiloc(x)$ is locally constant. Hence  the isomorphism class 
    of $\kiloc(x)$ does not depend on $x \in \calm$, and we will sometimes write $\kiloc({\mathcal M})$
     instead of $\kiloc(x)$.  Actually, on each component $\calm$, everything behaves as if the structure was analytic. 

  For $x \in M$, we consider  $\mathfrak{Is}(x)$,  the isotropy algebra at $x$, namely
          the Lie algebra of local Killing fields defined in a neighborhood of $x$ and vanishing at $x$.
  
  \begin{fact}
   \label{rk.isotropie}
   If $x \in \mint$, then the isotropy algebra $\mathfrak{Is}(x)$ is isomorphic to the Lie algebra of the stabilizer of $\dk(x)$ in $\OO(1,n-1)$.
  \end{fact}

  \begin{preuve}   
    Let us consider
    $\hx \in \hm$  in the fiber of $x$.  Every local Killing field $X$ around $x$ which  vanishes at $x$, lifts to a local Killing field  
     around $\hx$, still denoted $X$, which is vertical at $\hx$.   We call $ev_{\hx}$ the map $X \mapsto \omega(X(\hx))$.  The relation $\varphi_X^t.\, \hx=\hx.e^{tev_{\hx}(X)}$, available for $t$ in a neighborhood of $0$, together with the  invariance of $\dk$ under Killing flows, shows that $ev_{\hx}$ a linear embedding from ${\mathfrak I}_x$ to the Lie algebra of the stabilizer of $\dk(x)$ in $\OO(1,n-1)$ (this map is one-to-one because a local Killing field on $\hm$
        vanishing at a point must be identically zero).  Cartan's formula $L_X=\iota_X \circ d +d \circ \iota_X$ shows that whenever
         $X$ and $Y$ are two Killing fields around $\hx$, the relation  $\omega([X,Y])=K(X,Y)-[\omega(X),\omega(Y)]$ holds.  When $X$ or $Y$ is vertical, 
         $K(X,Y)=0$, proving that $ev_{\hx}$ is an anti-morphism of Lie algebras.  To see that $ev_{\hx}$ is onto, let us consider
          $\{ e^{t \xi} \}_{t \in \RR}$,  a
      $1$-parameter group of $\OO(1,n-1)$ fixing $\dk(\hx)$.  Clearly,  $\xi$ belongs to $\operatorname{Kill}^{\dk}(\hx)$, so that
        by Theorem \ref{thm.integrabilite}, $\omega^{-1}(\xi)$ is the evaluation at $\hx$ of a local Killing field.  
        This Killing field being vertical
         at $\hx$, it corresponds to a local Killing field  of ${\mathfrak{Is}}(x)$.  
  \end{preuve}

 \subsection{Nontrivial recurrence provides nontrivial Killing fields}
\label{sec.recurrence}

We still deal here with $(M^n,g)$ a closed $n$-dimensional Lorentz manifold ($n \geq 2$).  There is, as in Riemannian geometry, a notion of Lorentzian volume, 
which provides  a smooth, $\Iso(M^n,g)$-invariant measure on $M$.  This measure is finite under our asumption that $M$ is closed.
 When the group $\Iso(M^n,g)$ is noncompact, Poincar\'e's recurrence theorem applies and almost every point of $M$ is recurrent for 
 the action of $\Iso(M^n,g)$.  We are going to see that such a recurrence phenomenon is responsible for the existence of nontrivial continuous 
 local symetries.  The precise statement is:  

\begin{proposition}
\label{prop.dimension}

Let $(M^n,g)$ be a closed, $n$-dimensional  Lorentz manifold, and assume that $\Iso(M^n,g)$ is noncompact. 
Then 
\begin{enumerate}
 \item{For
 every $x \in \mint$, the isotropy
algebra $\mathfrak{Is}(x)$ generates a noncompact subgroup of ${\operatorname O}(T_xM)$.}
\item{For every component ${\mathcal M} \subset \mint$, the Lie algebra $\kil({\mathcal M})$ is at least $3$-dimensional.}
\end{enumerate}

\end{proposition}

 \begin{preuve}
 The proof of the first point is already contained in   \cite[Proposition 5.1]{frances.open}. 
 We summarize here the main arguments for the reader's convenience. Let $x$ be a recurrent point for $\Iso(M^n,g)$  and choose
 $\hx \in \hm$   in the fiber of $x$.  The recurrence hypothesis means that 
    there exists $(f_k)$ 
 tending to infinity in $\Iso(M^n,g)$, and $(p_k)$ a sequence of ${\operatorname O}(1,n-1)$ such that $f_k(\hx).p_k^{-1}$ tends to $\hx$. 
 By equivariance of the generalized curvature map $\dk: \hm \to \mathcal{W}$, we also have
 $$p_k.\dk(\hx) \to \dk(\hx).$$
 Observe that $(p_k)$ tends to infinity in ${\operatorname O}(1,n-1)$, because $\Iso(M^n,g)$ acts properly on $\hm$. 
 
    The ${\operatorname O}(1,n-1)$-orbits on $\mathcal{W}$ are locally closed, 
   because the action of ${\operatorname O}(1,n-1)$  is linear, hence algebraic. As a consequence, there exists 
   a sequence
   $(\epsilon_k)$ in ${\operatorname O}(1,n-1)$ with $\epsilon_k \to Id$ and  $\epsilon_k.p_k.\dk(\hx) = \dk(\hx)$.
  Since $(p_k)$ tends to infiny by properness of the action of $\Iso(M^n,g)$ on $\hm$,  so does  $(\epsilon_k.p_k)$, proving that the stabilizer
   $I_{\hx}$ of $\dk(\hx)$ in ${\operatorname O}(1,n-1)$ is noncompact. This group is  algebraic, hence the identity component
   $I_{\hx}^{o}$ is noncompact too. Fact \ref{rk.isotropie} then ensures that  
  $\mathfrak{Is}(x)$ 
 generates a noncompact subgroup of ${\operatorname O}(T_xM)$ (under the identification of $\mathfrak{Is}(x)$ with a subalgebra of $\mathfrak{o}(T_xM)$ under the
 isotropy representation), for  every 
 recurrent point $x \in \mint$.  The property is thus  true everywhere on $\mint$, by density of recurrent points on  $M$.
 
 To prove the second point, we start with  $x \in \mint$, and consider a simply connected neighborhood $U \subset \mint$ of $x$. Then, 
 every algebra  $\mathfrak{Is}(y)$, $y \in U$, 
 is realized as a Lie algebra of Killing 
fields defined on $U$.  The first point of proposition \ref{prop.dimension} says that there exists  $X$ a nontrivial Killing field on $U$,
such that $X(x)=0$.
The zero locus of $X$  is a nowhere dense set in $U$. We can thus pick $y \in U$ satisfying $X(y) \not =0$, and apply again the first point of 
the proof at $y$. We get a second nontrivial Killing 
field $Y$ defined on $U$, and vanishing at $y$. Let us now pick  a lightlike direction $u \in T_yM$ such that
$g_y(u,X(y)) \not =0$, and $u$ is transverse to the zero locus of $Y$.  Let $t \mapsto \gamma_u(t)$ be
the geodesic passing through $y$ at $t =0$ and such that $\dot{\gamma}_u(0)=u$.  
Then Clairault's equation ensures that $g(\dot{\gamma}_u,X)$ is constant on $\gamma_u$, hence  does 
not vanish on $\gamma_u$ by our choice of $u$.

On the other hand, by the same Clairault's equation, one has $g_{\gamma_u}(\dot{\gamma}_u,Y) =0$, 
and $Y(\gamma_u(t)) \not =0$ for small, nonzero,  values of $t$.  Hence, on some open subset of $U$, 
the orbits of the local Killing algebra  have dimension $\geq 2$, while for  every point $z \in U$, 
the dimension of  $\mathfrak{Is}(z)$ is $\geq 1$. It follows that the dimension of $\mathfrak{kill}^{loc}(U)$ 
is at least $3$.
 \end{preuve}

 \begin{remarque}
  \label{rem.pseudoRiemannian}
  The proof above does not use the fact that the metric $g$ is Lorentzian, and Proposition \ref{prop.dimension} actually holds for any
   pseudo-Riemannian (non Riemannian) metric.
 \end{remarque}

 \subsection{Components of the integrability locus, and their classification}
 
 \label{sec.components}
 We now stick to dimension $3$, and we consider a closed Lorentz manifold $(M,g)$.  We assume that the isometry 
 group $\Iso(M,g)$ is noncompact. 
 
 On each component ${\mathcal M} \subset \mint$, the discussion of Section \ref{sec.integration} shows that 
 there is a well-defined Lie algebra $\kiloc({\mathcal M})$ 
 of local Killing fields (beware that some monodromy phenomena may occur), which by Proposition \ref{prop.dimension} 
 is at least
  $3$ dimensional.  
   Because in a finite dimensional linear representation of 
   ${\operatorname O}(1,2)$, no point can have a stabilizer of dimension exactly $2$, and because of Fact \ref{rk.isotropie}, the dimension
    of $\mathfrak{Is}(x)$ is $1$ or $3$ for  every $x \in {\mathcal M}$, and the dimension  of $\kiloc({\mathcal M})$ is $6$, $4$ or $3$. 

     -  When this dimension is $6$, the component has constant sectional curvature.

     - When the dimension of $\kiloc({\mathcal M})$ is $4$, it is 
     not hard to check that $\calm$ is locally homogeneous (see for instance 
    \cite[Lemma 4]{sorin.karin}). The dimension of $\mathfrak{Is}(x)$ is then $1$ at each point $x \in \calm$.

    - When the dimension of $\kiloc({\mathcal M})$ is $3$, then the dimension of $\mathfrak{Is}(x)$ is $1$ or $3$ at each point $x \in \calm$.
    The $\kiloc$-orbits have dimension $0$ or $2$, and the component is  nowhere locally homogeneous.

    Proposition \ref{prop.dimension}
    ensures that whenever the dimension of 
    $\mathfrak{Is}(x)$ is $1$, then this algebra generates a hyperbolic or a parabolic flow in ${\operatorname O}(T_xM) \simeq {\operatorname O}(1,2)$.  In the first case, we say that $x$ is {\it 
    a hyperbolic point}, and in the second one we call $x$ {\it a parabolic point}.
    
    \begin{definition}[Hyperbolic and parabolic components]
     A component ${\mathcal M}$ of $\mint$ which is not of constant curvature  is  said to be hyperbolic when it contains a hyperbolic point. 
     Otherwise, it is called parabolic.
    \end{definition}
Observe that this defintion allows {\it a priori} a hyperbolic component to contain parabolic points 
(it will turn out later that it does not occur).
    
    To summarize, components of $\mint$ split into three (rough) categories.  
    \begin{enumerate}[a)]
    \item{The first category comprises all components having {\it constant 
    sectional curvature}. }
    \item{The second category comprises {\it hyperbolic components}.  Those in turn split into two subgategories:
    \begin{enumerate}[i)]
     \item The locally homogeneous ones,  for which 
     the dimension of $\kiloc({\mathcal M})$ is $4$.
     \item The non locally homogeneous ones,   for which 
     the dimension of $\kiloc({\mathcal M})$ is $3$.
     \end{enumerate}}
     \item{The remaining components are {\it  parabolic}. They can also
      be splitted into:
      \begin{enumerate}[i)]
      \item The locally homogeneous ones  for which $dim(\kiloc({\mathcal M}))=4$.
      \item The non locally homogeneous ones for which $dim(\kiloc({\mathcal M}))=3$. 
      \end{enumerate}}
    \end{enumerate}
    
    Let us notice that points belonging to a component with constant curvature are those for which the rank of $\dk$ is locally equal to $0$. 
    In the same way, points belonging to locally homogeneous components are those for which the rank of $\dk$ is locally equal to $2$.
     Finally, we prove:
    
    \begin{lemme}
     \label{lem.rank3}
     Points $x \in M$  belonging to a component of $\mint$ which is not locally homogeneous are exactly those at which the rank of $\dk$ is $3$.
     
    \end{lemme}

  \begin{preuve}
  Recall the generalized curvature map $\dk : \hm \to {\mathcal W}$ introduced in Section \ref{sec.integration}. We saw 
 in  Proposition \ref{prop.dimension} that for every $x \in \mint$, ${\mathfrak{kill}}^{loc}(x)$ has dimension $\geq 3$. Because $\dk$ is
 invariant along $\kiloc$-orbits in $\hm$, the corank of
  $\dk$ is a least $3$ on $\hmint$, and because $\hm$ has dimension $6$, the rank of $\hm$ is at most $3$ on the dense set $\hmint$, hence on $\hm$.
   The rank can only increase locally, 
   hence points where the rank of $\dk$ is $3$ actually  stay  in $\mint$.
    
  \end{preuve}

\section{Locally homogeneous Lorentz manifolds with noncompact isometry group}
\label{sec.homogeneous}

In this section, we prove Theorem \ref{thm.topological}  in the  case where all the components of $\mint$ are locally homogeneous, implying that
$(M,g)$ is locally homogeneous on a dense open set.
  We observed in the previous section that in this case, and under our standing assumption  that $\Iso(M,g)$ is noncompact,  the Lie algebra of 
  local Killing vector fields
  has dimension 
  $\geq 4$ on each components.  We can  apply the results of \cite{frances.open}, 
  saying that we must then have $\mint=M$, and the manifold $(M,g)$
   is locally homogeneous.
   
   \begin{theoreme}[\cite{frances.open}, Theorem B]
    \label{thm.quasihomogeneous}
    Let $(M,g)$ be a smooth $3$-dimensional Lorentz manifold.  
    Assume that on a dense open subset, the Lie algebra of local Killing fields is at least $4$-dimensional.
     Then $(M,g)$ is locally homogeneous.
   \end{theoreme}

There are a lot of homogeneous $3$-dimensional models for Lorentz manifolds.  Fortunately, very few of them can appear as 
the local geometry of a closed manifold with a noncompact isometry group.  We have indeed:

\begin{theoreme}\cite[Theorem 2.1]{sorin.ghani}
 \label{thm.ghanisorin}
 Let $(M,g)$ be a closed locally homogeneous Lorentz manifold.  Assume that at each point $x \in M$ the isotropy
  algebra $\mathfrak{Is}(x)$ generates a noncompact subgroup of $\OO(T_xM)$.  Then the metric $g$ is locally isometric to:
  \begin{enumerate}
   \item{A flat metric}
   \item{A  Lorentzian, non-Riemannian, left-invariant metric on $\widetilde{\PSL}(2,\RR)$.}
   \item{The Lorentz-Heisenberg metric $g_{LH}$ on the group $\Heis$.}
   \item{The Lorentz-Sol metric $g_{\operatorname{sol}}$ on the group $\SOL$.}
  \end{enumerate}

\end{theoreme}
The Theorem applies in our situation since Proposition \ref{prop.dimension} ensures that $\mathfrak{Is}(x)$  generates 
 a noncompact subgroup of $\OO(T_xM)$ for almost every (hence every) point $x$.
 
Lorentzian, non-Riemannian, left-invariant metrics $g_{AdS},g_u$ and $g_h$ on $\widetilde{\PSL}(2,\RR)$ were introduced in 
Section \ref{sec.exemplessl2}, while Lorentz-Heisenberg geometry was described in Section \ref{sec.lorentzheisenberg}.
 It remains to explain what is the Lorentz-Sol geometry.
 
 The Lie algebra $\mathfrak{sol}$ is the $3$-dimensional Lie algebra with basis $T,X,Y$ and nontrivial 
 bracket relations $[T,Z]=Z$, $[T,X]=-X$.  The corresponding connected, simply connected Lie group is denoted $\SOL$.
   On $\mathfrak{sol}$, we can consider the Lorentz scalar product such that $<T,Z>=<X,X>=1$, and all 
   other products are $0$.  After left-translating this scalar product on $\SOL$, we get a Lorentz metric
   $g_{\operatorname{sol}}$ on $\SOL$ which is called the Lorentz-Sol metric.  The isometry group of 
   $(\Sol,g_{\operatorname{sol}})$ contains $\Sol$  (acting by left-translations), but it is actually $4$-dimensional. The Lie 
   algebra of Killing fields is obtained by adding $Y$ to $T,X,Z$, with bracket relations $[T,Y]=2Y$ and $[X,Y]=Z$   
   (see \cite[Section 4.2]{sorin.ghani} for further
    details). 

\subsection{Ruling out Lorentz-$\SOL$ geometry}

We first establish:
\begin{proposition}
 \label{prop.compacite.sol}
 Let $(M,g)$ be a closed, $3$-dimensional Lorentz manifold locally modelled on  $(\SOL,g_{\operatorname{sol}})$.  Then $\Iso(M,g)$
  is a compact group.
\end{proposition}

The key property for proving Proposition \ref{prop.compacite.sol} is a  Bieberbach rigidity theorem for closed
manifolds modelled on 
 Lorentz-Sol geometry.
\begin{theoreme}\cite[Theorem 1.2 (iv), and Proof of Proposition 7.1]{sorin.ghani}
 \label{thm.bieberbach.sol}
 Let $(M,g)$ be a closed, $3$-dimensional Lorentz manifold locally modelled on  $(\SOL,g_{\operatorname{sol}})$.  Then $(M,g)$ 
 is isometric to the quotient of $(\SOL,g_{\operatorname{sol}})$ by a discrete subgroup 
   $\Gamma \subset \Iso(\SOL,g_{\operatorname{sol}})$.  Moreover, the intersection of $\Gamma$ with 
   $\Iso^o(\SOL,g_{\operatorname{sol}})$ is a lattice
    $\Gamma_0 \subset \SOL$ acting by left translations.
\end{theoreme}

 We know that $(M,g)$ is isometric to some quotient $\Gamma \backslash \Sol$, by Theorem \ref{thm.bieberbach.sol}.
 Let us denote by $L_{\SOL}$ the subgroup of $\Iso(\Sol,g_{\operatorname{sol}})$  comprising all left-translations by elements of $\Sol$.
   The group $\Gamma_0=\Gamma \cap \Iso^o(\Sol,g_{\operatorname{sol}})$  is Zariski-dense in $L_{\SOL}$ by Theorem 
   \ref{thm.bieberbach.sol}.
    It follows that $\Nor(\Gamma)$, the normalizer of $\Gamma$ in $\Iso(\Sol,g_{\operatorname{sol}})$, must normalize  $L_{\SOL}$.
     But the description of $\mathfrak{Iso}(\Sol,g_{\operatorname{sol}})$  made above shows  that $L_{\SOL}$ has finite index in its 
     normalizer.  We infer that $\Nor(\Gamma)/\Gamma$ is compact, what proves Proposition \ref{prop.compacite.sol}.

\subsection{Minkowski and  Lorentz-Heisenberg geometries with noncompact isometry group}
\label{sec.minkowski}

We now focus on closed Lorentz manifolds locally modelled on Minkowski space, or on Lorentz-Heisenberg geometry.  In both cases, 
 one has a Bieberbach type theorem.  This is very well known in the flat case, thanks to the works \cite{fried.goldman} and 
 \cite{goldman.kami}, and the completeness result of Carri\`ere for closed flat Lorentz manifolds \cite{carriere}.  For manifolds modelled on Lorentz-Heisenberg geometry, this is proved in 
 \cite[Proposition 8.1]{sorin.ghani}.  The precise statement is the following:
\newpage
 
\begin{theoreme}[Bieberbach's theorem for flat and Lorentz-Heisenberg manifolds]
 \label{thm.bieberbach.minkowski}
 Let $(M,g)$ be a closed, $3$-dimensional, Lorentz manifold. 
 \begin{enumerate}
 \item If $(M,g)$ is flat, there exists a discrete subgroup 
 $\Gamma \subset \Iso(\RR^{1,2})$ such that $(M,g)$ is isometric to the quotient $\Gamma \backslash \RR^{1,2}$.  Moreover, 
 there exists a connected $3$-dimensional Lie group $G \subset \Iso(\RR^{1,2})$, which is isometric to 
 $\RR^3$, $\Heis$ or $\Sol$, and which acts simply transitively on $\RR^{1,2}$, satisfying that $\Gamma_0=G \cap \Gamma$
  has finite index in $\Gamma$ and is a uniform lattice in $G$.
  \item If  $(M,g)$ is locally modelled on Lorentz-Heisenberg geometry, then it is
  isometric to the quotient of $(\Heis,g_{\Heis})$ by a discrete subgroup 
   $\Gamma \subset \Iso(\Heis,g_{\Heis})$.  Moreover, there exists a finite index subgroup
    $\Gamma_0 \subset \Gamma$  which  is a lattice
    $\Gamma_0 \subset \Heis$ acting by left translations.
\end{enumerate}
  \end{theoreme}

This theorem says that {\it up to finite cover}, a closed Lorentz manifold $(M,g)$ modelled on Minkowski, or Lorentz-Heisenberg 
geometry, is homeomorphic to $\TT^3$ or to $\TT_A^3$ for $A \in \SL(2,\ZZ)$ hyperbolic or parabolic. Noncompactness of the isometry
 group allows to be more precise.

\begin{proposition}
 \label{prop.topologie.minkowski}
 Let $(M,g)$ be a closed, $3$-dimensional Lorentz manifold, such that $\Iso(M,g)$ is noncompact.  We assume that 
 $(M,g)$ is orientable and time-orientable.
  \begin{enumerate}[i)]
 \item{ If $(M,g)$  is flat, then $M$ is diffeomorphic either to a torus $\TT^3$, or to a torus bundle $\TT_A^3$  with 
  $A \subset \SL(2,\ZZ)$  hyperbolic, or parabolic.}
\item{If $(M,g)$ is modelled on Lorentz-Heisenberg geometry, then $M$ is diffeomorphic to a torus bundle $\TT_A^3$  
with $A \in \SL(2,\ZZ)$ parabolic ($A \not =id$).}
 \end{enumerate}
\end{proposition}

\begin{preuve}
 The situation provided by Theorem \ref{thm.bieberbach.minkowski}  is the following 
 (both in the flat and Lorentz-Heisenberg case).  
 We have a $3$-dimensional Lie group $G$, which is either $\RR^3$, $\Heis$ or $\Sol$, as well as a left-invariant
 metric $\mu$ on $G$, and the manifold $(M,g)$ is isometric to a quotient of $(G,\mu)$ by a discrete subgroup 
 $\Gamma \subset \Iso(G,\mu)$.  Moreover, if we denote by $L_G$ the group of left-translations by elements of $G$, the 
  intersection $\Gamma_0=\Gamma \cap L_G$ has finite index in $\Gamma$, and is a uniform lattice in $L_G$.
  
  An important remark is that if $\Nor(\Gamma)$ denotes the normalizer of $\Gamma$ in $\Iso(G,\mu)$, then $\Nor(\Gamma)$
    normalizes $L_G$.  It is obvious in the case of Lorentz-Heisenberg geometry  $(G,\mu)=(\Heis,g_{\Heis})$. 
    In this case, the identity component $\Iso^o(G,\mu)$ is of the form $\RR \ltimes L_G$, and $L_G$ is thus 
       normalized by the full isometry group $\Iso(G,\mu)$.
     
     In the case of Minkowski geometry, one has to remember that the group $G$ (more accurately $L_G$) is the identity component of the 
     {\it crystallographic hull} of $\Gamma$ (see \cite[Section 1.4]{fried.goldman}).  The last part of 
     \cite[Theorem 1.4]{fried.goldman} ensures that in the case of Minkowski geometry, the crystallographic hull is unique.
      It follows that $\Nor(\Gamma)$ must normalize this crystallographic hull, as well as its identity component $L_G$.
      
      As a consequence, elements of $\Nor(\Gamma)$ (in particular elements of $\Gamma$) belong to the group
      $\Aut(G)_{\mu} \ltimes L_G$, where $\Aut(G)_{\mu}$
       denotes the automorphisms of $G$ preserving the metric $\mu$.  Let us  denote by $\Gamma_l$  the projection of
        $\Gamma$ on $\Aut(G)_{\mu}$.  If this projection is trivial, we get that $\Gamma \subset L_G$. The manifold 
        $M$ is obtained as a quotient of $\RR^3$, 
 $ \Sol$ or $ \Heis$ by a uniform lattice, and we are done.

        If this projection is nontrivial, we are going to get a contradiction.  Indeed, $\Gamma_l$
        must be a finite subgroup of $\Aut(G)_{\mu}$, because $\Gamma_0 \subset L_G$
          has finite index in $\Gamma$.  The isotropy representation $\rho$ of $\Aut(G)_{\mu}$ at $e$ identifies $\Gamma_l$
           with a nontrivial finite subgroup of $\SO(\mu)$, which is actually in $\SO^o(\mu)$ because $(M,g)$ is orientable and
           time-orientable.  Let $\Nor^o(\Gamma)$ be the intersection of $\Nor(\Gamma)$ with 
           $\Aut^o(G)_{\mu} \ltimes L_G$, and let $N_l$ the projection of $\Nor^o(\Gamma)$ on $\Aut^o(G)_{\mu}$.
             We get that $\rho(N_l)$ normalizes $\rho(\Gamma_l)$.  Now, nontrivial finite groups in $\SO^o(\mu)$ have a unique 
             fixed point in ${\mathbb H}^2$, so their normalizer in $\SO^o(\mu)$ are contained in a compact subgroup.
               It follows that $\Nor^o(\Gamma)$ is contained in a subgroup $K \ltimes L_G$, with $K$ compact in 
               $\Aut^o(G)_{\mu}$.  We thus see that $\Nor^o(\Gamma)/\Gamma_0$ is compact, implying the compactness of 
               $\Nor(\Gamma)/ \Gamma$. This in turns implies $\Iso(M,g)$ compact: Contradiction.
 
\end{preuve}

\subsection{Anti-de Sitter structures with noncompact isometry group}

It remains to study closed Lorentz manifolds $(M,g)$ modelled on a Lorentzian, non-Riemannian, left-invariant metric on 
 $\widetilde{\PSL}(2,\RR)$.  We observe that the identity components  $\Iso^o(\widetilde{\PSL}(2,\RR),g_u)$ and $\Iso^o(\widetilde{\PSL}(2,\RR),g_h)$ are actually included in
 $\Iso^o(\widetilde{\PSL}(2,\RR),g_{AdS})$.  Thus, if $(M,g)$ is a closed, orientable and time-orientable, Lorentz manifold modelled
  on $(\widetilde{\PSL}(2,\RR),g_u)$ or $(\widetilde{\PSL}(2,\RR),g_h)$, there exists an anti-de Sitter metric $g'$ on $M$ which is preserved by a finite
   index subgroup of $\Iso(M,g)$.  Hence, it will be enough for us to focus on the topology of closed anti-de Sitter
   manifolds with noncompact isometry group.  In the sequel, we will denote $\widetilde{\bf AdS}_3$  the space $(\widetilde{\PSL}(2,\RR),g_{AdS})$ 
   
   \begin{proposition}
    \label{prop.topologieads}
    Let $(M,g)$ be a closed, orientable and time-orientable, anti-de Sitter manifold. If $\Iso(M,g)$ is noncompact,
     then $M$ is homeomorphic to a quotient $\Gamma \backslash \widetilde{\PSL}(2,\RR)$, for a uniform lattice $\Gamma \subset \widetilde{\PSL}(2,\RR)$.
   \end{proposition}

It is worth noticing that all closed (orientable and time-orientable) $3$-dimensional anti-de Sitter manifolds are
 Seifert fiber bundles over hyperbolic orbifolds (a short proof of this fact can be found in 
 \cite[Corollary 4.3.6]{tholozan}).  Conversely, any Seifert fiber bundle over an hyperbolic orbifold, with nonzero Euler number, can
  be endowed with an anti-de Sitter metric (see \cite{scott}).  The  assumption that $\Iso(M,g)$
  is noncompact reduces 
   the possibilities for the allowed Seifert bundles.  For instance,  all nontrivial circle bundles 
   over a closed orientable 
   surface of genus $g \geq 2$ admit 
   anti-de Sitter metrics, but only those
    for which the Euler number divides $2g-2$ do occur in Proposition \ref{prop.topologieads}.
    
 It was shown in \cite{klingler} that closed anti-de Sitter manifolds are complete.  It follows that $(M,g)$ as in Proposition 
 \ref{prop.topologieads}  is a quotient of $\widetilde{\bf AdS}_3$ by a discrete subgroup $\tilde{\Gamma} \subset \Iso(\widetilde{\bf AdS}_3)$. Actually
   $\tilde{\Gamma} \subset \Iso^o(\widetilde{\bf AdS}_3)$ because $(M,g)$ is orientable and time-orientable.
 The center of 
  $\widetilde{\PSL}(2,\RR)$ is infinite cyclic, generated by an element $\xi$.   The group $\widetilde{\PSL}(2,\RR) \times 
  \widetilde{\PSL}(2,\RR)$ acts on ${\widetilde{\bf AdS}_3}$ by left and right translations: $(h_1,h_2).g=h_1gh_2^{-1}$.  This yields
   an epimorphism $\widetilde{\PSL}(2,\RR) \times 
  \widetilde{\PSL}(2,\RR) \to \Iso^o({\widetilde{\bf AdS}_3})$, with infinite cyclic kernel generated by $(\xi,\xi)$. The group $\Iso^o({\widetilde{\bf AdS}_3})$
   has a center $Z$ which is generated by the left action of $\xi$. Doing the quotient of $\Iso^o({\widetilde{\bf AdS}_3})$ by $Z$ yields an epimorphism
    $\pi: \Iso^o({\widetilde{\bf AdS}_3}) \to \PSL(2,\RR) \times \PSL(2,\RR)$.  Notice that $\PSL(2,\RR) \times \PSL(2,\RR)$ coincides with 
     the identity component of the isometries of $\PSL(2,\RR)$ endowed with its anti-de Sitter metric.

   An important result, known as {\it finiteness of level}, says that $\tilde{\Gamma} \cap Z \not = id$. This was first stated
    in \cite{kulkarni.raymond}.  A detailed proof can be found in \cite[Theorem 3.3.2.3]{saleinthese}.  Geometrically, this 
     theorem ensures that there exists a finite group of isometries $\Lambda \subset \Iso(M,g)$, which acts freely and 
     centralizes a finite index subgroup of $\Iso(M,g)$, such that the quotient manifold of $(M,g)$ by $\Lambda$
      is a quotient of $\PSL(2,\RR)$ by a discrete group $\overline{\Gamma} \subset \PSL(2,\RR) \times \PSL(2,\RR)$.
       Let us denote by $(\overline{M}^3,\overline{g})$ this new Lorentz manifold, and observe that 
       $(\overline{M}^3,\overline{g})$ is still noncompact. Observe also that the projection $\pi$ maps $\tilde{\Gamma}$
        onto $\overline{\Gamma}$.
        
        The structure of the group $\overline{\Gamma}$ is well 
       understood.  Up to conjugacy, there exists $\Gamma_0$ a uniform lattice in $\PSL(2,\RR)$, and a representation
        $\rho: \Gamma_0 \to \PSL(2,\RR)$ such that
        $$ \overline{\Gamma}=\{ (\gamma,\rho(\gamma)) \in \PSL(2,\RR) \times \PSL(2,\RR) \ | \ \gamma \in \Gamma_0 \}.$$
   This was established in \cite[Theorem 5.2]{kulkarni.raymond} when $\Gamma$ is torsion-free.  For a group with torsion, the adapted
    proof can be found in \cite[Lemma 4.3.1]{tholozan}.  
    
    Because the group $\Iso(\overline{M}^3,\overline{g})$ is noncompact, a result of Zeghib 
    (\cite[Theorem 1.2]{zeghibisom1}) ensures that $(\overline{M}^3,\overline{g})$ must admit a codimension one, lightlike, 
    totally geodesic foliation.  Such foliations ${\mathcal F}$ in $\PSL(2,\RR)$ endowed with the anti-de Sitter 
    metric are very well known.  Let $AG \subset \PSL(2,\RR)$ be the connected $2$-dimensional group corresponding to 
    the upper-triangular matrices (it is isomorphic to the affine group of the line). Then up to conjugacy in 
    $\PSL(2,\RR) \times \PSL(2,\RR)$, the leaves of ${\mathcal F}$ are given by $\{ gAG \ | \ g \in \PSL(2,\RR)\}$
     or $\{ AGg \ | \ g \in \PSL(2,\RR)\}$.  If $\overline{\Gamma}$ preserves such a foliation ${\mathcal F}$, we infer that
      the leaves are of the form $\{ gAG \ | \ g \in \PSL(2,\RR)\}$, and $\rho(\Gamma_0)$ normalizes $AG$, namely 
      $\rho(\Gamma_0) \subset AG$.  
      
    We now consider the normalizer $H$ of $\overline{\Gamma}$ in  $\PSL(2,\RR) \times \PSL(2,\RR)$.  The first projection $\pi_1(H)$
     must normalize $\Gamma_0$.  Since uniform lattices in $\PSL(2,\RR)$ are of finite index in their normalizer, we can 
     replace $H$ by a finite index subgroup and assume $\pi_1(H)=\Gamma_0$. Let us consider $h=(h_1,h_2)$ in 
      $H$, and $\gamma \in \Gamma_0$. Because $H$ normalizes $\overline{\Gamma}$,
      $(h_1 \gamma h_1^{-1}, h_2 \rho(\gamma)h_2^{-1}) \in \overline{\Gamma}$, which implies 
      $h_2 \rho(\gamma)h_2^{-1}=\rho(h_1)\rho(\gamma) \rho(h_1)^{-1}$. In other words, $h_2^{-1}\rho(h_1)$ centralizes 
      $\rho(\Gamma)$.  As a consequence, $\rho(\Gamma)$ can not be Zariski dense in $AG$.  Otherwise, $h_2^{-1}\rho(h_1)$
       should be trivial implying that $h=(h_1,h_2)$ actually belongs to $\overline{\Gamma}$.  We thus would get that 
       $\Iso(\overline{M}^3,\overline{g})$ is finite, a contradiction.
       
    As a result, $\rho(\Gamma)$ is included in a $1$-parameter subgroup of $AG$.  This implies that the group 
    $\tilde{\Gamma}$ is included in a product $\widetilde{\PSL}(2,\RR) \times \RR$, where $\widetilde{\PSL}(2,\RR)$ acts by 
    left translations, and $\RR \subset \widetilde{\PSL}(2,\RR)$ is a $\RR$-split or unipotent parameter group acting
     on the right.  We consider the projection $\pi_1: \tilde{\Gamma} \to \widetilde{\PSL}(2,\RR)$ on the left-factor.
      The group $\Gamma:=\pi_1(\tilde{\Gamma})$  projects surjectively on $\overline{\Gamma}$ by $\pi$, hence is a uniform lattice
       in $\widetilde{\PSL}(2,\RR)$.  Moreover, the kernel of $\pi_1$ must be trivial, otherwise some nontrivial element
        of $\tilde{\Gamma}$ would belong to $\{id \} \times \RR$, and $\Nor(\tilde{\Gamma})/\tilde{\Gamma}$ would be compact, contradicting
         the hypothesis $\Iso(M,g)$ noncompact. It follows that $\tilde{\Gamma}$ is isomorphic to $\Gamma$.
         
         In conclusion,  the two manifolds 
         $\Gamma \backslash \widetilde{\PSL}(2,\RR)$ and $\tilde{\Gamma} \backslash \widetilde{\PSL}(2,\RR)$ are two
          Seifert bundles over the hyperbolic orbifold $\Gamma_0 \backslash {\mathbb H}^2$.  Their Euler number is nonzero,
           so that
           they are both large Seifert manifolds (see \cite[p. 92]{orlik}).  Their fundamental
           groups are isomorphic (to $\Gamma$), hence by 
            \cite[Theorem 6, p 97]{orlik}  they are homeomorphic.  

       
    \subsection{Conclusions}
    The previous results show that closed, orientable and time-orientable, Lorentz $3$-dimensional manifolds which are locally
     homogeneous are homeomorphic to $\TT^3$, a torus bundle $\TT_A^3$ for $A$ hyperbolic or parabolic, or a quotient
      $\Gamma \backslash \widetilde{\PSL}(2,\RR)$.  This proves Theorem \ref{thm.topological} for manifolds such that
       all components of $\mint$ are locally homogeneous.  Our analysis shows moreover  that when $M$ is homeomorphic to an hyperbolic 
        torus bundle or to a $3$-torus, the metric must be flat.  When $M$ is homeomorphic to a parabolic torus bundle, the metric
         is either flat, or locally isometric to Lorentz-Heisenberg geometry.  When $M$ is homeomorphic to 
         $\Gamma \backslash \widetilde{\PSL}(2,\RR)$, the geometry is locally anti-de Sitter or locally modelled on a Lorentzian,
         non-Riemannian, left-invariant metric on $\widetilde{\PSL}(2,\RR)$.  This is in accordance to points $2,3,4$ of
          Theorem \ref{thm.geometrical}.

     

\section{Manifolds admitting a hyperbolic component} 
\label{sec.hyperbolic}

Our aim in this section is to prove Theorem  \ref{thm.topological} under the assumption that 
our Lorentz manifold $(M,g)$ admits  at least 
one hyperbolic 
component ${\mathcal M}$. Actually Theorem  \ref{thm.topological} will be implied by a more precise description
 provided by Theorem \ref{theo.hyperbolic} to be stated below.

\subsection{A first reduction}

\begin{proposition}
 \label{prop.hyperbolic}
 Assume that  the intergrability locus $\mint$ contains a hyperbolic component.  Then either $(M,g)$ is locally homogeneous, or
 there exists a hyperbolic component which is not locally homogeneous.
\end{proposition}

\begin{preuve}
 Assume  that all hyperbolic components in $\mint$ are locally homogeneous, and consider $\calm \subset \mint$ such a component.  
   It is open by definition, and we are going to show that the boundary $\partial \calm$ is empty, which will yield $\calm=M$ by connectedness of $M$. 
   Local homogeneity of $M$ will follow.
   
   Let us assume that  there exists $x \in \partial \calm$.  We are going to see that $x \in \mint$, which will yield a contradiction. 
   We pick $x_0 \in \calm$.  In the sequel, we use the notation $\dk(z)$ for the $\OO(1,2)$-orbit of $\dk(\hz)$ in ${\mathcal W}$  
   (where $\hz$ is any point in the fiber of $z$). 
   By assumption, the rank of $\dk$ is constant equal to $2$ on $\calm$.  Moreover, by local homogeneity, $\dk({\calm})$ 
   is a single $\OO(1,2)$-orbit in ${\mathcal W}$.  This orbit is $2$-dimensional because $\calm$ is locally homogeneous,
   but does not have constant curvature. The stabilizer of points in  $\dk(x_0)$  are hyperbolic $1$-parameter subgroups of $\OO(1,2)$
    by assumption that ${\mathcal W}$ is hyperbolic.  
  In every finite-dimensional representation of $\operatorname{O}(1,2)$, such hyperbolic orbits 
            are closed.  
            This is a standard fact, 
            the proof of which is, for instance, detailed in \cite[Annex B]{frances.open}. 
            It follows that $\dk(x)=\dk(x_0)$.  Hyperbolic $1$-parameter groups are open in the set of $1$-parameter groups of $\OO(1,2)$.
             It follows that there is a sufficiently small neighborhood $U$ of $x$ in $M$ such that the rank of $\dk$ on $U$ is $\geq 2$, 
             and for every $y \in U$, the orbit  $\dk(y)$ is $2$-dimensional
              with hyperbolic $1$-parameter groups as stabilizers of points (notice that $\dk(y)$ can not be $3$-dimensional because of 
              the first point of
               \ref{prop.dimension}).  If at some point $y \in U$, the rank of $\dk$ is $3$, 
              then $y$ belongs to a component of $\mint$  which is not locally homogeneous, by Lemma \ref{lem.rank3}. This component must 
              be hyperbolic because stabilizers in $\dk(y)$ are hyperbolic. Since we assumed that there are no such components, it follows
               that the rank of $\dk$ is constant equal to $2$ on $U$.  But then $x \in \mint$  leading to the desired contradiction.
               
\end{preuve}

The case of locally homogeneous manifolds was already settled in Section \ref{sec.homogeneous}, so that we will 
assume {\it in all the remaining of this section} that $\mint$ contains a component  which is hyperbolic but 
not locally homogeneous. We will call $\calm$ this component.
We are going to show that under these circumstances, $M$ is diffeomorphic to a hyperbolic  torus bundle.
More precisely,  the  geometry of $(M,g)$ can be described as follows:  
\begin{theoreme}
 \label{theo.hyperbolic}
 Assume that $(M,g)$ is a closed, orientable and time-orientable $3$-dimensional Lorentz manifold, 
 such that $\Iso(M,g)$ is noncompact. Assume that $(M,g)$ admits 
 a hyperbolic component which is not locally homogeneous.  Then
 \begin{enumerate}
  \item{The manifold $M$ is diffeomorphic to a $3$ torus ${\mathbb T}^3$, or a  
  torus bundle ${\mathbb T}_A^3$ where $A \in \SL(2,\ZZ)$ is  
  a hyperbolic matrix.}
  \item{The universal cover $(\tilde{M},\tilde{g})$ is isometric to $\RR^3$ endowed with the metric $dt^2+2a(t)dudv$ for 
  some positive nonvanishing, periodic,  smooth function 
  $a: \RR \to (0,+\infty)$.}
  \item{There is an isometric action of the Lie group $\Sol$ on $(\tilde{M},\tilde{g})$.}
 \end{enumerate}

\end{theoreme}

The proof of Theorem \ref{theo.hyperbolic}, will be the aim of Sections \ref{sec.anosovtori} to \ref{sec.conclusionhyperbolic}
 below.

\subsection{Existence of Anosov tori}
\label{sec.anosovtori}
 We first prove that under the hypotheses of Theorem \ref{theo.hyperbolic}, one can find an element $h$ in $\Iso(M,g)$ which acts 
 by an Anosov transformation of a flat
  Lorentz torus on $M$ (see Lemma \ref{lem.anosov.tore}).  
 \subsubsection{Facts about flat Lorentz surfaces}
 
 We begin by recalling elementary and well known facts about closed flat Lorentz surfaces. 
 Since Lorentz manifolds must have zero Euler characteristic, such a 
 surfaces are tori or Klein bottles. 
 \begin{lemme}
  \label{lem.lorentzsurface}
  A closed, flat Lorentz surface $(\Sigma,g)$ admitting a noncompact isometry group is a torus.  Moreover, any group $H \subset \Iso(\Sigma,g)$
   which does not have compact closure, contains an element $h$ acting on $\Sigma$ by a hyperbolic linear transformation.
 \end{lemme}

 \begin{preuve}
  Closed Lorentz manifolds with constant curvature are geodesically complete (\cite{carriere}, \cite{klingler}).  It follows that a closed, flat
  Lorentz
   surface $(\Sigma, g)$ is a quotient of the Minkowski plane $\RR^{1,1}$, by a discrete subgroup $\Gamma \subset \operatorname{O}(1,1) \ltimes \R^2$, which acts 
    freely and properly on $\RR^{1,1}$. Observe that nontrivial elements of $\operatorname{O}(1,1) $ are of two kinds.  Either they are hyperbolic 
    (namely have two real eigenvalues of modulus $\not = 1$), or they have order $2$ (an orthogonal symmetry with respect to a spacelike 
    (resp. timelike) line). It is readily checked that $\Gamma$ is either a lattice in $\R^2$, or admits a subgroup of index $2$ 
     which is such a lattice.  Assume we are in the first case, and let $H \subset \Iso(\Sigma,g)$ be a subgroup
   which does not have compact closure. We lift $H$ to $\tilde{H} \subset \operatorname{O}(1,1) \ltimes \R^2$.  The group $\tilde{H}$ has a nontrivial projection on 
    $\operatorname{SO}(1,1)$, otherwise $H$ would be compact.  Thus $\tilde{H}$ contains a conjugate of a hyperbolic element of 
    $\operatorname{O}(1,1)$, which acts as an Anosov diffeomorphism on $\Sigma$.
    
    In the second case, where the projection of $\Gamma$ on $\operatorname{O}(1,1)$ is an order $2$ subgroup,  the normalizer $Nor(\Gamma)$ must have
     trivial projection on $\operatorname{SO}(1,1)$, which implies that $\Gamma$ is cocompact in $Nor(\Gamma)$. This shows that flat
     Klein bottles have compact
      isometry group.
 \end{preuve}

 \subsubsection{Closed $\kiloc$-orbits}
 Our next aim is to exhibit some $\kiloc$-orbits which are closed surfaces.  
 
 \begin{lemme}
  \label{lem.anosov.tore}
  In every hyperbolic   component ${\mathcal M}$, there exists a $\kiloc$-orbit $\Sigma_0$ which 
  is a flat Lorentz $2$-torus, and such that
    there exists $h \in \Iso(M,g)$ leaving $\Sigma_0$ invariant, and acting on $\Sigma_0$ as a linear hyperbolic automorphism. 
   
 \end{lemme}


       We consider our distinguished component ${\mathcal M}$ that, we recall, is not locally homogeneous, and hyperbolic. 
       We  pick $x \in \calm$, and $\hx \in {\hat{\calm}}$  in the fiber of $x$.
        We already observed in the proof of Lemma \ref{lem.rank3} that the rank of $\dk$ is at most $3$ on $M$. 
        It is exactly $3$ at $\hx$, still by
          Lemma \ref{lem.rank3}, hence remains constant equal to $3$ 
          in a neighborhood of  $\hx$. Hence, if $U \subset \hcalm$ is a 
       small open set around
        $\hx$, $\dk(U)$ is a $3$-dimensional submanifold of ${\mathcal W}$.
         If $U$ is chosen small enough,  the 
        $\operatorname{O}(1,2)$-orbit of every point in $\dk(U)$ will be $2$-dimensional and will have hyperbolic $1$-parameter 
        groups as stabilizers of points. 
        Let us now call $\hat{\Lambda}$ the closed subset of $\hm$
         where the rank of $\dk$ is $ \leq2$. By Sard's theorem, the $3$-dimensional Hausdorff measure of $\dk(\hat{\Lambda})$ is zero.
         We infer the existence of
          $w \in \dk(U) \setminus \dk(\hat{\Lambda})$. Moving $\hx$ inside $U$, we assume that $w=\dk(\hx)$, and we denote by $\calo(w)$ the $\OO(1,2)$-orbit
           of $w$ in $\calw$.  
         By $\operatorname{O}(1,2)$-equivariance of $\dk$, the inverse image $\dk^{-1}({\mathcal O}(w))$ avoids $\hat{\Lambda}$, hence 
         the rank of $\dk$ is constant equal to $3$ on
         $\dk^{-1}({\mathcal O}(w))$. Lemma \ref{lem.rank3} then leads to the inclusion $\dk^{-1}({\mathcal O}(w)) \subset \hmint$.
         By the discussion right after Theorem \ref{thm.integrabilite}, the projection of $\dk^{-1}({\mathcal O}(w))$ on $M$ is a submanifold 
         $N$ of
          $M$. The stabilizer of
            $w$  in $\operatorname{O}(1,2)$ is hyperbolic, thus
            as mentioned in the proof of Lemma \ref{lem.rank3}, the orbit $\calo(w)$ is closed in $\calw$.
             It follows that $N$ is closed in $M$, hence compact.
           By (the proof of ) Theorem \ref{thm.integrabilite}, the $Is^{loc}$-orbit of $x$ is a union of connected components of $N$, and 
            the connected component of $x$ in $N$, denoted $\Sigma_{0}$, coincides with the $\kiloc$-orbit of $x$. 
            It is a connected compact surface in $M$. 
            Let us  show that this surface has Lorentz signature. 
            The Lie algebra $\mathfrak{Is}(x)$ is generated by a
local Killing 
field $X$ around $x$, vanishing at $x$, and  such that the flow $\{D_{x}\phi_X^t\} 
\subset \operatorname{O}(T_{x}M)$ is a hyperbolic flow. Linearizing $X$ around $x$ thanks to the exponential map, 
we see there are two distinct 
lightlike directions $u$ and $v$ in ${T}_{x}M$ such 
that the two geodesics $\gamma_u: s \mapsto \exp(x,su)$ 
and $\gamma_v: s \mapsto \exp(x,sv)$ are left invariant by $\phi_X^t$.
 In particular,  for $s \not =0$ close to $0$, $\dot{\gamma}_u(s)$ and 
$\dot{\gamma}_v(s)$ are 
colinear to $X$, hence  tangent to ${\mathcal O}(\gamma_u(s))$ and ${\mathcal O}(\gamma_v(s))$ respectively.  
By continuity, this property must still hold for $s=0$. 
We infer that ${T}_{x}({\mathcal O}(x))$ contains the two distinct lightlike 
directions $u$ and $v$, hence has Lorentz signature. By local homogeneity of the $\kiloc$-orbit $\Sigma_0$, we get that 
  $\Sigma_0$ is Lorentz, and moreover has
   constant Gauss curvature.  The only closed Lorentz surfaces of constant curvature are flat tori or Klein bottles.
             
             Now $\Iso(M,g)$ sends $\Sigma_0$ to components of the $Is^{loc}$-orbit of $x$,
             and there  are finitely many such components by compactness of $N$.  As a consequence
              the subgroup $H \subset \Iso(M,g)$ leaving $\Sigma_0$ invariant is noncompact. 
              Observe that  if $g_0$ is the metric induced by $g$ on 
              $\Sigma_0$, then the injection $H \to \Iso(\Sigma_0,g_0)$ is proper (see for instance \cite[Prop. 3.6]{zeghibflot}). 
              It follows that 
              $\Iso(\Sigma_0,g_0)$ is a noncompact group.  Lemma \ref{lem.lorentzsurface}  ensures that $(\Sigma_0,g_0)$ is 
              a flat Lorentz torus, and there exists $h \in \Iso(M,g)$ acting on $\Sigma_0$ by a hyperbolic linear automorphism.


 \subsection{Pushing  Anosov tori along the normal flow}
 \label{sec.anosov.mapping}
 From the $2$-torus   $\Sigma_0$ and the  diffeomorphism $h \in \Iso(M,g)$ given by Lemma \ref{lem.anosov.tore}, we are going to recover 
 the  topology of the whole manifold $M$, as well as its geometry.

 \subsubsection{Preliminary definitions}
 \label{sec.push}
 On the torus $\Sigma_0$, we choose a frame field $(E^-,E^+)$ with the property that $E^-$ and $E^+$ are future lightlike, satisfy $g(E^-,E^+)=1$, and 
 generate the strong stable and 
 unstable bundles of the Anosov diffeomorphism $h$.  Because $M$ is assumed to be orientable, this defines a smooth normal
 field $\nu: \Sigma_0 \to T\Sigma_0^{\perp}$ with the property that $(E^-,E^+,\nu)$ is a direct frame of $T_zM$ at each point  $z \in \Sigma_0$, 
 and $g(\nu,\nu)=+1$.

   In all the rest of the section, we pick once for all $z_0 \in \Sigma_0$  a periodic point of $h$ (recall that the set of periodic points is dense in $\Sigma_0$).
   This point has period $m_0$, and replacing $h$ by $h^{2m_0}$ if necessary, {\it we will assume henceforth that
 $h(z_0)=z_0$ and $h^* \nu =\nu$.}
 
 For every $z \in \Sigma_0,$  we will call $\gamma_z$ the oriented geodesic arc through $z$, with tangent $\nu(z)$ at $z$. 
  Observe that $\gamma_{z_0}$ is a closed spacelike geodesic.  Indeed, since $h$ is a Lorentz isometry, the  fixed points set ${\rm Fix}(h)$
  is  a  closed, totally geodesic submanifold of $M$. 
 The matrix of the differential  $D_{z_0}h$, expressed in the basis  $(E^-(z_0),E^+(z_0), \nu(z_0))$,
  is of the form 
 $\left( \begin{array}{ccc} \frac{1}{\lambda_0}&0&0\\ 0& \lambda_0&0\\ 0&0&1 \end{array} \right)$ with $|\lambda_0|>1$.  Linearizing  $h$ 
 around $z_0$ thanks to the exponential map,
  we see that  the component of  ${\rm Fix}(h)$ containing $z_0$,  is precisely $\gamma_{z_0}$.

 \subsubsection{The normal flow, and an auxiliary pseudo-Riemannian manifold}
 \label{sec.normalflow}
It will be usefull in the sequel  to consider the manifold $N=\RR \times \Sigma_0$.  On this manifold, we have the vector field
 $\frac{\partial}{\partial t}$.  Pushing the vector fields $E^-,E^+$ on $\{ 0 \} \times \Sigma_0$ by the flow of 
 $\frac{\partial}{\partial t}$, we get two more vector fields $\tilde{E}^-,\tilde{E}^+$ on $N$.  The frame field
  $(\tilde{E}^-,\tilde{E}^+,\frac{\partial}{\partial t})$ provides $N$ with an orientation.

Let us consider the map $f: (t,z) \mapsto \exp(z,t\nu(z))$.
It is well-defined and 
smooth  on some maximal open subset
$U_{max} \subset N$. An easy application of the inverse mapping theorem shows that $f: (-\epsilon,\epsilon) \times \Sigma_0 \to M$
is a one-to-one immersion  for small  $\epsilon >0$.

 A key property of the map $f$ is its equivariance with respect to the action of $h$, namely :
 \begin{equation}
  \label{eq.equivariance}
 h \circ f(t,z) = f(t,h(z)),
 \end{equation}
  which is available for $(t,z) \in U_{max}$ (observe that $f(U_{max})$ is left invariant by $h$).  Relation (\ref{eq.equivariance}) just follows 
  from the fact that $h$ is an isometry preserving the normal field $\nu$.  
 
 In the following, we are going to introduce  
$$\tau_m := \sup\{ s \in (0,\infty) \ | \  f:(0,s) \times \Sigma_0 \subset U_{max} \to M 
{{\rm \ is \ an \ injective }{\rm  \ immersion}}\}.$$
 It will be sometimes more suggestive  to restrict $f$ to $\{0\} \times \Sigma_0$, and consider 
 {\it the normal flow of $\Sigma_0$}, $\phi^t: \Sigma_0 \to M$ defined by 
 $\phi^t(z):=\exp(z,t \nu(z))$.     By what we said before, 
 $\phi^t$ is at least defined on $(0,\tau_m)$, and for all $t \in (0,\tau_m)$, $\phi^t: \Sigma_0 \to M$ is a proper embedding, with image 
 $\Sigma_t \subset M$.  Equivariance relation (\ref{eq.equivariance}) shows that $h$ preserves  $\Sigma_t$ and acts on it as an Anosov diffeomorphism.
   In particular, the stable and unstable bundles must be lightlike for the metric $g$, showing that $\Sigma_t$ is a Lorentz torus.  We denote by $g_t$
    the restriction of $g$ to $\Sigma_t$. 

 The map $t \mapsto \phi^t(z_0)$ provides a (cyclic) parametrization of the closed geodesic $\gamma_{z_0}$ at speed $+1$.
 Hence for  every $t \in \RR$, the map $z \mapsto \phi^t(z)$ is defined and smooth on some small neighborhood ${\mathcal U_t} \subset \Sigma_0$
  containing 
  $z_0$.  We call ${\mathcal E}^{\pm}(t):=D_{z_0}\phi^t(E^{\pm}(z_0))$, and  
   the formula
  $a(t):=g_{\gamma_{z_0}(t)}({\mathcal E}^{-}(t),{\mathcal E}^{+}(t))$ defines a smooth function $a: \RR \to \RR$. 
  This in turns defines on  the manifold $N$ a
  symmetric $(2,0)$-tensor $\tilde{g}=dt^2+a(t)g_{0}$  (here $g_0$ is the metric induced by $g$ on $\Sigma_0$).  Observe that 
  $\tilde{g}$ is not a metric on $N$ because $a(t)$ might vanish for some values of $t$.

  \subsubsection{First return time}
  
  For any point  $z \in \Sigma_0$, the geodesic $\gamma_z$ is defined on some maximal interval $[0,\tau_z^*)$.  If 
 $\gamma_z((0,\tau_z^*)) \cap \Sigma_0 \not = \emptyset$, then there exists a smallest $\tau(z) \in (0,\tau_z^*)$ such that 
 $\gamma_z(\tau(z)) \in \Sigma_0$.  When $\gamma_z((0,\tau_z^*)) \cap \Sigma_0  = \emptyset$, we just put $\tau(z)=+ \infty$.
   We introduce {\it the first return set:}
 $$\Omega_r= \{  z  \in \Sigma_0\ | \ \tau(z)<+ \infty {\rm{ \ and\  }} \gamma_z^{'}(\tau(z)) {\rm{ \ is \ transverse \ to \ }}T\Sigma_0 \}$$
    It is clear that $\Omega_r$ is an open set, and it is nonempty because any periodic point of
  $h$ belongs 
to $\Omega_r$ (for such points, $\gamma_z^{'}(\tau(z))$ is actually orthogonal to $T\Sigma_0$).
  Now the map $\varphi: z \mapsto (\tau(z), \gamma_z^{'}(\tau(z)))$ is continuous on $\Omega_r$.  When $z$ is periodic for $h$, 
  $\gamma_z'(\tau(z))=\epsilon(z) \nu(\gamma_z(\tau(z)))$, where $\epsilon(z):=\pm 1$.  
  By density of such periodic points, 
  we get that $\varphi$ maps continuously $\Omega_r$ to $\RR_+ \times \{-1,+1\}$.  
  Observe that $\Omega_r$, as well as $\varphi$ are $h$-invariant.
  Because $h$
  is topologically transitive on $\Omega_r$, this implies that $\varphi$ is actually a constant map  
  $z \mapsto (\tau_r,\epsilon)$. 
  We call $\tau_r$ the {\it first return time} of the normal flow, and $\phi^{\tau_r}: \Omega_r \to \Sigma_0$
   {\it the first return map}.

 \subsection{First geometric properties of the normal flow}  
 \label{sec.geometric.properties}
  
  The previous section shows that the normal flow $\phi^t$ is defined on $(0,\tau_m)$.  We detail here its main geometric properties.
  
 \begin{proposition}
   \label{fact.lorentz}
   \begin{enumerate}
   \item{For each $t \in (0,\tau_m)$, $\phi^t$ is an homothetic transformation from $(\Sigma_0,g_0)$ to $(\Sigma_t,g_t)$. More precisely
    $a(t)>0$ and $(\phi^t)^*g_t=a(t)g_0$.}
   \item The tensor $\tilde{g}$ is a Lorentz metric on $(0,\tau_m) \times \Sigma_0$, and $f: ((0,\tau_m) \times \Sigma_0, \tilde{g}) \to (M,g)$
     is a one-to-one, orientation preserving, isometric immersion.
   \end{enumerate}

  \end{proposition}
  \begin{preuve}
   Equivariance relation (\ref{eq.equivariance}) implies that $h$ acts as an Anosov diffeomorphism on $\Sigma_t$, and $\mathcal{E}^{-}(t)$ 
   (resp. $\mathcal{E}^{+}(t)$) generates the stable (resp. unstable) bundle of $h$ at $\gamma_{z_0}(t)$. It follows that $\mathcal{E}^{\pm}(t)$
    are lightlike, and linearly independant since $D_{z_0}\phi^t$ is one-to-one. 
     This implies $a(t)=g(\mathcal{E}^{-}(t),\mathcal{E}^{+}(t))\not = 0$.  Because $a(0)=1$, we get 
   $a(t)>0$ for $t \in (0,\tau_m)$.  
   Relation (\ref{eq.equivariance}) shows that  $\varphi^t$ maps the stable (resp. unstable) foliation of $h$ on ${\Sigma}_0$ to the stable 
  (resp. unstable) foliation
   of $h$ on ${\Sigma}_t$. Hence the differential $D_z\varphi^t$ maps at each point $z$ of ${\Sigma}_0$ the lightcone of ${T}_z{\Sigma}_0$ 
   to the lightcone of 
   ${T}_{\varphi^t(z)}{\Sigma}_t$, which means that 
    $(\varphi^t)^*g_t=\sigma_tg_0$, for some smooth function $\sigma_t : {\Sigma}_0 \to \RR_+^*$.  Now,  
     because of (\ref{eq.equivariance}),  the function  $\sigma_t$ is $h$-invariant, hence constant since  $h$ admits dense orbits
     on ${\Sigma}_0$.  This constant is given by
     $g(D_{z_0}\varphi^t(E^-(z_0)), D_{z_0}\varphi^t(E^+(z_0)))$, namely 
       ${a}(t)$.

   The second point follows easily.  Indeed, we already noticed that $a(t)>0$ for $t \in (0,\tau_m)$, which ensures that $\tilde{g}$ is Lorentzian on 
   $(0,\tau_m) \times \Sigma_0$. By the first point, $f$ will be isometric if we prove that $T_{\gamma_z(t)}\Sigma_t$ is orthogonal to 
   $\gamma_z'(t)$ for all $t \in (0,\tau_m)$.  Now, observe   that for a linear Lorentz transformation 
      $L=\left( \begin{array}{ccc} \lambda&0&0\\ 0& \frac{1}{\lambda}&0\\ 0&0&1 \end{array} \right)$ with $|\lambda|>1$, the only
       Lorentz plane invariant by $L$ is the one generated by the two first basis vectors, namely the orthogonal to the line of fixed point
       of $L$.  This remark shows that  if  $z \in \Sigma_0$ is a periodic point for $h$, $T_{\gamma_z(t)}\Sigma_t \perp \gamma_z'(t)$ holds.
       By density of periodic points of $h$ on ${\Sigma}_0$, the property actually holds  for all
          $z \in {\Sigma}_0$.  Finally $f$ is orientation preserving because $(E^-,E^+,\nu)$ is positively oriented, and we 
          defined the orientation of $N$ to make $(\tilde{E^-},\tilde{E}^+,\ddt)$ positively oriented.
  \end{preuve}

 Let us also say a few words about the surfaces $\Sigma_t$. We observe that generally, $\Sigma_t$ are  not totally geodesic submanifolds of $M$.  However, they enjoy the weaker condition:
 \begin{fact}
  \label{fait.lightlike.geodesics}
  For every $t \in (0,\tau_m)$, the parametrized lightlike geodesics 
  of $\Sigma_t$ for the metric $g_t$ are parametrized geodesics for the metric $g$.
 \end{fact}

 This can be checked directly by computation for the metric $\tilde{g}$ on $(0,\tau_m) \times \Sigma_0$.  
%
  From Fact \ref{fait.lightlike.geodesics}, we infer 
  the following relation, available for all $z \in \Sigma_0$, $t \in (0,\tau_m)$ and  $s \in \RR$~:
  \begin{equation}
   \label{eq.equivariance2}
   \phi^t(\exp(z,sE^{\pm}(z)))=\exp(\phi^t(z),sD_z\phi^t(E^{\pm}(z)))
  \end{equation}

 \subsection{Completeness of the normal flow}
 Thanks to the previous section, we understand pretty well the behaviour of the normal flow for $t \in (0,\tau_m)$.  Our next step is to show 
 that the flow can be
  extended for $t \geq \tau_m$.
 \subsubsection{Extension of the normal flow  at $t=\tau_m$}
 An important step toward extending the normal flow for $t=\tau_m$ is to show
  that the geometry of the Lorentz manifold $((0,\tau_m) \times \Sigma_0, \tilde{g})$ does not degenerate at the boundary.
  Precisely, we prove :
  
  \begin{lemme}
   \label{lem.nondegeneracy}
   There exists $\epsilon >0$ such that $\tilde{g}$ is a Lorentz metric on the open set $(- \epsilon, \tau_m+\epsilon) \times \Sigma_0 \subset N$.
  \end{lemme}

  \begin{preuve}
   We just have to show that $a(\tau_m) \not =0$.  
   
   Recall the  point $z_0 \in {\Sigma}_0$ we introduced at the  begining of Section \ref{sec.push}. 
   This point is fixed by $h$, and we already 
   observed that $f(t,z_0)$ 
  exists for all $t \in \RR$. 
  Hence, on  ${\mathcal U} \subset {\Sigma}_0$, a small convex neighborhood of $z_0$ (convex relatively to 
  the metric $g_0$),
  we have an extended flow  $\phi^t: {\mathcal U} \to M$  defined for $t \in (0,\tau_m+\delta)$, $\delta>0$.  Saturating 
  ${\mathcal U}$ by the action of $(h^m)_{m \in \ZZ}$, we get a dense open set ${\mathcal V}$  on which $\phi^t$ is 
  defined for $t \in (0,\tau_m+\delta)$.  Observe that  ${\mathcal V}$ contains de stable and unstable manifolds of $h$ at 
  $z_0$, namely $W^{\pm}:=\{ \exp(z_0,sE^{\pm}) \ | \ s \in \RR  \}.$
  
  Along the geodesic $t \mapsto \gamma_{z_0}(t)$, we define two vector fields $(E^-(t),E^+(t))$  by parallel transporting 
       $(E^-(z_0),E^+(z_0))$. For each 
    $t \in (0,\tau_m)$, relation (\ref{eq.equivariance}) yields the existence of two nonzero reals $\lambda_t^{\pm}$ such
     that ${\mathcal E}^{\pm}(t)=D_{z_0}\varphi^t(E^-(z_0))=\lambda_t^{\pm}E^{\pm}(t))$.
     Observe that ${a}(t)=\lambda_t^+ \lambda_t^-$.  Proving  $a(\tau_m) \not =0$ amounts to show that $\lambda_t^{\pm}$ are bounded away from
      $0$ in $(0,\tau_m)$.  
      
      Assume for a contradiction that there exists 
      some
      sequence $(t_k)$ in $(0,\tau_m)$, such that $t_k \to \tau_m$ and  $\lambda_{t_k}^- \to 0$. 
      Relation (\ref{eq.equivariance2})  says that for $s \in \RR$, and 
       $k \in \NN$
      $$ \phi^{t_k}(\exp(z_0,sE^{-}(z_0)))=\exp(\phi^{t_k}(z_0),s\lambda_{t_k}^-E^-(t_k)).$$
      This implies $\phi^{\tau_m}(\exp(z_0,sE^-(z_0)))=\phi^{\tau_m}(z_0)$ for all $s \in \RR$. In particular, because the unstable
       manifold $W^-$ is dense in ${\mathcal V}$, we get that $\phi^{\tau_m}(z)=\phi^{\tau_m}(z_0)$ for every
       $z \in {\mathcal V}$.  Let us choose a $h$-periodic point  $z_1 \in {\mathcal V}$, $z_1 \not = z_0$, of period 
       $q \in \NN^*$ (such a point exists by density of 
       $h$-periodic points in $\Sigma_0$). By what we just said, $\gamma_{z_0}(\tau_m)=\gamma_{z_1}(\tau_m)$. We observe that  
     $\gamma_{z_0}^{\prime}(\tau_m)$ and $\gamma_{z_1}^{\prime}(\tau_m)$ can not be linearly independent, otherwise  
       $D_{\gamma_{z_0}(\tau_m)}h^q$ would fix pointwise a $2$-dimensional space in $T_{\gamma_{z_0}(\tau_m)}M$,  
       implying that at $\gamma_{z_0}(\tau_m)$, 
        $Dh^q$ is trivial or has order $2$. 
      A Lorentz isometry being completely determined by its first jet at a given point, this situation would lead to 
      $h^{2q}=id$, a contradiction. We infer that $\gamma_{z_0}'(\tau_m)=-\gamma_{z_1}'(\tau_m)$, and 
      $z_1=\gamma_{z_0}(2 \tau_m)$.  Applying the same argument to a periodic point $z_2 \in {\mathcal V}$ different from $z_0$ and
      $z_1$, we get a contradiction.
      Interverting the role of $W^-$ and $W^+$, the same argument holds if $\lambda_{t_k}^+ \to 0$ for some sequence $t_k \to \tau_m$, 
      and the lemma follows.

  \end{preuve}

  We have shown that the Lorentz metric $\tilde{g}$ on $(0,\tau_m) \times \Sigma_0$ extends to a Lorentz metric
   on $(-\epsilon, \tau_m+\epsilon) \times \Sigma_0$.  Our next goal is to extend our isometric embedding $f$
   to a map $\overline{f}: [0,\tau_m] \times \Sigma_0$.  This will be done thanks to the following general
   extension result, which is of independent interest.
 \begin{proposition}
  \label{prop.extension}
  Let $(L,\tilde{g})$ be a Lorentz manifold, and $\Omega \subset L$ an open subset such that the closure $\overline{\Omega}$ is a manifold 
  with boundary. Assume 
   that the boundary $\partial \Omega$ is a smooth Lorentz hypersurface of $L$.  
  If $(M,g)$ is another Lorentz manifold having same dimension as $N$, and if $f: (\Omega,\tilde{g}) \to (M,g)$ is a one-to-one  isometric 
  immersion, then $f$ extends to
   a smooth isometric immersion  $\overline{f} : \overline{ \Omega} \to M$. 
\end{proposition}

In the previous proposition, smooth isometric immersion means that $\overline{f}: \overline{\Omega} \to M$ admits a well defined 
differential  
$D_z{\overline{f}}: T_zN \to T_{\overline{f}(z)}M$ for every $z \in \overline{\Omega}$, which is isometric with respect 
to $\tilde{g}$ and $g$, 
and varies smoothly with $z$. 

\begin{preuve}
 The main part of the proof is to show the following:
 
 \begin{lemme}
 \label{lem.extension.locale}
 Each point $x \in \partial \Omega$ admits an open neighborhood $U_x \subset N$ such that
 \begin{enumerate}[i)]
  \item  The sets $U_x \cap \Omega$ and ${\mathcal U}_x:=U_x \cap \partial \Omega$ are connected.
  \item There exists a smooth injective immersion $\tilde{f}_x: U_x \to M$  such that $\tilde{f}_x$ and $f$ coincide on $U_x \cap \Omega$.
 \end{enumerate}

\end{lemme}
 
 \begin{preuve}
  We consider, at $x$, the vector $\nu$ which is normal to $T_x(\partial \Omega)$, and points toward $\Omega$.  We consider
   $\gamma$ a small  geodesic segment starting from $x$ and satisfying $\gamma'(0)=\nu$, as well 
   as a sequence $(x_k)$ of points of $\gamma \cap \Omega$ converging to $x$.  Since $M$ is compact, we may assume 
   that $f(x_k)$  converges to a point $y \in M$.  In small neighborhoods of $x$ and $y$, we choose two orthonormal 
   frame fields, which yield at each points $z,z'$ of those neighborhoods, isometric identifications 
   $i_z: \RR^{1,n-1} \to (T_zL,\tilde{g})$, $i_{z'}:\RR^{1,n-1} \to (T_{z'}M,g)$ (here, $\RR^{1,n-1}$ 
   stands for $n$-dimensional Minkowski space). Obviously, one can choose our orthonormal frame fields such that
   $i_{\gamma(t)}^{-1}(\gamma'(t))$ is a constant vector $\xi \in \mathcal{U}$. Also, there are ${\mathcal U}, {\mathcal V}$ neighborhoods of the origin
    in $\RR^{1,n-1}$ such that $u \mapsto \exp(z,i_z(u)$, $u \in {\mathcal U}$,  and $v \mapsto \exp(z',i_{z'}(v))$, $v \in {\mathcal V}$,
    make sense and are diffeomorphisms on their images.  In the trivialization given by the frame fields, the 
     sequence of differentials 
     $(D_{x_k}f)$  becomes a sequence of matrices $(A_k)$ in $\OO(1,n-1)$.  Since $f$ is an isometry, we have the relation
     \begin{equation}
      \label{eq.exponentielle}
      f(\exp(x_k,u))=\exp(f(x_k),A_k(u))
     \end{equation}
for every $u \in {\mathcal U}$.
   
   We can prove the lemma if we show that  the sequence $(A_k)$ is contained in a compact set of $\OO(1,n-1)$. 
   For if it is
   the case, we may assume $A_k \to A_{\infty}$, and shrinking maybe $\mathcal{U}$, we will have $A_k(\mathcal{U}) \subset \mathcal{V}$ 
   for all $k \in \NN$.  Then, we  choose ${\mathcal C} \subset \RR^{1,n-1}$ an open cone with vertex $0$, 
   containing $-\xi$ and contained in 
   ${\mathcal U}$.  For $k_0$ large enough, $U_x=\exp(x_{k_0}, i_{x_{k_0}}({\mathcal C}))$  contains $x$, and if
   ${\mathcal C}$ is chosen connected and narrow enough around $-\xi$, $U_x \cap \Omega$ and $U_x \cap \partial \Omega$ 
   are connected.  The map $f_x : U_x \to M$ given by $f_x(\exp(x_{k_0},i_{x_{k_0}}(u)))=\exp(y_{k_0},i_{y_{k_0}}(u))$, 
   $u \in {\mathcal C}$ is a one-to-one immersion which coincides with $f$ on $U_x \cap \Omega$.
   
   It remains to explain why the sequence $(A_k)$ must be bounded.  If not, we apply the $KAK$ decomposition of
   $\OO(1,n-1)$ to the sequence $(A_k)$, and after considering a subsequence we can write 
    $A_k$ as a product $M_kD_kN_k$ with $M_k \to M_{\infty}$  (resp. $N_k \to N_{\infty}$) in 
    $\OO(1,n-1)$, and $D_k=\left(  \begin{array}{ccc}  \lambda_k & & \\
     & \ddots & \\
      & & \lambda_k^{-1}\\ \end{array} \right)$, $|\lambda_k| \to \infty$.
      We see that there exists a lightlike hyperplane ${\mathcal H} \subset \RR^{1,n-1}$ (namely the image by
      $N_{\infty}^{-1}$ of $\operatorname{Span}(e_2,\ldots,e_n)$) with the following dynamical 
      property: For every $u \in {\mathcal H}$, there exists $u_k \to u$ such that after extracting a subsequence,  
      $A_k(u_k) \to u_{\infty}$.  Moreover, we see that if $v \not \in {\mathcal H}$, one can find 
      a sequence of reals $s_k \to 0$ such that $A_k(s_kv) \to v_{\infty} \not = 0$.
      
      Because ${\mathcal H}$ is lightlike while $T_x(\partial \Omega)$ has Lorentz signature, one can find a nonzero 
      $u \in {\mathcal H} \cap {\mathcal U}$ such that $i_{x}(u) \not \in T_x(\partial \Omega)$ and 
      $i_{x}(u)$ points  toward $\Omega$.  We choose a sequence $(u_k)$ in $\mathcal{U}$ converging to $u$ such that
       $A_k(u_k)$ tends to $u_{\infty}$  (after extraction).
       We can also pick some  $v \in {\mathcal{U}}\setminus { \mathcal{H}}$ such that $i_x(v)$ points toward $\Omega$.
       Then we can find $(s_k)$ a sequence of reals tending to $0$ such that $A_k(s_k v)$ converges to $v_{\infty} \not =0$, and $(u_k)$.
        Observe that $\exp(x_k,u_k)$ and $\exp(x_k, u_k + s_kv)$ belong to $\Omega$ for $k$ large.   Now,
        $f(\exp(x_k,u_k))$ tends to $f(\exp(x,u))$, and relation (\ref{eq.exponentielle}) shows 
        that $f(\exp(x,u))=\exp(y,u_{\infty})$.  On the other hand, $f(\exp(x_k, u_k+s_k v))$ should also converge to
         $f(\exp(x,u))$, because $s_k \to 0$.  But relation (\ref{eq.exponentielle}) says that this sequence actually 
         converges to
          $\exp(y,u_{\infty}+v_{\infty})$.  Since $v_{\infty} \not =0$, and because we can rescale $u$ and $v$ so that 
          $u_{\infty}$ and $u_{\infty}+v_{\infty}$ belong to ${\mathcal V}$, we have 
          $\exp(y,u_{\infty}+v_{\infty}) \not = \exp(y,u_{\infty})$, and we get  a contradiction.
 \end{preuve}

 Lemma \ref{lem.extension.locale} easily provides a smooth extension of $f$, $\overline{f}: \overline{\Omega} \to M $,
 putting $\overline{f}(x)=f(x)$ if $x \in \Omega$, and 
 $\overline{f}(x):=\tilde{f}_x(x)$ for every $x \in \partial \Omega$.  This extension map $\overline{f}$  is well defined because if $U_x \cap U_{x'} \not = \emptyset$, 
 then $\tilde{f}_x$ and $\tilde{f}_{x'}$ are equal to $f$ on $U_x \cap U_{x'}\cap \Omega$, hence on  $U_x \cap U_{x'}\cap \partial \Omega$.
  Observe that the relation $f^*g=\tilde{g}$ which is available on $U_x \cap \Omega$ must still hold 
  on $U_x \cap \overline{\Omega}$.  This proves that $\overline{f}$  is an isometric immersion.
\end{preuve}

\subsubsection{The normal flow at $\tau_m$ realizes the first return map}
We apply Proposition \ref{prop.extension} choosing for $L$ the Lorentz manifold 
$((- \epsilon, \tau_m + \epsilon) \times \Sigma_0, \tilde{g})$ and for ${\Omega}$ the product 
$(0,\tau_m) \times \Sigma_0$.  We get a smooth, orientation preserving,  extension 
$\overline{f} : ([0,\tau_m] \times \Sigma_0, \tilde{g}) \to (M,g)$ 
which is an isometric immersion coinciding with $f$ on $(0,\tau_m) \times \Sigma_0$.
\begin{proposition}
 The extension ${\overline{f}}$ maps $\{ \tau_m\} \times {\Sigma}_{0}$  diffeomorphically and isometrically onto $\Sigma_0$.  In other words,
 the first return time $\tau_r$ coincides with $\tau_m$,
  and $\phi^{\tau_m}: \Sigma_0 \to \Sigma_0$ realizes the first return map.
\end{proposition}

\begin{preuve}  In the proof, we are going to write $\tilde{\Sigma}_{\tau_m}$ (resp.  $\tilde{\Sigma}_0$) instead of 
$\{\tau_m \} \times \Sigma_0$ (resp. $\{0\} \times \Sigma_0$).
  We recall the notations $\frac{\partial}{\partial t}, \tilde{E}^-,\tilde{E}^+$ from Section \ref{sec.normalflow}. We first show that the restriction of ${\overline{f}}$ to $\tilde{\Sigma}_{\tau_m}$ is one-to-one.  Assume for a 
contradiction that it is not the case. We get two points $\tilde{z}_1=(\tau_m,z_1)$ and $\tilde{z}_2=(\tau_m,z_2)$ on $\tilde{\Sigma}_{\tau_m}$ 
  such that ${\overline{f}}(\tilde{z}_1)={\overline{f}}(\tilde{z}_2)$.  Observe that 
  $D_{\tilde{z}_1}{\overline{f}}(T_{\tilde{z}_1}\tilde{\Sigma}_{\tau_m})=D_{\tilde{z}_2}{\overline{f}}(T_{\tilde{z}_2}\tilde{\Sigma}_{\tau_m})$ because a transverse intersection 
   of those two subspaces would not be compatible with injectivity of ${\overline{f}}$ on $(0,\tau_m) \times \Sigma_0$.
     Looking at orthogonal subspaces, and because ${\overline{f}}$ is an isometric immersion, we get
     $D_{\tilde{z}_1}{\overline{f}}(\frac{\partial}{\partial t})=\pm D_{\tilde{z}_2}{\overline{f}}(\frac{\partial}{\partial t})$.
      Again, $D_{\tilde{z}_1}{\overline{f}}(\frac{\partial}{\partial t})= D_{\tilde{z}_2}{\overline{f}}(\frac{\partial}{\partial t})$ would violate the injectivity of ${\overline{f}}$ on  $(0,\tau_m) \times \Sigma_0$.
        We can reformulate equality $D_{\tilde{z}_1}{\overline{f}}(\frac{\partial}{\partial t})=- D_{\tilde{z}_2}{\overline{f}}(\frac{\partial}{\partial t})$, saying that 
        $\gamma_{z_1}(\tau_m)=\gamma_{z_2}(\tau_m)$, and 
  $\gamma_{z_1}^{'}(\tau_m)=-\gamma_{z_2}^{'}(\tau_m)$. 
  It follows that $z_1$  and $z_2$ are in the first return set $\Omega_r$, $z_2=\theta_r(z_1)$, 
  and the first return time $\tau_r$ 
  equals
   $2 \tau_m$.  As a consequence, we get for every $z \in \Omega_r$, the identity:
   $$ {\overline{f}}(\tau_m,z)={\overline{f}}(\tau_m,\theta_r(z)).$$
   Defining $\tilde{\theta}_r(t,z):=(t,\theta_r(z))$, this identity yields:
   $$ D_{(\tau_m,z)}{\overline{f}}(\tilde{E}^{-})=D_{\tilde{\theta}_r(\tau_m,z))}{\overline{f}}(D_{(\tau_m,z)\tilde{\theta}_r}(\tilde{E}^-)).$$
   Because $\theta_r$ commutes with $h$, there exists $\alpha(z) \not = 0$ such that 
   $D_{(\tau_m,z)\tilde{\theta}_r}(\tilde{E}^-)=\alpha(z)\tilde{E}^-$.  We thus obtain:
   $$  D_{(\tau_m,z)}{\overline{f}}(\tilde{E}^{-})= \alpha(z) D_{\tilde{\theta}_r(\tau_m,z))}{\overline{f}}(\tilde{E}^-).$$
   For the same reasons, there exists $\beta(z) \not =0$ such that
   $$  D_{(\tau_m,z)}{\overline{f}}(\tilde{E}^{+})= \beta(z) D_{\tilde{\theta}_r(\tau_m,z))}{\overline{f}}(\tilde{E}^+).$$
   
   Going back to $z=z_0$, $\theta_r(z)=z_1$, we see that $(D_{(\tau_m,z_1)}{\overline{f}})^{-1} \circ D_{(\tau_m,z_0)}{\overline{f}}$ is a 
   linear isometry, preserving the orientation, and sending the direct frame $(\tilde{E}^-,\tilde{E}^+,\ddt)$ at 
   $(\tau_m,z_0)$ to the frame $(\alpha(z_0) \tilde{E}^-, \beta(z_0) \tilde{E}^+,-\ddt)$ at $(\tau_m,z_1)$.  
   The isometric condition yields $\alpha(z_0) \beta(z_0)=1$ and the orientation-preserving condition yields
    $- \alpha(z_0) \beta(z_0)$.  This provides the desired contradiction.

  Once we know that ${\overline{f}}$ is one-to-one in restriction to $\tilde{\Sigma}_{\tau_m}$, we get that
  ${\overline{f}}(\tilde{\Sigma}_{\tau_m})$ is a Lorentz surface of $M$, to which we can again apply the normal flow.  
  This results into 
  an extension  of ${\overline{f}}$ to a smooth immersion defined on a domain $(0,\tau_m+\epsilon) \times \Sigma_0$. 
  If ${\overline{f}}(\tilde{\Sigma}_{\tau_m})$ does not meet $\Sigma_0$, it is easily checked that for $\epsilon>0$ 
  small enough, ${\overline{f}}$ is
   one-to-one on $(0,\tau_m+\epsilon) \times \Sigma_0$, contradicting the definition of $\tau_m$.
   
  We infer that there exist $z_1$ and $z_2$ in $\Sigma_0$ such that ${\overline{f}}(\tau_m,z_1)={\overline{f}}(0,z_2)$.  
  Observe that 
  $D_{(\tau_m,z_1)}{\overline{f}}(T\tilde{\Sigma}_{\tau_m})$ and $D_{(0,z_2)}{\overline{f}}(T \tilde{\Sigma}_0)$ can not intersect  transversely, otherwise
   ${\overline{f}}$ would not be one-to-one on $(0,\tau_m) \times \Sigma_0$.  We can recast this property saying that
   $\gamma_{z_1}(\tau_m)=z_2$, and $\gamma_{z_1}'(\tau_m)$ is orthogonal to $T_{z_2} \Sigma_0$.  In other words, 
   $z_1$ belongs to the first return set $\Omega_r$, and the return time is $\tau_r=\tau_m$.  It means in particular that for all 
    $z \in \Omega_r$, ${\overline{f}}(\tau_m,z) \in \Sigma_0$.  By density of $\Omega_r$ in $\Sigma_0$, we finally get that
     ${\overline{f}}$ maps $\tilde{\Sigma}_{\tau_m}$ isometrically and diffeomorphically onto $\Sigma_0$.  
  

\end{preuve}

\subsection{End of proof of Theorem \ref{theo.hyperbolic}}
\label{sec.conclusionhyperbolic}
  Let us just recollect what we did so far.  First, showing that the normal flow $\phi^t$ is defined on 
  $(-\epsilon,\tau_m+\epsilon)$ with 
  $\phi^{\tau_m}(\Sigma_0)=\Sigma_0$ immediately implies that $\phi^t$ is defined for every   $t \in \RR$.  Equivalently, the map $f$ is defined
   on all of $N=\RR \times \Sigma_0$.  
   
   Next, Proposition \ref{fact.lorentz}  implies that $(\phi^{\tau_m})^*g_0=a(\tau_m)g_0$, with $a(\tau_m)>0$.  Because the global
   Lorentz volume of $\Sigma_0$ must be preserved, we get $a(\tau_m)=1$.  The transformation $\phi^{\tau_m}$ is a Lorentz isometry of 
   $(\Sigma_0,g_0)$ commuting with $h$: It must be either $\pm Id$ or a linear hyperbolic transformation. 
    The possibility $\phi^{\tau_m}=- Id$ is ruled out by the assumption that $(M,g)$ is time-orientable.
    
    In the following, we denote by $A$ 
   the transformation $\phi^{\tau_m}$.  We just showed that $t \mapsto a(t)$ is $\tau_m$-periodic, and thanks to Propositions \ref{fact.lorentz}
    and \ref{prop.extension}, we get that $f: (N,\tilde{g}) \to (M,g)$ is an isometric immersion.  Let us call $\varphi : N \to N$ the transformation
     $\varphi(t,x)=(t+1,A^{-1}x)$.  Then $\varphi$ acts isometrically for $\tilde{g}$, and $f \circ \varphi = f$.  Calling  $\Gamma$ the cyclic
      group generated by $\varphi$, we finally see that $f$ induces an isometry between $\Gamma \backslash N$ (endowed with the metric induced by ${\tilde g}$) and 
      $(M,g)$.  This shows the topological part of Theorem \ref{theo.hyperbolic}.
      
      Since $\Sigma_0$ is a flat torus, the universal cover $(\tilde{M},\tilde{g})$ is isometric to $\RR^3$ endowed with
       the metric $dt^2+2a(t)dudv$.  Affine transformations preserving the planes $t=t_0$ and acting by Lorentz isometries
        on the Minkowski $(u,v)$-plane, provide an isometric action of $\Sol$ on $(\tilde{M},\tilde{g})$.
         This shows points $2)$ and $3)$ of Theorem \ref{theo.hyperbolic}.


\section{The local geometry of manifolds with no hyperbolic component}
\label{sec.local.geometry}

We keep going in our study of closed $3$-dimensional Lorentz manifolds $(M,g)$, such that $\Iso(M,g)$ is not compact.
Thanks to sections \ref{sec.homogeneous} and \ref{sec.hyperbolic}, we can prove Theorem \ref{thm.topological}   when
 all the components of the integrability locus $\mint$ are locally homogeneous, or when there exists at least 
 one hyperbolic component.  Looking at the posibilities
  for the different components listed in Section \ref{sec.components}, it only remains to investigate the case 
  where all the components are either 
  parabolic or of constant curvature, and there is at least one non locally homogeneous component.  This section is devoted 
  to a careful geometric study of such manifolds, and our aim is to prove the

  \begin{theoreme}
\label{thm.conf.flat}
Let $(M,g)$ be a $3$-dimensional Lorentz manifold. If all the components of $\mint$ are   of constant curvature or parabolic, and if  $(M,g)$ is not locally homogeneous, 
then
 $(M,g)$ is conformally flat.
\end{theoreme}

Recall that $(M,g)$ is said to be conformally flat if each sufficiently small open neighborhood of $M$ is conformally 
diffeomorphic to an open subset of Minkowski space.  Hence, Theorem \ref{thm.conf.flat} tells us that {\it 
at the conformal level}, 
our structure $(M,g)$ is locally homogeneous.  This local information will be decisive to recover the global properties of
 $(M,g)$, both topologically and geometrically, a  task that  will be carried over in Section \ref{sec.global}.


\subsection{More on the geometry of parabolic components}
Parabolic components split into two categories, the locally homogeneous ones for wich the Killing algebra is $4$-dimensional,
 and the others for which it is $3$-dimensional.
 
 \subsubsection{Locally homogeneous parabolic components}
 \label{sec.lochomoparabolic}
 The study of locally homogeneous parabolic components was made in \cite{frances.open}, 
 and can be summarized as follows:

\begin{proposition}\cite[Proposition 4.3]{frances.open}
  \label{prop.lochomoparabolique}
  Let $\calm$ be a component of the integrability locus $\mint$ of a $3$-dimensional Lorentz manifold $(M,g)$, which is 
   locally homogeneous and parabolic.  
   Then:
  \begin{enumerate}
   \item{If the scalar curvature of $\calm$ is $0$, the Lie algebra $\kiloc(\calm)$ is isomorphic to a 
   semi-direct product
   $\RR \ltimes \heis,$ where $\heis$ stands for the $3$-dimensional Heisenberg Lie algebra.}
   \item{If the scalar curvature is nonzero on $\calm$, then  $\kiloc(\calm)$ is isomorphic to $\sld \oplus \RR$.}
  \end{enumerate}
\end{proposition}

Actually, the statement of \cite{frances.open} is slightly more precise since it describes which semi-direct products 
 $\RR \ltimes \heis$ can occur.  However, we won't need this extra information here.  For the sequel, it will be important to notice
  that in the first case of  Proposition \ref{prop.lochomoparabolique}, the Lie subalgebra $\heis$ contains the 
  isotropy algebra at each points, hence 
  acts with $2$-dimensional pseudo-orbits  (this follows from the computations 
  done in \cite[Section 4.3.3]{frances.open}).
 
 \subsubsection{Parabolic components which are not locally homogeneous}
We  investigate now the geometry of parabolic components which are not locally homogeneous. 
 \begin{proposition}
 \label{prop.heis}
 Let $\calm \subset \mint$ be a parabolic component which is not locally homogeneous.
 \begin{enumerate}
 \item{The Lie algebra $\kiloc(\calm)$ is isomorphic to the $3$-dimensional
  Heisenberg algebra $\heis$. }
 \item{The $\kiloc$-orbits on $\calm$ are totally geodesic, lightlike surfaces.}
 \item{The scalar curvature $\sigma$ vanishes on $\calm$.}
 \end{enumerate}
 \end{proposition}

 The proof of Proposition \ref{prop.heis}  involves quite a bit of computations, that we defer to Annex A, at the end 
 of the text.
 Let us mention an important corollary which will be crucial later on.
 
 \begin{corollaire}
  \label{coro.islocorbites}
  Let $\calm$ be a parabolic component which is not locally homogeneous.  Let $x \in \calm$.  Then the $Is^{loc}$-orbit 
  of $x$ is a submanifold of $\mint$ which is closed in $\mint$.
 \end{corollaire}
\begin{preuve}
 We already know (Theorem \ref{thm.integrabilite}) that the $Is^{loc}$-orbit 
  of $x$ is a submanifold $\Sigma$ of $\mint$.  We have to show that $\Sigma$ is closed in $\mint$.  We thus consider 
  a sequence  $(x_k)$ of $\Sigma$ converging to a point $x_{\infty} \in \mint$.  Let $\hx$ be a lift of $x$ in $\hm$.
   We recall the generalized curvature map $\dk: \hm \to \mathcal{W}$  (see Section \ref{sec.curvaturemap}).  Let us call
    $w=\dk(\hx)$, and $\mathcal{O}.w$ the orbit of $w$ under the action of $\OO(1,2)$ on $\mathcal{W}$.  Since $\calm$ is a 
    parabolic component, $\mathcal{O}.w$ is $2$-dimensional and the isotropy at $w$ is a $1$-parameter unipotent subgroup 
     of $\OO(1,2)$.
    Let $(\hx_k)$ be a sequence of $\hm$ lifting $(x_k)$, such that $\hx_k \to \hx_{\infty}$.  Since $x_k \in \Sigma$ for all $k$,
    we have $\dk(\hx_k) \in \mathcal{O}.w$ for all $k$, and in particular $\dk(\hx_{\infty})$ belongs to the closure 
    $\overline{\mathcal{O}.w}$. The action of $\OO(1,2)$ on $\mathcal{W}$ is algebraic, hence orbits of 
    $\overline{\mathcal{O}.w} \setminus \O.w$ have dimension $<2$.  Since there are no $1$-dimensional orbits in 
    finite dimensional representations of $\OO(1,2)$, we conclude that if $\dk(\hx_{\infty})$ belongs to 
    $\overline{\mathcal{O}.w} \setminus \O.w$, then the stabilizer of $\dk(\hx_{\infty})$ is $3$-dimensional. Since 
    $\hx_{\infty} \in \hmint$, this means that the isotropy algebra $\mathfrak{Is}(x_{\infty})$ is $3$-dimensional, 
    isomorphic to  $\oo(1,2)$
    (Fact \ref{rk.isotropie}).  Let $\calm'$ the component containing $x_{\infty}$.  The points $x_k$ belong to $\calm'$
     for $k$ large enough, what shows $\kiloc(\calm') \simeq \kiloc(\calm)\simeq \heis(3)$.  This contradicts 
     $\mathfrak{Is}(x_{\infty}) \simeq \oo(1,2)$.
     
     We infer that $\dk(\hx_{\infty}) \in \O.w$.  Hence, replacing the sequence $\hx_k$  by $\hx_k.p_k$  for a 
     bounded sequence $(p_k)$ of $P$, we may assume that $\dk(\hx_k)=w$ for all $k$.  By the discussion following Theorem 
     \ref{thm.integrabilite}, $\dk^{-1}(w) \cap \hmint$  is a submanifold, the connected component of which are 
     $\kiloc$-orbits.  We conclude that for $k$ large enough, $\hx_{\infty}$ and $\hx_k$ are in the same $\kiloc$-orbit.
      The same is thus true for $x_{\infty}$ and $x_k$, and the corollary is proved.

\end{preuve}

\subsection{Conformal flatness}
\label{sec.conformal.flatness}
Under the standing assumptions stated at the begining of Section \ref{sec.local.geometry}, 
the only non locally homogeneous components in $M$ are parabolic.
  It follows from Proposition \ref{prop.heis} that the scalar curvature of $g$ is constant on each component, 
  and equal to zero on 
  the non locally homogeneous ones.  As a consequence,  the scalar curvature  vanishes identically on $M$, which implies that
    components of constant 
  sectional curvature are actually flat, hence conformally flat.  It thus remains to show that all parabolic components 
  (locally homogeneous or not) are conformally flat. Observe that conformal flatness is given by a tensorial condition, 
  namely the vanishing of the Cotton-York tensor in dimension $3$, 
  so that $(M,g)$ will be conformally flat as soon as a dense open subset of $M$ is.  Observe also 
  that the vanishing of the scalar
  curvature  says that 
   locally homogeneous parabolic components are exactly those described 
  by the first point
  of  Proposition \ref{prop.lochomoparabolique}. This fact, together with Proposition \ref{prop.heis} and 
  the remark after Proposition \ref{prop.lochomoparabolique} reduces the proof of Theorem \ref{thm.conf.flat} to
  the following general observation:

\begin{proposition}
 \label{prop.heisenberg.flat}
 Let $(N,h)$ be a $3$-dimensional Lorentz manifold.  
 Assume that there exists on $N$ a Lie algebra $\lien$ of Killing fields which is isomorphic to 
 $\heis(3)$, and whose   pseudo-orbits have dimension $\leq 2$.  Then all pseudo-orbits are $2$-dimensional and lightlike, and
  $(N,h)$ is conformally flat.
\end{proposition}

\begin{preuve}
Our hypothesis that all pseudo-orbits have dimension $\leq 2$ implies that the 
isotropy algebra at each point $x \in N$ is a nontrivial subalgebra of $\lien$.  This isotropy algebra is thus isomorphic to $\RR$, $\RR^2$
 or $\heis(3)$.  Because there is no subalgebra of ${\oo}(1,2)$ isomorphic to $\RR^2$ or ${\heis}(3)$, 
 the isotropy algebra of $\lien$ is $1$-dimensional at each point, and all pseudo-orbits of $\lien$ have dimension $2$. 
Let us consider  $X,Y,Z$ three Killing fields generating $\lien$, and satisfying  the relations $[X,Y]=-Z$  and 
$[X,Z]=[Y,Z]=0$.  We are going 
 to look at the subalgebra ${\mathfrak a}$ spanned by $Y$ and $Z$. Because no subalgebra of ${\oo}(1,2)$ is isomorphic to $\RR^2$, pseudo-orbits 
 of ${\mathfrak a}$ have dimension $1$ or $2$. We claim that the open subset $\Omega$ where the pseudo-orbits of $\liea$ 
 are $2$-dimensional  is dense in $N$.   To see this, let us consider $\Delta$  a $1$-dimensional  pseudo-orbit of 
 $\liea$, and 
  let $x \in \Delta$.  The isotropy, in  ${\mathfrak a}$,  of the point $x$ is spanned by an element  $U=aY+bZ$. 
  There is another 
  vector vector field $V=cY+dZ$ such that $v:=V(x) \not =0$. 
    Since $U$ and $V$ commute, $U$ actually vanishes at each
  point of $\Delta$.  Let $t \mapsto \gamma(t)$  be a geodesic for the metric $h$, satisfying $\gamma(0)=x$,  
  $h(\gamma'(0),v) \not =0$, and $\gamma'(0) \not \in \RR.v$.  
  Clairault's equation ensures that for $t>0$ small enough, $h(\gamma'(t),U(\gamma(t)))=0$ and $h(\gamma'(t),V(\gamma(t))) \not =0$.
   Observe that $U(\gamma(t)) \not = 0$, because locally, the zero set of a nontrivial Killing field on a $3$-dimensional
   manifold is a submanifold of dimension $\leq 1$.  We thus get that $\gamma(t) \in \Omega$ for $t>0$ small, ensuring 
   the density of $\Omega$.
   
   This density property shows that we will be done if we show that $\Omega$ is conformally flat. To this aim, we consider a point $x_0 \in \Omega$. 
    Since, $[Y,Z]=0$ and
 $Y,Z$ span a $2$-dimensional space at each point of $\Omega$, there exist local coordinates $(x_1,x_2,x_3)$ around $(0,0,0)$
  such that $Z=\delx$ and $Y=\dely$.
  Because the orbits of $\lien$ are $2$-dimensional, $X$ is of the form $\lambda \delx + \mu \dely$ for some functions $\lambda$ and $\mu$.
    The bracket relations $[X,Z]=0$ and $[X,Y]=-Z$ lead to 
    $0=\frac{\partial \lambda}{\partial x_1}=\frac{\partial \mu}{\partial x_1}=\frac{\partial \mu}{\partial x_2}$ and
     $\frac{\partial \lambda}{\partial x_2}=1$. Hence we can write
     $$ X=(x_2+a(x_3))\delx + b(x_3)\dely.$$
     Observe that replacing $X$ by $X-a(0)Z-b(0)Y$ won't affect the bracket relations between $X$, $Y$ and $Z$, so that we
      will assume in the following that $a(0)=b(0)=0$.
      
   Let us consider a point $p=(p_1,p_2,p_3)$. The vector field $U=X-(p_2+a(p_3))Z-b(p_3)Y$ is nonzero and vanishes at $p$.
   We compute that at $p$: 
    $$[U,\delx]=0, \ [U,\dely]=- \delx,\  [U,\delz]=-a'(p_3)\delx-b'(p_3)\dely.$$
    Since $U$ belongs to $\lien$, hence 
    is Killing for 
    the metric $h$, we infer that the matrix $A=\left(  \begin{array}{ccc}
                                                       0&-1&a'(p_3)\\
                                                       0&0&b'(p_3)\\
                                                       0&0&0\\
                                                      \end{array}
 \right) $ must be antisymmetric for the Lorentz scalar product $h_p$.  It is readily checked that a rank $1$ nilpotent 
 matrix never has this property (basically because $\exp(tA)$ would be a nontrivial $1$-parameter group in (a conjugate of) $\OO(1,2)$
  fixing pointwise a $2$-plane, which is impossible).  We thus infer that the derivative $b'$ is nowhere $0$.  Because 
  $b(0)=0$, this implies that $\frac{b(x_3)}{x_3}$ (or $b'(0)$ if $x_3=0$) is nowhere $0$. 
  The transformation
   $$ \varphi: (x_1,x_2,x_3) \mapsto (\frac{b(x_3)}{x_3}x_1,\frac{b(x_3)}{x_3}x_2 -a(x_3),x_3).$$
    thus yields  a local diffeomorphism fixing the origin.  Applying $\varphi^*$ to $X,Y,Z$, we get
   $$ X=x_2 \delx+x_3 \dely, \ Y=\frac{x_3}{b(x_3)}\dely, \ {\rm and} \ Z=\frac{x_3}{b(x_3)}\delx.$$
   
   Let again  $p=(p_1,p_2,p_3)$ be a point in our coordinate chart.  The vector field $U=X-p_2 \delx-p_3\dely$ vanishes at 
    $p$, and is a Killing field for $h$ because $U=X-p_2\frac{b(p_3)}{p_3}Z-b(p_3)Y$.  A straigthforward computation yields
     $$ [U,\delx]=0, \ [U,\dely]=-\delx, \ {\rm and}\  [U,\delz]=-\dely.$$
     It follows that the matrix $\left(  \begin{array}{ccc}
                                                       0&1&0\\
                                                       0&0&1\\
                                                       0&0&0\\
                                                      \end{array}
 \right) $ must be antisymmetric with respect to $h_p$.  This allows to see that the matrix of $h_p$ in the frame 
 $(\delx,\dely,\delz)$ is of the form 
 $\left(  \begin{array}{ccc}
                                                       0&0&-\beta(p)\\
                                                       0&\beta(p)&0\\
                                                       -\beta(p)&0&\gamma(p)\\
                                                      \end{array}
 \right) $, with $\beta(p) >0$.  Now, $Z$ and $Y$ being Killing fields for $h$, we see that $\beta$ and $\gamma$ only depend on the variable 
 $x_3$, and we conclude that the metric $h$ writes as :
 $$ -2 \beta(x_3)dx_1dx_3+\beta(x_3)dx_2^2+\gamma(x_3)dx_3^2.$$
  Now, if $x_3 \mapsto \zeta(x_3)$ is a primitive of $\frac{-\gamma(x_3)}{2 \beta(x_3)}$, a change of coordinates
  $$ (x_1,x_2,x_3) \mapsto (x_1+\zeta(x_3),x_2,x_3)$$ shows that  $h$ is locally isomorphic to
  $-2\beta(x_3)dx_1dx_3+\beta(x_3)dx_2^2$, hence is conformally flat.
     
\end{preuve}

\section{Geometry on Einstein's universe}
\label{sec.geomeinstein}

Lorentz conformally flat structures in dimension $n=3$ are examples of $(G,X)$-structures in the sense of Thurston.
 In particular, there is a universal space among those structures, called Einstein's universe $\ein$, such that if 
 $(M,g)$ is Lorentz and conformally flat, there exists  a conformal immersion $\delta: \tilde{M} \to \ein$, which is 
 equivariant under a representation of $\pi_1(M)$ into $\operatorname{Conf}(\ein)$  (see Section \ref{sec.developing} below).  The 
  proof of  Theorem \ref{thm.topological} for manifolds $(M,g)$ satisfying hypotheses of Theorem
  \ref{thm.conf.flat} and with a noncompact isometry group,  will rely in a crucial way on the study of this developing
  map $\delta$. 
   This study will be carried over in the next section \ref{sec.global}, and it will require  a deeper knowledge of 
   the geometry 
   of $\ein$.  That's why we dedicate the present section to studying $\ein$ in more details.   
 The reader eager to learn more about the geometry of $\ein$ is refered to \cite{charlesthese} or \cite{primer}.
\subsection{Basics on  Einstein's universe}
\label{sec.description.einstein}
Einstein's universe is the Lorentz analogue of the Riemannian conformal sphere.  We recall its construction, sticking 
 to dimension $3$, which is the relevant one for our purpose. 

Let $\RR^{2,3}$ be the space $\RR^{5}$ endowed with the quadratic form
$$ Q^{2,3}(x_0, \ldots, x_{4}) = 2x_0x_{4} + 2x_{1}x_{3} + x_2^2$$
 We consider the null cone  
$$\mathcal{N}^{2,3} = \{  x \in \RR^{2,3} \ | \  Q^{2,3}(x)=0 \}$$  
and denote by $\mathcal{\widehat{N}}^{2,3}$ the cone  $\mathcal{N}^{2,3}$
with the origin removed.
The projectivization ${\BP}(\widehat{\mathcal{N}}^{2,3})$ is a smooth submanifold of $\R\BP^{4}$, 
and inherits from the pseudo-Riemannian structure of  $\RR^{2,3}$ a Lorentz
 conformal class (more details can be found in \cite{charlesthese}, \cite{primer}).  
 We call  the $3$-dimensional {\it Einstein universe}, denoted $\ein$ 
 this compact manifold  ${\BP}(\widehat{\mathcal{N}}^{2,3})$ with this conformal structure. 
One can check that a $2$-fold cover of $\ein$ is conformally diffeomorphic to 
 the product $({\bf S}^1 \times {\bf S}^{2},-g_{{\bf S}^1} \oplus g_{{\bf S}^2})$. 

The orthogonal group of $Q^{2,3}$, isomorphic to $\mbox{O}(2,3)$,  
acts naturally on the $4$-dimensional projective space, preserving  $\ein$ and its conformal structure. It turns out (see Theorem 
\ref{th.liouville} below) that $\mbox{PO}(2,3)$
 is the full conformal group of $\ein$.  Observe that $\ein$ is homogeneous under the action of $\mbox{PO}(2,3)$.

\subsubsection{Photons and  lightcones}

It is a remarkable fact of Lorentz geometry 
that all the metrics of a given conformal class have 
the same lightlike geodesics (as sets but not as parametrized curves). 
In the case of Einstein's universe, the lightlike geodesics 
are the projections on $\ein$ of totally isotropic 2-planes $P \subset \RR^{2,3}$ (namely planes $P$ on which   
$Q^{2,3}$ vanishes identically). We will rather use the term {\it photon} for the lightlike geodesics of Einstein's universe.
Observe that all photons  of $\ein$ are simple closed curves, and are endowed with a natural class of 
projective parametrizations.

Given a point $p$ in $\ein$, the {\it lightcone with vertex $p$}, denoted by $C(p)$, 
is the union of all photons containing $p$. 
If $p \in \ein$ is the projection of $u \in {\mathcal N}^{2,3}$,
the lightcone $C(p)$ is just $\BP({ u}^{\bot} \cap \mathcal{N}^{2,3})$.  
The lightcone $C(p)$ is singular (from the differentiable viewpoint) at its vertex $p$, 
and 
$C(p) \backslash \{p \}$ is 
topologically a cylinder.  The entire cone $C(p)$ has the topology of  a $2$-torus pinched at $p$.

\subsubsection{Stereographic projection}
\label{sec.stereographic}
There is for $\ein$ a  generalized notion of 
stereographic projection, which shows 
that $\ein$ is a conformal compactification of the Minkowski space.

Let us call $\RR^{1,2}$ the space $\RR^3$ endowed with the quadratic form 
 $Q^{1,2}(x,x)=2x_1x_3+x_2^2$.
Consider $\varphi : \RR^{1,2} \to \ein$ given in projective coordinates of $\BP(\RR^{2,3})$ by
\begin{equation}
\label{eq.stereographic}
\varphi: x=(x_1,x_2,x_3) \mapsto [-\frac{1}{2}Q^{1,2}(x,x): x_1 : x_2 : x_3 : 1]
\end{equation}
 Then $\varphi$ is a conformal embedding of $\RR^{1,2}$ into $\ein$,
called the inverse {\it stereographic projection } with respect to $p_0:=[e_0]$. 
The image $\varphi(\RR^{1,2})$ is a dense open set of $\ein$ 
with boundary the lightcone $C(p_0)$.  Observe that this proves the fact (rather hard to visualize): The complement of a
lightcone $C(p)$ in $\ein$ is connected.
 
\subsubsection{Developing conformally flat structures into Einstein's universe}
\label{sec.developing}
It is a standard fact that Einstein's universe satisfies an analogue of the classical Liouville's theorem on the sphere.  
Namely:

\begin{theoreme}[Liouville's theorem for $\ein$]
 \label{th.liouville}
 Let $U \subset \ein$ be a connected nonempty open set.  Let $f: U \to \ein$ be a conformal immersion.  Then $f$ is 
 the restriction to $U$ of an unique element of $\mbox{PO}(2,3)$.
\end{theoreme}

The existence of the stereographic projection (\ref{eq.stereographic}), and the transitivity of the action of 
$\mbox{PO}(2,3)$ on $\ein$ shows that $\ein$ is conformally flat.  Liouville's theorem \ref{th.liouville} shows that any
  $3$-dimensional, conformally flat Lorentz structure $(M,g)$ is actually a $(\mbox{PO}(2,3),\ein)$-structure, in the sense of Thurston.
 
 As a consequence, for every conformally flat Lorentz structure $(M,g)$, there exists a conformal immersion
 $$ \delta: (\tm,\tilde{g}) \to \ein$$
 called {\it the developing map} of the structure. Here, $\tm$ is the universal cover of the manifold $M$, and $\tilde{g}$ is the
  lifted metric. 
 This developing map comes with a {\it holonomy morphism}  $\rho: {\mbox{Conf}(\tm,\tilde{g})} \to \mbox{PO}(2,3)$ satisfying
  the equivariance relation:
  
  \begin{equation}
  \label{eq.equiholonomy}
   \delta \circ h = \rho(h) \circ \delta
  \end{equation}
available for every $h \in {\mbox{Conf}(\tm,\tilde{g})}$.


\subsection{More geometry on $\ein$}
\label{sec.more.geometry}

\subsubsection{The foliation $\calfd$}
\label{sec.foliation-einstein}
We refer here to the notations introduced in Section \ref{sec.geomeinstein}.
 Let $P$ be the plane in $\RR^{2,3}$ spanned by the vectors $e_0$ and $e_1$.  The form $Q^{2,3}$ vanishes identically on 
 $P$, hence the projection of $P$ on $\ein$ defines a photon that we will denote by $\Delta$. The 
 open subset obtained by removing $\Delta$ to  $\ein$ will be  called $\Omega_{\D}$.
 
  Given a point $p \in \D$, we consider the lightcone $C(p)$ with vertex $p$.
Since $\D$ is a photon, we have $\D \subset C(p)$. 
Now, the intersection of $C(p)$ with $\Omega_{\D}$, namely $C(p) \setminus \Delta$ is a lightlike hypersurface of 
$\Omega_{\D}$, diffeomorphic to a plane. We call it $F_{\D}(p)$. We now make the observation that in $\ein$, there is no 
nontrivial lightlike triangle, namely if two  photons $\D_1$ and $\D_2$ intersect $\Delta$ transversely 
at two distinct points,
 then $\D_1 \cap \D_2 = \emptyset$.  This is the geometric counterpart of the following algebraic fact: In $\RR^{2,3}$, 
 there are no $3$-dimensional spaces on which 
$Q^{2,3}$ vanish identically.  It follows that if $p \not= p'$ are points of $\D$, $C(p) \cap C(p')=\Delta$, or 
in other words
$F_{\Delta}(p) \cap F_{\Delta}(p^{\prime}) = \emptyset$.  This shows that $\{ F_{\Delta}(p) \}_{p \in \D}$ are the leaves of 
a codimension $1$ lightlike foliation of 
 $\Omega_{\D}$, that we will call $\calfd$. The space of leaves of $\calfd$ is naturally identified with $\Delta$.
 For $x \in \Omega_{\Delta}$, we will adopt the notation $F_{\Delta}(x)$ for 
  the leaf of $\calfd$ containing $x$.
 
 \subsubsection{Symetries of the foliation $\calfd$}
  \label{sec.symetries}
  Let us call $G_{\Delta}$ the stabilizer of $\D$ in $\PO$. Obviously, $G_{\Delta}$  preserves $\Omega_{\Delta}$ and the foliation 
    $\calfd$.  
    
    It is readily checked that this group 
  is a semi-direct product 
  $$G_{\Delta} \simeq \operatorname{PGL}(2,\RR) \ltimes N, $$ where the group $N$ is isomorphic to the $3$-dimensional 
  Heisenberg group $\Heis(3)$, and given in $\PO$ by the matrices:
  
  \begin{equation}
   \label{eq.formule.heisenberg}
  N(x,y,z):=\left( \begin{array}{ccccc}
     1&0&-x&\, \, {-(z+xy)}&\, -\frac{x^2}{2}\\
     0&1&-y&-\frac{y^2}{2}&z\\
     0&0&1&y&x\\
     0&0&0&1&0\\
     0&0&0&0&1\\
     \end{array} \right) \  \ x,y,z \in \RR.
\end{equation}

The factor $\operatorname{PGL}(2,\RR)$  is the subgroup of $\PO$ corresponding to matrices:
\begin{equation}
  R_A:= \label{eq.formule.pgl2}
   \left( \begin{array}{ccc}
           A &0&0\\
           0&1&0\\
           0&0&\frac{A}{\operatorname{det}(A)}\\
          \end{array} \right) \ \ A\in \operatorname{PGL}(2,\RR).
\end{equation}  

 Observe that the action of $G_{\Delta}$ on the space of
leaves of $\calfd$ corresponds to the projective action of the factor $\operatorname{PGL(2,\RR)}$ on $\Delta$.
The  subgroup $S_{\Delta} \subset G_{\Delta}$   which preserves individually all the leaves of $\calfd$ is a semi-direct product 
$$ S_{\Delta} \simeq \RR_+^* \ltimes N,$$
 where the factor $\RR_+^*$ corresponds to matrices:
 
\begin{equation}
  R_{\lambda}:=  \label{eq.facteurR}
  \left( \begin{array}{ccccc}
     \lambda&0& & & \\
     0& \lambda & & & \\
      & & 1 & & \\
      & & & \frac{1}{\lambda}&0\\
      & & & 0& \frac{1}{\lambda}\\
     \end{array} \right), \ \lambda \in \RR_+^*.
\end{equation}

Let us end this algebraic parenthesis by giving more details about the action of the group $N$. Obviously, $N$ fixes
 the point $p_0=[e_0] \in \ein$, hence if we perform a stereographic projection given by formula (\ref{eq.stereographic}),
  the group $N$ becomes a subgroup of conformal transformations of $\RR^{1,2}$.  These transformations are affine, given by
  \begin{equation}
   \label{eq.formuleaffine}
  N(x,y,z)=  \left(  \begin{array}{ccc} 
  1&-y&-\frac{y^2}{2}\\ 0&1&y\\0&0&1\\ \end{array}\right) + \left(  \begin{array}{c} z\\ x\\ 0 \\ \end{array}\right)
  \end{equation}

Inside the group $N$, there is a $2$-dimensional subgroup of translations, denoted $T$, comprising all transformations 
of the form
$$ T(x,z):=Id + \left(  \begin{array}{c} z\\ x\\ 0 \\ \end{array}\right), \ x,z \in \RR.$$ 
In $\PO$, such transformations take the matricial form:  
$$T(x,z)=\left( \begin{array}{ccccc}
     1&0&-x&\,  {-z}& -\frac{x^2}{2}\\
     0&1&0&0&z\\
     0&0&1&0&x\\
     0&0&0&1&0\\
     0&0&0&0&1\\
     \end{array} \right).$$
     
     From this matricial representation, it is straigthforward to check the following
     
     \begin{fact}
      \label{fact.action}
      \begin{enumerate}
       \item The set of fixed points for the action of the group $N$ (resp. $T$) on $\ein$ is exactly $\Delta$.
       \item For every $x \in \Omega_{\Delta}$, the $N$-orbit of $x$ is the leaf $F_{\Delta}(x)$
       \item The action of $T$ is free on $\Omega_{\Delta} \setminus F_{\Delta}(p_0)$, 
        and orbits of $T$ on this open set coincide with leaves of $\calfd$.
       \item On $ F_{\Delta}(p_0)$, orbits of $T$ are $1$-dimensional and coincide with the photons of 
       $C(p_0)$, with $p_0$ removed.
       
      \end{enumerate}

     \end{fact}

 In the rest of the paper, we will adopt the notations $\lieg_{\Delta}, \lies_{\Delta}, \lien, \mathfrak{t}$ for the Lie 
  subalgebras of $\oo(2,3)$ corresponding to the groups $G_{\Delta}, S_{\Delta}, N, T$.

\subsection{Standard Heisenberg algebras in $\oo(2,3)$}
\label{sec.standard}

The Lie group $N$ admits a Lie algebra  $\lien \subset \oo(2,3)$ that will be called
the {\it standard Heisenberg algebra } of 
 $\oo(2,3)$.  

It is not true that all subalgebras of $\oo(2,3)$ which are isomorphic to  $\heis(3)$ 
are conjugated to 
the standard algebra $\lien$.  There is however the following  useful characterization:

\begin{lemme}
 \label{lem.standard}
 Let $\lieh \subset \oo(2,3)$ be a Lie subalgebra isomorphic to $\mathfrak{heis}(3)$, and $H \subset \PO$ the corresponding 
 connected Lie subgroup.  Assume there exists 
 a nonempty open set of $\ein$ where the orbits of $H$ are $2$-dimensional and lightlike.  Then $\lieh$ 
  is conjugated in $\PO$ to the standard Heisenberg algebra $\lien$.
\end{lemme}

\begin{preuve}
 As any solvable Lie subalgebra of $\oo(2,3)$, $\lieh$ 
 must leave invariant a line $\RR.v$ or a 2-plane $P$ in $\RR^{2,3}$.  Such a  vector $v$ can not be
 timelike or spacelike, otherwise the decomposition 
  $\RR^{2,3}=\RR.v \oplus v^{\perp}$  would lead to an embedding of $\lieh$ in one of the Lie algebras 
  $\RR \oplus \oo(1,3)$ or $\RR \oplus \oo(2,2) \simeq \RR \oplus \mathfrak{sl}(2,\RR) \oplus \mathfrak{sl}(2,\RR)$.
   But none of those algebras contains a subalgebra isomorphic to $\heis(3)$.  Similarly,  $P$ can not be of
   signature $(+,+)$, $(+,-)$ or $(-,-)$, otherwise the decomposition $\RR^{2,3}=P \oplus P^{\perp}$ would lead to an embedding
     of $\lieh$ into $\oo(2) \oplus \oo(2,1) \simeq \RR \oplus \oo(1,2)$, 
     $\oo(1,1) \oplus \oo(1,2) \simeq \RR \oplus \oo(1,2)$ or $\oo(2) \oplus \oo(3) \simeq \RR \oplus \oo(3)$. One checks 
     as above that  this is not possible. The only possibilities are then:
     \begin{enumerate}[a)]
      \item The vector $v$ is lightlike or $P$ has signature 
     $(0,+)$ (resp. $(0,-)$).  This means that $H$ has a global fixed point in $\ein$, that 
     we can assume to be $p_0$ after conjugating into $\PO$.
     \item The form $Q^{2,3}$ vanishes identically on $P$, in which case $H$ has an invariant photon that we can assume to be $\Delta$.
     \end{enumerate}

 We first deal with  case $a)$.  After considering a stereographic projection of pole $p_0$,
 $\lieh$ becomes a subalgebra of $\mathfrak{Conf}(\RR^{1,2}) \simeq (\RR \oplus \oo(1,2)) \ltimes \RR^3$.  Here the normal
 subalgebra $\RR^3$ integrates into the subgroup of translations. Let us consider  the projection
  $\pi: (\RR \oplus \oo(1,2)) \ltimes \RR^3 \to \oo(1,2)$.  Since $\oo(1,2)$ does not have any subalgebra isomorphic to 
   $\heis(3)$ or $\RR^2$, the rank of $\pi_{|\lieh}$ is $0$ or $1$. Because $\RR \ltimes \RR^3$  (with $\RR$ acting by homothetic 
   transformations on $\RR^3$) does not contain a copy of $\heis(3)$, this  rank is actually $1$, hence the kernel of 
   $\pi_{|\lieh}$, denoted $\liea$,
   has dimension 
    $2$ in $\lieh$, hence is abelian.  The only subalgebras isomorphic 
    to $\RR^2$ in  $\RR \ltimes \RR^3$ are actually contained in $\RR^3$.  
     
     Our hypothesis on the orbits of the group $H$  implies that the translation vectors in $\liea$ span 
     a lightlike hyperplane,
      hence after conjugating into $\Conf(\RR^{1,2})$, we can assume $\liea=\mathfrak{t}$, where  $\mathfrak{t}$ was introduced 
      at the end of Section \ref{sec.symetries}.

     The first point of Fact \ref{fact.action} implies that since $H$ centralizes $\mathfrak{t}$, $H \subset G_{\Delta}$.
      The hypothesis on the orbits of $H$ says that on some open set, $H$-orbits and $T$-orbits coincide.  Points $3$ and $4$ of
       Fact \ref{fact.action} imply
      that the action of $H$ on $\Delta$ is trivial on some nonempty open set, hence trivial. This yields 
      $H \subset S_{\Delta}$.
        Because the normalizer of $\mathfrak{t}$ in $S_{\Delta}$ is $N$, we finally get $H=N$, and the proof is completed in this case.

 Consider now  case $b)$. Because $H$ leaves $\Delta$ invariant, $H$ is a subgroup of $G_{\Delta}$. As above, we can
  look at  the morphism 
  $$\pi: \lieg_{\Delta} \simeq (\RR \oplus \mathfrak{sl}(2,\RR)) \ltimes \lien \to \mathfrak{sl}(2,\RR).$$
  The same arguments as above show that
   the kernel of $\pi_{| \lieh}$ is a $2$-dimensional  abelian Lie subalgebra $\liea \subset \lieh$.  Observe that 
   $\liea \subset \lies_{\Delta}$, and the only $2$-dimensional 
   abelian subalgebras of 
   $\lies_{\Delta}$ are included in $\lien$. After conjugating into $G_{\Delta}$, we can ensure $\liea={\mathfrak t}$.  We then 
   finish the proof as in the first case.  
 
\end{preuve}



\section{The global geometry of manifolds without hyperbolic components}
\label{sec.global}

This section is devoted to establishing Theorem \ref{thm.topological}  in the only remaining case to be studied, 
namely that of closed $3$-dimensional Lorentz manifolds $(M,g)$ which are not locally homogeneous, such that $\mint$ does not admit any hyperbolic component, and 
 with a noncompact isometry group $\Iso(M,g)$.  By Theorem \ref{thm.conf.flat}, those manifolds are conformally flat, 
 and we consider $\delta: \tilde{M} \to \ein$ the corresponding developing map.  We  also recall the
  holonomy morphism $\rho: \operatorname{Conf}(\tilde{M}, \tilde{g})  \to \PO$. 
  
What we will really show in this section is:
\begin{theoreme}
\label{th.theoplat}
 Let $(M,g)$ be a closed, orientable and time-orientable, $3$-dimensional Lorentz manifold, 
 such that $\Iso(M,g)$ is noncompact.  We assume that
  $(M,g)$ is not locally homogeneous, and that $\mint$ does not admit any hyperbolic component.  Then:
  \begin{enumerate}
   \item The manifold $M$ is homeomorphic to a $3$-torus, or a parabolic torus bundle ${\bf T}_A^3$.
   
    \item There exists a metric $g'=e^{2\sigma}g$ in the conformal class of $g$ which is flat, and which is preserved by
   $\Iso(M,g)$.                                                                                    
   \item There exists a  
  smooth, positive, periodic 
   function $a:\RR \to (0,\infty)$ such that the universal cover 
  $(\tilde{M},{\tilde{g}})$ is isometric to $\RR^3$ endowed with the metric
  $$ {\tilde{g}}=a(v)(dt^2+2dudv).$$                                                                                                                  
   \item There is an isometric action of $\Heis$ on $(\tm,\tilde{g})$.                                                                                                                  
  \end{enumerate}

\end{theoreme}

This result clearly implies Theorem \ref{thm.topological} in the case under study.  Its proof will be the aim of 
Sections \ref{sec.approximately}
 to \ref{sec.conclusion} below.  In all those sections, $(M,g)$ satisfies the asumptions of Theorem \ref{th.theoplat}.

\subsection{Approximately stable foliation on ${M}$}
\label{sec.approximately}
So far, we saw that $(M,g)$ is an agregate of (possibly infinitely many)  components, the local geometry of 
which we understand fairly well.  But we need a global object which allows to understand how those component fit together. 
 This global object turns out to be a foliation provided  by the noncompactness of $\Iso(M,g)$ as follows.

  Consider a sequence $(f_n)$ in $\Iso(M,g)$  which tends to infinity,
   and call $AS(f_n)$  the subset of $TM$ comprising all vectors $v \in TM$ for which there exists a sequence $(v_n)$ in $TM$ converging to $v$, such
    that $|Df_n(v_n)|$ is bounded (where $|.|$ is the norm associated to an auxiliary Riemannian metric on $M$).  
    In \cite{zeghibisom1}, A. Zeghib proved
     the following result :
  \begin{theoreme}\cite[Theorem 1.2]{zeghibisom1}
   \label{theo.approximatively}
   Let $(M,g)$ be a closed Lorentz manifold, and $(f_n)$ a sequence of $\Iso(M,g)$ tending to infinity.  Replacing if necessary
    $(f_n)$ by a subsequence, the set  $AS(f_n)$ is a codimension $1$, lightlike, Lipschitz distribution in $TM$, which 
   integrates into a codimension $1$, totally geodesic, lightlike foliation.
  \end{theoreme}
   
   The foliation given by Theorem \ref{theo.approximatively} is called {\it the approximately stable foliation of $(f_n)$}.
   
   In the particular case of a $3$-dimensional manifold, codimension $1$, totally geodesic, lightlike foliations 
    have very nice 
   properties that were studied by A. Zeghib in \cite{zeghibfoliation}.  He proved in particular:
   
   \begin{theoreme}\cite[Theorem 11]{zeghibfoliation}
    \label{th.zeghibfol}
    Let $(M,g)$ be a $3$-dimensional closed Lorentz manifold.  Let $\calf$ be a $C^0$, codimension $1$, totally geodesic, 
    lightlike  foliation of $M$.  Then:
    \begin{enumerate}
     \item A leaf of $\calf$ is homeomorphic to a  plane, a  cylinder or a torus.
     \item The foliation $\calf$ has no vanishing cycles.
    \end{enumerate}

   \end{theoreme}

 We now choose a sequence $(f_n)$ tending to infinity in $\Iso(M,g)$, and after considering a suitable subsequence, we denote
 the approximatively stable foliation of $(f_n)$ by $\calf$. 
 By Theorem \ref{th.zeghibfol}, the leaves of $\calf$ are planes, 
 cylinders or tori.  Our main aim, and a decisive step to prove Theorem \ref{th.theoplat} will be to show that 
 {\it all leaves of $\calf$ are tori}, yielding the
   torus bundle structure of $M$.  It will be convenient in the sequel to consider the lift of $\calf$ to the universal cover $\tm$.  We will call
    $\tcalf$ this lifted foliation.

\subsection{Leaves of $\calf$ coincide with $\kiloc$-orbits on non locally homogeneous components}
\label{sec.calfetkill}

The aim of this section is to show: 

\begin{proposition}
 \label{prop.coincidencecalf}
 Let $\calm \subset \mint$ be a component which is not locally homogeneous.  Then $\kiloc$-orbits in $\calm$ coincide with leaves 
 of $\calf$.  In particular, any such component is saturated by leaves of $\calf$.
\end{proposition}

\subsubsection{The pullback foliation $\tcalfd$ and its geometric properties}
\label{sec.pulbackfolia}

We consider the developing map  $\delta : \tm \to \ein$, and take the pullback by $\delta$
of the foliation ${\calf}_{\Delta}$  defined in Section \ref{sec.foliation-einstein}.  We get in this way 
  a (singular) foliation $\tcalfd$ on $\tm$.  Actually, $\tcalfd$ is a genuine foliation by lightlike hypersurfaces
 on the open set $\tilde{\Omega}_{\Delta}=\delta^{-1}(\Omega_{\Delta})$.  
 Singularities occur on the  complement of $\tilde{\Omega}_{\Delta}$ in $\tm$, namely  
 $\tilde{\Delta}:=\delta^{-1}(\Delta)$. This singular set is either empty (in which case $\tcalfd$ is a regular foliation on $\tm$), 
  or  a $1$-dimensional
  lightlike manifold.  
  
  Let us emphasize the fact  that {\it a priori}, we don't have any  invariance property for  $\tcalfd$ under the action of 
  the fundamental group $\pi_1(M)$. 
   In particular, there is no reason for  $\tcalfd$  to define any foliation on $M$.
   
   In the following, we will identify $\oo(2,3)$ with the Lie algebra of conformal vector fields of $\ein$ (see Theorem 
   \ref{th.liouville}).
   We can pull back the vector fields of the Lie algebra  $\lien$  by the developing map $\delta:\tm \to \ein$, 
   getting a Lie algebra $\tlien$ 
 of conformal vector fields on $(\tm,{\tilde{g}})$.  By Fact \ref{fact.action}, the pseudo-orbits
 of $\tlien$ coincide with the leaves of 
  $\tcalfd$.
  \subsubsection{Foliation $\tcalfd$ and $\kiloc$-orbits}
  \label{sec.folorbits}
  A first important feature of the foliation  $\tcalfd$ is its relation to the $\kiloc$-orbits in $\tm$.  
  To see this, let us
   fix a parabolic component $\calm \subset \mint$ as in Proposition \ref{prop.coincidencecalf}.
We lift this component to the universal cover $\tm$ of $M$,
 and call $\tcalm$ a connected component of this lift.  
 
 \begin{lemme}
  \label{lem.coincidence}
  Replacing if necessary the developing map $\delta$ by $g \circ \delta$, for some $g \in \PO$, 
  the restriction to $\tcalm$ of any vector field of $\tlien$ is a Killing field for $\tilde{g}$. Conversely, any local
  Killing field defined on some open set $U \subset \tcalm$  is the restriction of a vector field in $\tlien$.
 \end{lemme}
  
  \begin{preuve}
  We pick $x \in \tcalm$ and $U$ a $1$-connected neighborhood of $x$ on which the 
 developing map $\delta$ is injective.  If $U$ is chosen small enough, the Lie algebra $\mathfrak{k}_U$ of 
 Killing fields on $U$ coincides with
  $\kiloc(x)$. Einstein's universe $\ein$ satisfies a generalization of Liouville's theorem:  Any  conformal Killing field 
  defined on some connected open set of $\ein$ is the restriction of a global one.  Thus the
   algebra $\delta_*(\liek_U)$ is a subalgebra of $\oo(2,3)$ isomorphic to the $3$-dimensional Heisenberg algebra.  
    The pseudo-orbits of $\delta_*(\liek_U)$ on $\delta(U)$ are $2$-dimensional and lightlike by the second point of Proposition 
    \ref{prop.heis}.  Lemma \ref{lem.standard} applies and says that post-composing $\delta$ by an element of 
    $\PO$, we may assume $\delta_*(\liek_U)=\lien$.  We thus get that any Killing field on $U$ is the 
    restriction of a vector field of $\tlien$. Since $\tlien$ and $\liek_U$ have same dimension, the restriction to $U$ 
    of {\it any} vector field 
    $X \in \tlien$  must be Killing.  Let us call  $Y=X_{|U}$.  
    Let us pick an arbitrary $y \in \tcalm$, 
   and  draw a
    simple curve $\gamma$ joining $y$ to $x$ inside $\tcalm$.  Let us consider $V$ a $1$-connected open neighborhood of $\gamma$ contained in 
    $\tcalm$ and containing $U$. Because the dimension of $\kiloc(z)$ is constant on $\tcalm$,  the vector field $Y$ can be extended 
    by  
    analytic continuation to  a Killing field (still denoted $Y$) defined on $V$. But now, $Y$ and $X_{|V}$ are two conformal
    Killing fields on $V$, which coincide on $U$. They must then coincide on $V$, 
    showing that $X$ is Killing in a neighborhood of $y$. The same dimentional argument as above shows that conversely, a Killing field
     defined in a sufficiently small neighborhood of $y$ is the restriction of a field in $\tlien$.
    
    
  \end{preuve}
 \begin{corollaire}
  \label{coro.nonsingulier}
  The component $\tcalm$ is included in $\tilde{\Omega}_{\Delta}$.
 \end{corollaire}
\begin{preuve}
 Points of $\tilde{\Delta}$ are singularities for the vector fields of $\tlien$.  
 Hence if a point $x \in \tilde{\Delta}$ belongs to $\tcalm$,
  Lemma \ref{lem.coincidence} will provide a Lie subalgebra of
  Killing fields vanishing at $x$ and isomorphic to ${\mathfrak{heis}}(3)$.  The isotropy representation then  
   yields  an embedding of Lie algebras ${\mathfrak{heis}}(3) \to \oo(1,2)$. This is impossible.
\end{preuve}

We conclude this paragraph with the following important lemma.

\begin{lemme}
 \label{lem.saturation1}
 Let $x$ be a point of  $\tcalm$, and  $\tilde{F}_{\Delta}(x)$ the leaf of $\tcalfd$ through $x$.  Then $\tilde{F}_{\Delta}(x)$ 
  is included in $ \tcalm$, and  coincides
 with the $\kiloc$-orbit of $x$.
\end{lemme}

\begin{preuve}
 Let us consider a leaf $\tFd$ having a nonempty intersection with $\tcalm$.  Assume for a contradiction that 
 $V=\tFd \cap \tcalm$ is not all of  $\tFd$. It means that $V$ is an open subset of $\tFd$ having a nontrivial boundary 
 $\partial V$ inside $\tFd$.  Of course, $\partial V \subset \partial \tcalm$  (this last boundary is taken in $\tm$). Since 
 $\tFd$ is a pseudo orbit of $\tlien$, it is easy to show that there exists $y \in \partial {V}$, a vector field 
 $X \in {\tlien}$ and a point $x \in {V}$ such that the local orbit  $t \mapsto \phi_{X}^t.x$ is defined on $[0,1]$,
  $\phi_{X}^t.x$ belongs to ${V}$ for $t \in [0,1/2)$ but $\phi_{X}^{1/2}.x \in \partial {V}$. We denote by $\hat{R}$ the bundle of 
  frames on $\tm$, and exceptionnaly in this proof, we adopt the notation $\hm$ for the bundle of orthonormal
  frames of $\tm$ (and not of $M$).
  The local action of $\phi_{X}^t$ lifts naturally to $\hat{R}$.
  We  pick $\hx \in \hm$  in the fiber of $x$, and 
  look at the orbit $t \mapsto \phi_{X}^t.\hx$  in $\hat{R}$. 
  Because $X$ is Killing on $\tcalm$ (Lemma \ref{lem.coincidence}), 
      this orbit is contained
       in $\hm$ for $t \in [0,1/2)$, and the same is true for  $t \in [0,1/2]$ because $\hm$ is closed in $\hat{R}$.
         We now look at the generalized curvature map $\dk : \hm \to {\mathcal{W}}$, and its derivative that we see as 
         a map $D \dk : \hm \to Hom(\lieg, \mathcal{W})$. 
         The map $t \mapsto D \dk(\phi_{X}^t.\hx)$ makes sense for $t \in [0,1/2]$, 
         and is constant on this interval because $X$ is Killing on $\tcalm$.  In particular, 
         the kernel of $D \dk(\phi_{X}^t.\hx)$ is the same for all $t \in [0,1/2]$, hence the rank of $\dk$ is the same
          at $\hx$ and at $\phi_{X}^{1/2}.\hx$.  We get  that the rank of $\dk$ at $\phi_{X}^{1/2}.x$ is $3$, but
           we  already observed in the proof of Lemma \ref{lem.rank3}, that all points where $\dk$ has rank $3$ are
         contained in $\tm^{\rm int}$.
          We infer $\phi_{X}^{1/2}.x \in \tm^{\rm int}$,  contradicting $\phi_{X}^{1/2}.x \in \partial \calm$.
   
   The last part of the lemma follows easily. Lemma \ref{lem.coincidence}, together with Corollary \ref{coro.nonsingulier}
    ensures that for every  $x \in \tcalm$, the $\kiloc$-orbit of $x$ coincides with $\tFd(x) \cap \tcalm$.
     But $\tFd(x) \cap \tcalm=\tFd(x)$ by the first part of the proof.  
\end{preuve}

\subsubsection{Proof of Proposition \ref{prop.coincidencecalf}}     
     
 We keep the notations of the previous paragraph.  
 We also lift the foliation $\calf$ to a foliation $\tcalf$ on the universal cover $\tm$.  For each $x \in \tm$, we 
 denote by $\tF(x)$ the leaf of $\tcalf$ containing $x$.
  
Thanks to Lemma \ref{lem.saturation1}, Proposition  \ref{prop.coincidencecalf} will be a simple consequence of:



\begin{lemme}
 \label{lem.saturation2}
 For every $x \in \tcalm$, one has $\tFd(x)= \tF(x)$.  
 
\end{lemme}

\begin{preuve}
 We work on $\tcalm$, and we consider the two $1$-dimensional lightlike distributions 
 $\tcald_{\Delta}=T\tcalfd^{\perp}$  and  $\tcald=T\tcalf^{\perp}$.
 Our aim is to show that those distributions  coincide on $\tcalm$.  For every $x \in \tcalm$, 
 let us introduce the set 
 ${\mathcal C}(x)$, comprising all lightlike directions $u \in {\mathbb P}(T_x\tm)$ such that there exists a lightlike totally geodesic 
 hypersurface $\Sigma$ through $x$, with $T_x\Sigma^{\perp}=u$.  Let us recall a key observation made in \cite{zeghibisom2}:
 
 \begin{lemme}\cite[Proposition 2.4]{zeghibisom2}
 \label{lem.sec.curvature}
 If the set ${\mathcal C}(x)$ spans $T_x\tm$, then the sectional curvature at $x$ is constant.
 \end{lemme}

 If at some point $x$ of $\tcalm$, the directions $\tcald_{\Delta}$ and $\tcald$ do not coincide, then
  Lemma \ref{lem.sec.curvature} implies that they must both be fixed by the local flow generated by the isotropy algebra 
  ${\mathfrak Is}(x
  )$.  But a  nontrivial parabolic $1$-parameter flow in $\OO(1,2)$ has only one invaraint direction: Contradiction.

    We are thus led to  $\tFd(x) \subset \tF(x)$ for every $x \in \tcalm$.  
     This  inclusion can not be proper, otherwise $\tcalfd$ could be extended in a smooth way to points of $\tilde{\Delta}$.
\end{preuve}

  \begin{remarque}
   \label{rem.codimension1}
   The previous proof  shows actually that on a non locally homogeneous parabolic component $\tcalm$, any lightlike, 
   totally geodesic, codimension $1$ 
    foliation has to coincide with $\tcalf$. 
  \end{remarque}

\subsection{Existence of toral leaves for ${\mathcal F}$}
\label{sec.toral}
  We keep going in the geometric study of the foliation $\calf$  (and simultaneously 
  in the understanding of the developing map $\delta$), by proving the existence of one toral leaf for $\calf$.
\subsubsection{Injectivity properties of $\delta$}
\label{sec.injectivity-leave}
We keep the notations of Section \ref{sec.pulbackfolia}, and recall the open set $\tilde{\Omega}_{\Delta}$ 
where the foliation $\tcalfd$ is defined.
\begin{lemme}
 \label{lem.injective}
 Let $x \in \tilde{\Omega}_{\Delta}$, and assume  that the leaves $\tF(x)$ 
 and $\tFd(x)$ coincide.  Then $\delta$ is injective in restriction to $\tF(x)$.
\end{lemme}

\begin{preuve}
 Considering if necessary a finite cover of $M$ (what won't change $\tm$), there exists  $W$ a vector field on $M$,
 tangent to $\calf$  and satisfying $g(W,W)=1$.
   We lift $W$ to a  vector field $\tilde{W}$ on $\tm$, which is tangent to $\tcalf$.  Notice that $\tilde{W}$ is complete.
    By assumption,  $\tilde{W}$ is tangent to $\tF(x)$. The proof follows now closely the arguments of 
    \cite[Proposition 6.5]{zeghibflot}. Let us $\tilde{\cald}$ the $1$-dimensional foliation integrating $\tcalf^{\perp}$.
     A fundamental remark made in \cite[Proposition 2]{zeghibfoliation} is that because 
     $\tF(x)$ is totally geodesic, any vector field $U$ tangent to $\tcald$ acts as a Killing field on the degenerate
      surface $(\tF(x),\tilde{g})$.  It follows easily that if $\gamma$ and $\eta$ are two curves on $\tF(x)$ parametrized 
      by $[0,1]$, such
       that $\gamma(0)$ and $\eta(0)$ belong to the same leaf of $\tcald$, and if $\gamma$ and $\eta$ have the same length with respect to 
       $\tilde{g}$, then $\gamma(1)$ and $\eta(1)$ also belong to a same leaf of $\tcald$.
         Applying this remark to the integral curves of the flow  $\{ \psi^t \}$ generated by $\tilde{W}$, 
          we obtain that $\psi^t$ maps leaves of $\tilde{\cald}$ to leaves of $\tilde{\cald}$.
            Given $\tD_0$ a leaf of
     $\tcald$ in $\tF(x)$, the union ${\mathcal U}(\tD_0)=\bigcup_{t \in \RR} \psi^t(\tD_0)$ is open in $\tF(x)$, 
     and two such open sets either coincide, or are disjoint, 
     so that $\tF(x)={\mathcal U}(\tD_0)$. By hypothesis, the leaves 
     $ \tFd(x)$ and $\tF(x)$ coincide, so that the developing map $\delta$ sends $ \tF(x)$ to 
     $F_{\Delta}=F_{\Delta}(\delta(x)) \subset \ein$. Let $\gamma: I \to \tF(x)$ be an injective parametrization of the leaf of 
     $\tcald$ through $x$.  We observe that for every $s \in I$, $\delta$ is injective on the curve 
     $t \mapsto \psi^t(\gamma(s))$, because in $F_{\Delta}$, there is no closed curve transverse to photons of $F_{\Delta}$.
       Also, for every $t \in \RR$, $\delta$ is injective in restriction to $s \mapsto \psi^t(\gamma(s))$, because no 
       photon in $F_{\Delta}$ is closed.  The injectivity of $\delta$ on $\tF(x)$ follows.
   \end{preuve}

\subsubsection{The group $\Iso(M,g)$ is not a torsion group}
\label{sec.notorsion}

 Our noncompactness hypothesis on the group $\Iso(M,g)$ does not prevent {\it a priori} $\Iso(M,g)$
 from being a torsion group. 
  In particular, we still don't know if  there exists a single element 
 $h \in \Iso(M,g)$ such that $\{ h^k \}$ is infinite discrete. The aim of this paragraph is to show it is indeed the 
 case, and 
 to prove the stronger statement:
 
 \begin{proposition}
  \label{prop.pastorsion}
  Let $\calm \subset \mint$ be a component which is not locally homogeneous.  Let $F \subset \calm$ be a leaf of $\calf$ 
  containing at least one recurrent point. Let $S_F$ be the stabilizer of $F$ in $\Iso(M,g)$.
 There exists $h \in S_{F}$ such that the group $\{ h^k \}$
  is not relatively compact in $\Iso(M,g)$.
 \end{proposition}

Recall (Proposition \ref{prop.coincidencecalf}) that non locally homogeneous components of $\mint$ are saturated by 
the leaves of $\calf$.

 The proof of Proposition \ref{prop.pastorsion} will require the  intermediate lemmas 
 \ref{lem.stab-ferme} and \ref{lem.holonomy} below.
  We  lift $F$ to a leaf $\tF \subset \tcalm$ and call 
  ${S}_{\tF}$ the stabilizer of 
  $\tF$ in $\Iso(\tm,\tilde{g})$.  Observe that ${S}_{\tF}$  projects surjectively on $S_F$ under the epimorphism 
  $\Iso(\tm, \tilde{g}) \to \Iso(M,g)$.
 \begin{lemme}
  \label{lem.stab-ferme}
  For every leaf $F \subset \calm$ containing recurrent points, the groups $S_F$ and $S_{\tF}$ are closed, 
  noncompact subgroups of \,$\Iso(M,g)$
  and \,$\Iso(\tm,\tilde{g})$ respectively.
 \end{lemme}

 \begin{preuve} We  first prove that $S_F$ is closed in $\Iso(M,g)$.  
  If $(f_k)$ is a sequence of $S_F$ which converges to $f_{\infty} \in \iso(M,g)$, then for $k$ very large, 
  $f_k^{-1}f_{\infty}$ belongs to the identity component $\Iso^o(M,g)$.
   Because $F$ coincides with a $\kiloc$-orbit of 
 $\calm$ (Lemma \ref{lem.saturation2}), we thus have $f_k^{-1}f_{\infty}(F)=F$, which in turns implies 
 $f_{\infty}(F)=f_k(f_k^{-1}f_{\infty}(F))=f_k(F)=F$. 
  
  Let us now check that  $S_F$ is noncompact. By assumption on $F$, there is a recurrent point  $x$ in $F$. It means that
    there exists
 a sequence $(f_k)$ tending to infinity in  $\Iso(M,g)$ such that $f_k(x) \to x$. 
 Because the $\operatorname{Is}^{loc}$-orbit of $x$ is a $2$-dimensional submanifold, 
 the connected components of which are $\kiloc$-orbits (see Theorem \ref{thm.integrabilite}), we get that 
 $f_k(x) \in F$ for $k$ large enough.  In particular, 
   $S_F$ is a noncompact subgroup of  $\Iso(M,g)$.  
  
  The corresponding assertions on $S_{\tF}$ are then straigthforward.
\end{preuve}

  \begin{lemme}
   \label{lem.holonomy}
  Let $F \subset \calm$ be a leaf of $\calf$, and $\tF$ a lift of $F$ to $\tm$.
  \begin{enumerate}
   \item The holonomy morphism $\rho$ maps the group ${S}_{\tF}$ into the group $S_{\Delta}$.  In particular, any element
    of $S_{\tF}$ leaves invariant the leaves  of $\tcalf$ which are sufficiently close to $\tF$.
   \item The morphism $\rho : S_{\tF} \to S_{\Delta}$ is injective and  proper.
  \end{enumerate}

  \end{lemme}

\begin{preuve}
 We heavily use the notations introduced in Section \ref{sec.more.geometry}. 
 We choose a  transversal $I \subset \tcalm$ to the foliation $\tcalf$, that cuts $\tF$ at $x$. We assume
  that $I$ is small enough, so that  $\delta$ sends $I$ 
 injectively on a transversal $J$ of ${\mathcal F}_{\Delta}$.  Shrinking $I$ if necessary, $J$ meets each leaf of 
 ${\mathcal F}_{\Delta}$ at most once, so that by Lemma \ref{lem.injective},  $I$ meets each leaf of $\tcalf$ at most once.
  We call $\mathcal{V}$ the open subset obtained by saturating $I$ by leaves of $\tcalf$.  By what we just said, $I$ is 
   the space of leaves of $\mathcal{V}$, and 
  the map $\varphi:=\pi_{\Delta} \circ \delta: I \to \Delta$ gives an identification of $I$ with $J'=\pi_{\Delta}(J)$, 
  the space of leaves of the foliation induced by ${\mathcal F}_{\Delta}$ on $\delta({\mathcal V})$.
    Under this identification, the point $x$ is sent to a point $p \in J'$. 
    
    Let $\tth$ be an element of $S_{\tF}$.  Its action on the space of leaves of $\tcalf$ 
    yields a germ $\overline{h}$ of diffeomorphism of $I$ fixing $x$. The equivariance relation
     $\delta \circ \tth=\rho(\tth) \circ \delta$  shows that $\rho(\tth)$  permutes leaves of $\calfd$  near $\delta(\tF)$.
      In particular, $\rho(\tth)$  maps $J'$ to an interval of $\Delta$ containing $p$.  We infer that $\rho(\tth)$ 
      preserve $\Delta$, what yields $\rho(\tth) \in G_{\Delta}$.  Moreover, denoting 
      $l : G_{\Delta}\simeq \PGL(2,\RR) \ltimes N \to \PGL(2,\RR)$, we get the equivariance relation
     $\varphi \circ \overline{h}=l(\rho(\tth)) \circ \varphi$.  Now, $l(\rho(\tth))$
      acts as an element of $\PGL(2,\RR)$ on $\Delta$, admitting  $p$ as fixed point.  We know the local dynamics of a 
      M\"obius  transformation around one of its fixed points:  If $l(\rho(\tth))$ is nontrivial, we can choose
       $q \in J'$, $q \not = p$,  such that $l(\rho({\tth}^k))(q)$  belongs to $J'$ for all $k \geq 0$, and 
       $\lim_{k \to \infty}l(\rho({\tth}^k))(q)=p$.  This means that if $\tF'$ is a leaf corresponding to $\varphi^{-1}(q)$, 
        the iterates ${\tth}^k(\tF')$ will accumulate on $\tF$. But $\tF'$ is a $\kiloc$-orbit 
        by Proposition \ref{prop.coincidencecalf}, and closeness of the $Is^{loc}$-orbit of $\tF'$ in $\tm^{\rm{int}}$ 
        (Corollary 
        \ref{coro.islocorbites}) says that $\tF$ and all the $h^k(\tF')$, $k \in \NN$, belong to the same $Is^{loc}$-orbit.
         This accumulation phenomenon then contradicts  the fact that $\operatorname{Is}^{loc}$-orbits are  
  submanifolds in $\tm^{\rm{int}}$ (see Theorem \ref{thm.integrabilite}).  We conclude that $l(\rho(\tth))$ is trivial, which implies
   that $\rho(\tth) \in S_{\Delta}$.  Moreover, $\overline{h}$ is trivial, 
   which means that all leaves of $\tcalf$ close to 
    $\tF$ are left invariant by $\tth$.
    
 We now prove the second point of the Lemma.   Let $\tth \in S_{\tF}$ such that $\rho(h)=id$.  Equivariance relation
  $\delta \circ \tth = \rho(\tth) \circ \delta$, together with Lemma \ref{lem.injective}, shows that the action of $\tth$
  on $\tF$.
  The following  fact then implies $\tth=id$.
   \begin{fact}
\label{fact.injectiveproper}
Let $N$ be a Lorentz manifold and $\Sigma \subset N$ a lightlike hypersurface.  denote by $S_{\Sigma}$ the stabilizer of 
$\Sigma$ in $\Iso(N)$. Then the restriction map $r: S_{\Sigma}\to \mathrm{Homeo}(\Sigma)$ is injective and proper. 
 \end{fact}

\begin{preuve}
The proof relies on the fact that the map, which to an isometry associates its $1$-jet  at a given 
point, is injective and proper, and that 
restricting elements of $\OO(1,n-1)$ to a lightlike hyperplane is also injective and proper.  
Details can be found  
in \cite[Prop 3.6]{zeghibflot}, for instance.
\end{preuve}

Properness of the map $\rho:S_{\tF} \to S_{\Delta}$ follows the same lines.  If $(\tth_k)$ is a sequence of $S_{\tF}$  such that
 $\rho(\tth_k)$ is relatively compact in $S_{\Delta}$.  Then $\rho(\tth_k)_{|\delta(\tF)}$ is relatively compact, hence 
 the retriction of ${\tth_k}$ to $\tF$ is relatively compact by Lemma \ref{lem.injective}.  Fact \ref{fact.injectiveproper} yields that $(\tth_k)$
  is relatively compact in $S_{\tF}$.

\end{preuve}

We can now proceed to the proof of Proposition \ref{prop.pastorsion}.  
We know from Theorem \ref{th.zeghibfol}, that the leaves of $\calf$ are discs, cylinders or tori, 
and there are no vanishing circles.  It means that 
 the leaf $\tF$ is a disc, and the stabilizer $\Gamma_F$ of $\tF $ in $\pi_1(M)$ is either trivial, or 
 a discrete subgroup isomorphic to $\ZZ$ or 
 $\ZZ^2$.  On the other hand, we also know that $\tilde{S}_F/\Gamma_F$ is noncompact, because of Lemma 
 \ref{lem.stab-ferme}.

 \begin{enumerate}[a)]
  \item {\it Case where $F$ is a disk.}
  We choose a nontrivial $\tth \in S_{\tF}$.  It restricts
  to a nontrivial transformation of $\tF$ (Fact \ref{fact.injectiveproper}), hence $\rho(\tth)$ 
  is a nontrivial element of $S_{\Delta}$, by Lemma
   \ref{lem.holonomy}.  Every nontrivial element of $S_{\Delta}$ generates an infinite discrete group.  
    In particular, $\{\rho(\tth)^k\}_{k \in \ZZ}$
   is  not  relatively 
   compact in 
   $S_{\Delta}$.  Second point  of Lemma \ref{lem.holonomy} says that $\{h ^k \}$ is not relatively compact in 
   $\Iso(\tm,\tilde{g})$. Fact \ref{fact.injectiveproper} thus implies that $\{{\tth}^k_{|\tF}\}$ is not relatively compact
    in $\operatorname{Homeo}(\tF)$,
    hence  the same is true  for $\{h_{|F}^k \}$, because the projection $\pi: \tF \to F$ is  a diffeomorphism in the case
     we are considering. 
     Finally, $\{ h^k \}$ is not relatively compact in $\Iso(M,g)$.
   
  \item {\it Case where $F$ is a cylinder.}  Because $\calf$ does not have vanishing cycles, $\Gamma_F:=S_{\tF} \cap \Gamma$
    is nontrivial, generated by a single element ${\gamma}$.  After considering an index $2$ subgroup of 
    $S_{\tF}$, we may assume that $\gamma$ is centralized by all elements of $S_{\tF}$. We observe that $\rho(\gamma)$ is nontrivial
     by Lemma \ref{lem.holonomy}, and consider its  centralizer in $S_{\Delta}$.
      Two cases can then occur: 
     
     - The group $\rho(S_{\tF})$ is contained in a $1$-paramater subgroup of $S_{\Delta}$. In this case,
      $\rho(S_{\tF})/<\rho(\gamma)>$ is relatively compact in $S^1$. This implies that
      $(S_F)_{|F}$ is relatively compact, hence $S_F$ is relatively compact in $\Iso(M,g)$ 
      (again Fact \ref{fact.injectiveproper}).  This is ruled out by Lemma \ref{lem.stab-ferme}.
      
      - If we are not in the previous case, we can find $\tth \in S_{\tF}$ such that the group generated by $\rho(\tth)$ and $\rho(\gamma)$
       is discrete isomorphic to $\ZZ^2$.  As above, applying second point of Lemma \ref{lem.holonomy} and Fact 
       \ref{fact.injectiveproper}, one gets that $\tth$ projects to   $h \in S_F$, such that 
        $\{ h^k \}$ is infinite discrete in $\Iso(M,g)$.
        
  \item {\it Case where $F$ is a torus.}  This time, $\Gamma_{F}$ is isomorphic to $\ZZ^2$ and generated by $\gamma_1$ 
  and $\gamma_2$.  Lemma \ref{lem.holonomy} ensures that $\tau_1:=\rho(\gamma_1)$ and   $\tau_2:=\rho(\gamma_2)$
   generate a discrete subgroup of $S_{\Delta}$ isomorphic to $\ZZ^2$.  Such a subgroup must be included in  $N$, and after 
   conjugating 
   $\rho$ into $G_{\Delta}$ (what amounts to post-compose $\delta$ by some element of $G_{\Delta}$), 
   we have that $<\tau_1,\tau_2> \subset T$. We must have  $\rho(S_{\tF}) \subset N$ because $S_{\tF}$ normalizes
   $\Gamma_F$, and because $S_F$ is noncompact, 
   $\rho(S_{\tF}) \not \subset T$.  Picking $\tth \in S_{\tF}$ such that $\rho(\tth) \not \in T$, we get an element $h \in S_F$ which,
    by similar 
   arguments as above, 
      generates an infinite discrete group $\{  h^k \} \subset \Iso(M,g)$.
   
 \end{enumerate}

\subsubsection{Existence of a toral leaf} 
 \label{sec.toralproof}
 We now consider an element $h \in \Iso(M,g)$ given by Proposition \ref{prop.pastorsion}, namely $\{ h^k \}$ is not relatively 
 compact in $\Iso(M,g)$.  Theorem \ref{theo.approximatively} provides an approximately stable foliation 
 ${\mathcal F}_h$ associated to
  a subsequence of $\{ h^k \}$, and since all what we did before did not assume anything special on ${\mathcal F}$, we can 
  decide that now $\calf={\mathcal F}_h$.
 
 \begin{proposition}
  \label{prop.existence.tori}
  Every  leaf $F$ of $\calf$ which is included in $\calm$  
  is a torus.
 \end{proposition}

 Let $F$ be a leaf of $\calf$ included in $\calm$, such that almost every point of $F$ is recurrent for $\{  h^k \}$.
 We lift $F$ to $\tF \subset \tcalm$, and we will also assume that $\delta(\tF)$ is not included in the leaf 
 $F_{\Delta}(p_0)$ (see Fact \ref{fact.action}). 
 Observe that since $M$ is closed, Poincar\'e recurrence ensures that 
 almost every point of $(M,g)$ is recurrent for $\{ h^k \}$.  It follows that for almost every leaf of $\calf$, almost 
 every point is recurrent (leaves of $\calm$ coincide with $\kiloc$-orbits hence are locally closed, so there is nothing
 tricky in 
 desintegrating the volume form in $M$  along those leaves).  Hence almost every leaf of $\calf$ in $\calm$ is an $F$ with 
 the properties above.
 We claim that $F$ is a torus. To see this, we lift  $h$ to an element 
  $\tth \in S_{\tF}$, and our assumption is that the set of  points in $\tF$, which are recurrent under 
   the group $<\tth, \Gamma_{F} >$  have  full measure in $\tF$. By Lemma \ref{lem.injective}, almost all points of 
   $\delta(\tF)$ must be recurrent  
    for the group $\rho(<\tth,\Gamma_{F}>)$. We now go back to the analysis made at the end of 
    Section \ref{sec.notorsion}.  
   If $F$ is not a torus,  we are in cases $a)$ or $b)$ 
    of this discussion.  In case $a)$, $\Gamma_F$ is trivial.  The action of $\rho(\tth)$ on $\delta(\tF)$  is
     conjugated to that of an affine transformation on (an open subset of) the plane (see Section \ref{sec.more.geometry} and formula 
     (\ref{eq.formuleaffine})).  The set of recurrent points of $\rho(\tth)$ has
      thus zero measure on $\delta(\tF)$. A contradiction.
      
      In case $b)$, we saw that the group 
     $\rho(<\tth,\Gamma_{F}>)$ is conjugated to a lattice in the closed subgroup $T$. Since by assumption, $\delta(\tF)$ is
     not included in  the leaf 
     $F_{\Delta}(p_0)$, point $3)$ of Fact \ref{fact.action} shows that $T$ has no recurrent point on $\delta(\tF)$. 
      We reach a new contradiction.

     The arguments above show that almost all leaves $F \subset \calm$ are tori.  Now, for a codimension $1$ foliation 
     on a closed manifold, the union 
     of all compact leaves is itself compact (see \cite[Chap II, Corollary 3.10]{godbillon}). 
     Proposition \ref{prop.existence.tori} follows.


\subsection{All leaves of ${\mathcal F}$ are tori}

We keep the notations of the last section.  We still consider $h \in \Iso(M,g)$ such that $\{ h^k \}$ is 
not relatively compact.  We consider the  approximately stable  foliation $\calf$ associated to some 
sequence $(h^{n_k})$, with 
$n_k \to \infty$.  Proposition \ref{prop.existence.tori} and its proof  provide 
us with a  leaf $F_0$ which is an $h$-invariant torus, and is included in a non locally homogeneous component $\calm$.
We lift $F_0$ to $\tF_0 \subset \tm$, and $h$ to $\tth \in \Iso(\tm,\tilde{g})$
 preserving $\tF_0$. The group $\Gamma_0=S_{\tF_0} \cap \pi_1(M)$ is discrete and 
 isomorphic to $\ZZ^2$, generated by two elements $\gamma_1$ and $\gamma_2$. Lemma \ref{lem.holonomy} ensures that 
 $\tau_1:=\rho(\gamma_1)$ and   $\tau_2:=\rho(\gamma_2)$ 
   generate a discrete subgroup of $S_{\Delta}$ isomorphic to $\ZZ^2$. After conjugating into
   $G_{\Delta}$, we may assume that  $\tau_1$
    and $\tau_2$ are elements of $T$. 
    
    We call in the sequel $\tilde{H}$ the subgroup generated by $\tth, \gamma_1$ and 
    $\gamma_2$. Lemma \ref{lem.holonomy} ensures that $\rho(\tth) \in S_{\Delta}=\RR_+^* \ltimes N$  
    (see Section \ref{sec.symetries}). Elements of $S_{\Delta}$ which are not in $N$ act on $N$ with nontrivial dilation.
     They can not preserve any lattice in $T$.    This says that because $\tth$ normalizes  $\Gamma_0$, 
     we must have $\rho(\tilde{H}) \subset N$.  
    There are thus three elements $X,Y,Z$
     in the Lie algebra $\lien$ such that $\tau_1=e^X$, $\tau_2=e^Z$ and $\rho(\tth)=e^Y$. 
     The center of $\lien$ is included in $\operatorname{Span}(X,Z)$, and we pick $Z_0 \not =0$ in this center.  
    
    We now pull back the four
      vector fields $X,Y,Z,Z_0$ on $\ein$ by the developing map $\delta: \tm \to \ein$.  This way, we get  vector fields
       $\tX,\tY, \tZ, \tZ_0$ on $\ein$.  Observe that $\tau_1^*X=X$, $\tau_2^*Z=Z$, $\rho(\tth)^*Y=Y$, and 
       $\rho(\tth)^*Z_0=\tau_1^*Z_0=\tau_2^*Z_0=Z_0, \ \tau_1^*Z=\tau_2^*Z=Z$
       imply the relations
       
       \begin{equation}
        \label{eq.vectorfields}
        \gamma_1^*\tX=\tX,\, \gamma_2^*\tZ=\tZ,   \ \tth^*\tY=\tY, \ {\rm and} \ 
        \tth^*\tZ_0=\gamma_1^*\tZ_0= \gamma_2^*\tZ_0=\tZ_0, \  \gamma_1^*\tZ= \gamma_2^*\tZ=\tZ
       \end{equation}
on $\tm$.
 After introducing those notations, we can state what will be the last technical step of our study:
 
 \begin{proposition}
  \label{prop.ensembleE}
For every $x \in \tm$, we have:
\begin{enumerate}[i)]
 \item{The point $x$ belongs to $\tilde{\Omega}_{\Delta}$ and $\tilde{F}(x)=\tilde{F}_{\Delta}(x)$.}
 \item{The leaf $\tilde{F}(x)=\tilde{F}_{\Delta}(x)$ is $\tilde{H}$-invariante.}
 \item{The restriction of the $3$ vector fields $\tX,\tY,\tZ$ are complete on $\tilde{F}_{\Delta}(x)$, and the equalities
  $\phi_{\tilde{X}}^1=\gamma_1, \ \phi_{\tilde{Z}}^1=\gamma_2, \ \phi_{\tilde{Y}}^1=\tilde{h}$ hold on $\tilde{F}_{\Delta}(x)$.}
\end{enumerate}
\end{proposition}
The proposition will show that $\tcalf$ coincide $\tcalfd$ on $\tm$, and
 $\Gamma_{0}$-invariance of the leaves  of $\tcalf$ will easily imply that leaves of $\calf$ are all tori.  In the next section 
 \ref{sec.conclusion} we will derive more consequences from this equality ${\mathcal E}=\tm$, and prove Theorem 
 \ref{th.theoplat}.

We are going to consider    the set $\mathcal{E} \subset \tm$, 
 comprising all points $x \in \tm$  satisfying the three conditions of Proposition \ref{prop.ensembleE}, and show that 
  $\mathcal{E}$ is nonempty, open and closed in $\tm$, yielding  $\mathcal{E}=\tm$.

 \subsubsection{The set $\mathcal{E}$ is nonempty}
 \label{sec.nonempty}
 We check here that every point $x \in \tF_0$ belongs to $\mathcal{E}$. Recall $\rho: \operatorname{Conf}(\tm) \to \PO$
  the holonomy morphism.
  
\begin{lemme}
 \label{lem.feuillecocompacte}
 Let $x \in \tilde{\Omega}_{\Delta}$ such that $\tF(x)=\tFd(x)$. Assume moreover  that $\tFd(x)$ is invariant by a 
  subgroup  $\Lambda \subset \pi_1(M)$, isomorphic to $\ZZ^2$ and such that $\rho(\Lambda) \subset T$.
  Then the map $\delta$ is a diffeomorphism from $\tFd(x)$ to $F_{\Delta}(\delta(x))$.  Moreover 
  $F_{\Delta}(\delta(x)) \not =F_{\Delta}(p_0)$.
\end{lemme}
  
  \begin{preuve}
    Lemma \ref{lem.injective} ensures that 
  $\delta$ is a diffeomorphism from $\tFd(x)$ to an open subset $U \subset F_{\Delta}(\delta(x))$. 
  The group $\Lambda$ is isomorphic 
   to $\ZZ^2$ and 
   acts properly discontinuously on the disk $\tF(y)=\tFd(y)$.  
   By a cohomological dimension argument, the action must be cocompact.  
    The group $\rho(\Lambda)$ is thus a lattice in $T$, and must act  properly and cocompactly on $U$. Last point of Fact 
  \ref{fact.action} says that the action of $\rho(\Lambda)$ can not be  proper on any open subset of $F_{\Delta}(p_0)$.
   We thus infer that $F_{\Delta}(\delta(x)) \not = F_{\Delta}(p_0)$.  In particular, again by Fact \ref{fact.action}, 
    the action of $\rho(\Lambda)$ is proper and
    cocompact on $F_{\Delta}(\delta(x))$.  We then must  have $U=F_{\Delta}(\delta(x))$. 
  \end{preuve}

  The completeness of $\tX,\tY$ and $\tZ$ on
     $\tF_0$ follows from Lemma \ref{lem.feuillecocompacte}, applied for $\Lambda=\Gamma_0$, because $X,Y,Z$ are complete on leaves of $\calfd$.
     The relations $\tau_1=e^X$, $\tau_2=e^Z$ and $\rho(\tth)=e^Y$ imply 
     that the relations $\phi_{\tilde{X}}^1=\gamma_1, 
     \ \phi_{\tilde{Z}}^1=\gamma_2, \ \phi_{\tilde{Y}}^1=\tilde{h}$ hold on $\tF_0$. 
     We infer that $\tF_0 \subset {\mathcal E}$.
 
 \subsubsection{The set $\mathcal{E}$ is open}
 \label{sec.open}
 
 We begin by stating a lemma that we will use repeatedly in the sequel.
 
 \begin{lemme}
  \label{lem.coincide}
  Let $U \subset \tilde{\Omega}_{\Delta}$ be a connected open set. 
  Let $f,g:U \to \tm$ two conformal immersions.  Assume that for some $x \in \tilde{\Omega}_{\Delta}$,
   $\tF_{\Delta}(x) \cap U \not = \emptyset$, and that $f$ and $g$ coincide on $\tF_{\Delta}(x) \cap U$, then $f$ and $g$ coincide on $U$.
 \end{lemme}

 \begin{preuve}
  Shrinking $U$ if necessary and looking at $\delta(U) \subset \ein$, we are reduced to the situation 
  of two transformations $g_1$ and $g_2$ of $\PO$  which coincide on some open subset of a lightcone in $\ein$.  At level of 
  linear algebra, it means that those two transformations of $\PO$  must  coincide on a lightlike hyperplane of 
  $\RR^{2,3}$.  This easily implies $g_1=g_2$.  
 \end{preuve}

 Let us start with $x \in \mathcal{E}$. Vector fields $\tX,\tY$ and $\tZ$ are complete 
 in restriction to $\tFd(x)=\tF(x)$, so given $\epsilon >0$, 
   we can choose $U \subset \tilde{\Omega}_{\Delta}$  a  small neighborhood of $x$ such that $\phi_{\tX}^t.y$ is defined 
   on $[-\epsilon, 1+\epsilon]$ for every $y \in U$, and for every $t \in [- \epsilon, 1+\epsilon]$ and every $y \in U$, 
   $\phi_{\tZ}^s.\phi_{\tX}^t.y$ is defined for every $s \in [-\epsilon, 1+ \epsilon]$.  Lemma \ref{lem.coincide} says
    that identity $\phi_{\tX}^t.y=\gamma_1.y$  holds on $U$, because it holds on $U \cap \tFd(x)$.
      It follows easily from the property $\gamma_1^*\tX=\tX$ that for 
      every $y \in U$, $\phi_{\tX}^t.y$ is defined for 
      $t \in \RR$.  Relation $\gamma_1^*\tZ=\tZ$  now implies that $\phi_{\tZ}^s.\phi_{\tX}^t.y$ makes sense for 
       every $t \in \RR$,  $s \in [- \epsilon, 1+\epsilon]$, and $y \in U$.  Let us call 
       $\mathcal{U}=\{  \phi_{\tX}^t.y \ | \ t \in \RR, y \in U\}$.  This is an open set on which $\phi_{\tZ}^1$ 
       is defined.  Relation $\phi_{\tZ}^1=\gamma_2$ holds on ${\mathcal U} \cap \tFd(x)$, hence on ${\mathcal U}$ by 
       Lemma \ref{lem.coincide}.  Together with the property $\gamma_2^* \tZ=\tZ$, this implies that 
       $\phi_{\tZ}^s.\phi_{\tX}^t.y$ makes sense for every $y \in U$, and $s,t \in \RR$.  
       
       Now, Lemma \ref{lem.feuillecocompacte} says that $F_{\Delta}(\delta(x))\not = F_{\Delta}(p_0)$.  If $U$ was 
       chosen  small enough, $\delta(y) \not \in F_{\Delta}(p_0)$ for every $y \in U$.  It follows that
        $(t,s) \mapsto e^{sZ}.e^{tX}.\delta(y)$  is a diffeomorphic parametrization of $F_{\Delta}(\delta(y))$.
          In other words, for every $y \in U$, $\{ \phi_{\tZ}^s.\phi_{\tX}^t.y \ | \ (s,t) \in \RR^2  \}$
          coincides with the leaf $\tFd(y)$, and the developing map
          $\delta: \tFd(y) \to F_{\Delta}(\delta(y))$ is a diffeomorphism. 
          Completeness of vector fields $\tX,\tY,\tZ$  on $\tFd(y)$  follows, because vector fields $X,Y,Z$ are complete on leaves 
          of $\calfd$.
          
          Moreover, Lemma \ref{lem.coincide} says that relations $\phi_{\tX}^1=\gamma_1, \ \phi_{\tZ}^1=\gamma_2$ and 
          $\phi_{\tY}^1=\tth$  hold on 
          ${\mathcal U}$ because they hold on $\mathcal{U} \cap \tFd(x)$.  In particular, for 
          $y \in U$, the leaf $\tFd(y)$ is stable by $\gamma_1, \ \gamma_2$ and $\tth$, hence by $\tilde{H}$.

 To conclude that $U \subset {\mathcal E}$, it remains to check that $\tFd(y)$ coincides with $\tF(y)$ for every $y \in U$. 
 A first observation is that $\tFd(y)$ is diffeomorphic to $F_{\Delta}(\delta(y))$,  hence is a disk.  It follows from a 
 cohomological dimension argument that because $\Gamma_{0} \simeq \ZZ^2$, the quotient $\Sigma=\Gamma_{0} \backslash 
 \tFd(y)$ is a torus 
 in $M$.
 Recall
   from (\ref{eq.vectorfields}) the relation $\tth^*\tZ_0=\gamma_1^*\tZ_0=\gamma_2^* \tZ_0=\tZ_0$. Remember also that $\tZ_0$
    is a linear combination of $\tX$ and $\tZ$, hence tangent to $\tFd(y)$. 
  On the torus  $\Sigma$, $\tZ_0$ thus induces
  a  vector field $\overline{Z}_0$ which is $h$-invariant. Notice that $\overline{Z}_0$ is lightlike because $Z_0$ is lightlike 
  on $\ein$, and  $\overline{Z}_0$ is nonsingular 
  because
    the singularities of $Z_0$ are exactly the points of $\Delta$, and $\tFd(y) \subset \tilde{\Omega}_{\Delta}$.
    Hence, for every $z \in \Sigma$, $\overline{Z}_0(z)^{\perp}=T_z \Sigma$.  On the other hand, 
    equality $D_zh^{n_k}(\overline{Z}_0(z))=\overline{Z}_0(h^{n_k}(z))$ shows that $\overline{Z}_0(z)$ belongs to 
  the aproximately stable distribution of $h^{n_k}$ (see the definition of this distribution 
  in Section \ref{sec.approximately}). The aproximately stable distribution has codimension $1$ and is lightlike, so that 
   it coincides with $\overline{Z}_0(z)^{\perp}=T_z\Sigma$ for all $z \in \Sigma$. We conclude that $\Sigma$ is a leaf 
   of $\calf$, what proves $\tF(y)=\tFd(y)$.

\subsubsection{The set $\mathcal{E}$ is closed}
\label{sec.closed}

 We consider a sequence $(x_k)$ of ${\mathcal E}$ converging to $x_{\infty} \in \tm$.
  The leaf $\tF(x_{\infty})$ is accumulated by the sequence of leaves
  $\tF(x_k)=\tF_{\Delta}(x_k)$.  In particular, the vector fields $\tX,\tY,\tZ$ being  tangent to $\tF_{\Delta}(x_k)$ for all 
  $k$, they  are also tangent to 
  $\tF(x_{\infty})$.  Point $2)$ of Fact \ref{fact.action} then says that 
  $\tF(x_{\infty}) \setminus \tilde{\Delta}$ is a union of leaves of $\tcalfd$.  If the set 
  $\tF(x_{\infty}) \cap \tilde{\Delta}$ is not empty, those leaves of $\tcalfd$ might be prolongated smoothly accros the singular set 
  $\tilde{\Delta}$, a contradiction.  We infer that $\tF(x_{\infty}) \subset \tilde{\Omega}_{\Delta}$, and 
    $\tF(x_{\infty})=\tF_{\Delta}(x_{\infty})$.
 
 The union of the  compact  leaves of $\calf$ is a compact subset of  $M$ (see \cite[Chap II, Corollary 3.10]{godbillon}). 
 Since $\calf$ has no vanishing cycles,  $\tF(x_{\infty})$ is left invariant by a discrete 
 subgroup  $\Lambda_1 \subset \pi_1(M)$ which is isomorphic to $\ZZ^2$. 
 We choose $I \subset M$, a small transversal to the foliation $\calf$ containing the point $\pi(x_{\infty})$.
 Following the loops of $F(\pi(x_{\infty}))$ defining $\Lambda_1$ in the neighboring leaves, we 
 get  a corresponding holonomy morphism:
 $$ \operatorname{hol}: \Lambda_1 \to Diff(I).$$
 If $\gamma \in \Lambda_1$, and if for $z \in I$ we  have ${\operatorname{hol}}(\gamma^2).z \not = z$, then the pseudo-orbit 
 ${\operatorname{hol}}(\gamma^n).z$ is infinite and the 
  leaf $F(z)$ cannot be closed.  Because $F(\pi(x_k))$ is a torus for each $k \in \NN$, it follows that replacing $\Lambda_1$ by some index $2$ subgroup, we may assume that 
   $\Lambda_1$ leaves invariant $\tF(x_k)$  for $k$ large enough.
   
 This property also  shows that $\rho(\Lambda_1)$, which is {\it a priori} not a subgroup of $G_{\Delta}$,
 leaves  invariant infinitely many leaves of $\calfd$, hence infinitely many points on $\Delta$.  It follows that
  $\rho(\Lambda_1)$ leaves $\Delta$ invariant: $\rho(\Lambda_1) \subset G_{\Delta}$. A M\"obius transformation 
  fixing infinitely many points on the circle must be trivial, hence 
 $\rho(\Lambda_1) \subset S_{\Delta}$.  
 
 The group generated by $\Lambda_1$ and $\Gamma_0$ is included in $\pi_1(M)$, hence
 must act properly discontinuously on each $\tF(x_k)$.  It
  follows that $\Lambda=<\Lambda_1,\Gamma_0>$ is a discrete group isomorphic to $\ZZ^2$.  In particular $\Lambda_1$ commutes with $\Gamma_0$,
   and $\rho(\Lambda_1) \subset T$. We can then apply 
    Lemma \ref{lem.feuillecocompacte} to $\Lambda$, and we get that   
   $\delta$ is a diffeomorphism from $\tF(x_{\infty})=\tFd(x_{\infty})$ to $F_{\Delta}(x_{\infty})$, and 
   $F_{\Delta}(\delta(x_{\infty})) \not = F_{\Delta}(p_0)$.  It follows that $\tX,\ \tY$ and $\tZ$ are complete in 
   restriction to $\tFd(x_{\infty})$.

   Let $y_{\infty} \in \tF(x_{\infty})$, and let $U \subset \tilde{\Omega}_{\Delta}$ be a small neighborhood of 
   $y_{\infty}$. By completeness of $\tX,\ \tY$ and $\tZ$ in restriction to $\tFd(x_{\infty})$, and shrinking $U$ if necessary,
    the local diffeomorphisms  $\phi_{\tX}^1$, $\phi_{\tY}^1$ and $\phi_{\tZ}^1$ are defined on $U$.  For $k$ large, 
    $U \cap \tFd(x_k) \not = \emptyset$, and identities  $\phi_{\tX}^1=\gamma_1$, $\phi_{\tY}^1=\tth$ and 
    $\phi_{\tZ}^1=\gamma_2$ hold on $U \cap \tFd(x_k)$.  Lemma \ref{lem.coincide} says that those identities hold on $U$.
      Finally $y_{\infty}$ was arbitrary in $\tFd(x_{\infty})$ so that these identities hold on $\tFd(x_{\infty})$. This proves
      $x_{\infty} \in {\mathcal E}$.
        
  

\subsection{Proof of Theorem \ref{th.theoplat}}
\label{sec.conclusion}

Let us draw further conclusions from Proposition \ref{prop.ensembleE}.  The coincidence of the foliations 
$\tcalf$ and $\tcalfd$ implies that $\tcalfd$ is $\pi_1(M)$-invariant.  Moreover, the
$\Gamma_0$-invariance of each leaf $\tFd$, together with Lemma \ref{lem.feuillecocompacte} implies that $\delta$ is injective
 on each leaf $\tFd$, and that $\delta(\tm) \subset \Omega_{\Delta} \setminus F_{\Delta}(p_0)$.  
 
 Also, it follows from Proposition \ref{prop.ensembleE} that $\Gamma_0$ is exactly the subgroup of $\pi_1(M)$ leaving each
  leaf of $\tcalf$ invariant.  It follows that $\Gamma_0$ is normal in $\pi_1(M)$.  We claim that $\Gamma_0$ is also
   normalized by $N_{\pi_1}$, the normalizer of $\pi_1(M)$ in 
  $\Iso(\tm,\tilde{g})$.   Indeed, if $f \in N_{\pi_1}$, then $f\Gamma_0 f^{-1}$ leaves each leaf of $f(\tcalf)$ invariant.
   Now $f(\tcalf)$ is a lightlike,  totally geodesic, codimension $1$ foliation.  Remark \ref{rem.codimension1} 
   ensures that
    on any 
   non locally homogeneous component $\tcalm$, $f(\tcalf)$ coincides with $\tcalf$.  In particular, $f\Gamma_0 f^{-1}$ 
   coincide with $\Gamma_0$ on $\tcalm$, hence $f\Gamma_0f^{-1}=\Gamma_0$.

  The group  $\rho(N_{\pi_1})$ normalizes 
    $\rho(\Gamma_0)$, hence $T$ since $\rho(\Gamma_0)$ is Zariski-dense in $T$. By fact \ref{fact.action}, 
    the lightcone $C(p_0)$ can be characterized as the set of points where the orbits of $\rho(\Gamma_0)$ 
    are contained 
      in a photon of $\ein$. It follows that $C(p_0)$ is left invariant by 
     $\Nor(T)$, the normalizer of $T$ in $\PO$.   Applying the stereographic projection $\varphi$ of pole $p_0$  (see Section \ref{sec.stereographic})
       we can see $\rho(N_{\pi_1})$ as a subgroup of $\operatorname{Conf}(\RR^{1,2})$.  We then show:
       
  \begin{lemme}
   \label{lem.normalisateur}
  Seen in $\operatorname{Conf}(\RR^{1,2})$, the elements of  $\rho(N_{\pi_1})$ are included in the group
  
$$G:=\left\lbrace \left(  \begin{array}{ccc} 
  \epsilon_1&- \epsilon_1 y&-\frac{\epsilon_1}{2 }y^2\\ 0&\epsilon_2 & \epsilon_2 y\\0&0& \epsilon_1\\ 
  \end{array}\right) + 
  \left(  \begin{array}{c} z\\ x\\ t \\ \end{array}\right), \ \ x,z,y,t \in \RR, \ \epsilon_i=\pm 1 \right\rbrace.$$
  
  In particular, we have the inclusion $\rho(N_{\pi_1}) \subset \Iso(\RR^{1,2})$.
  \end{lemme}
     
  \begin{preuve}
   After performing the stereographic projection $\varphi$, the foliation $\calf_{\Delta}$ restricted 
   to $\ein \setminus C(p_0)$ becomes a foliation of $\RR^{1,2}$.  Formula (\ref{eq.stereographic}) for $\varphi$ readily shows 
   that this is the foliation by affine planes 
   of direction $\operatorname{Span}(e_1,e_2)$.  Recall  (see Section \ref{sec.symetries}) that the group $T$ corresponds to the 
   group of translations of vectors $v \in \operatorname{Span}(e_1,e_2)$.
    Since $\Nor(T)$, hence $\rho(N_{\pi_1})$, must preserve 
    this foliation (this is a consequence of  Fact \ref{fact.action}), we see that elements of
     $\rho(N_{\pi_1})$ belong to the subgroup  $G' \subset  \operatorname{Conf}(\RR^{1,2})$ comprising all 
     elements of the form:
     \begin{equation}
   \label{eq.formuleG}
     \left(  \begin{array}{ccc} 
  \lambda \mu&-\lambda \mu y&-\frac{\lambda \mu}{2}y^2\\ 0&\mu & \frac{\mu}{\mu}y\\0&0&\frac{\mu}{\lambda}\\ 
  \end{array}\right) + 
  \left(  \begin{array}{c} z\\ x\\ t \\ \end{array}\right), \ \ x,y,z,t \in \RR \ \ \lambda, \mu \in \RR^*.
  \end{equation}
  
   If a matrix $  \left(  \begin{array}{ccc} 
  \lambda \mu&-\lambda \mu y&-\frac{\lambda \mu}{2}y^2\\ 0&\mu & \frac{\mu}{\mu}y\\0&0&\frac{\mu}{\lambda}\\ 
  \end{array}\right)$  preserves a lattice in $\operatorname{Span}(e_1,e_2)$, then the determinant of
   its restriction to $\operatorname{Span}(e_1,e_2)$ is $\pm 1$.  It follows that $\mu=\pm \frac{1}{\sqrt{|\lambda|}}$.
   
   We saw that   $\rho(\tth)$ belongs to the group $N$, hence has the form:
    $$ \rho(\tth)= \left(  \begin{array}{ccc} 
  1&- y&-\frac{ y^2}{2 }\\ 0&1 & y\\0&0& 1\\ 
  \end{array}\right) + 
  \left(  \begin{array}{c} z\\ x\\ 0 \\ \end{array}\right), \ y \not = 0.$$
  
  In particular, because $\tth$ normalizes $\Gamma_0$, if $\tau \in \rho(\Gamma_0)$, and if we see $\tau$ as a translation of vector 
  $v \in \operatorname{Span}(e_1,e_2)$, then $\rho(\Gamma_0)$ will also contain $v'=v-A.v$, where $A=\left( 
  \begin{array}{cc} 1&y\\
  0&1 \end{array}\right)$.  In other words $\rho(\Gamma_0)$ contains a  translation of vector $\alpha e_1$, 
  $\alpha \not =0$.  The fact that $\rho(N_{\pi_1})$ normalizes the discrete group $\rho(\Gamma_0)$, leads to
   the relation $\lambda \mu =\pm 1$ in (\ref{eq.formuleG}). Together
   with the relation $\mu= \pm \frac{1}{\sqrt{|\lambda|}}$, this  leads to $|\mu|=|\lambda|=1$, and the Lemma follows.
    
  \end{preuve}
    
   Lemma \ref{lem.normalisateur} says that our $(\ein,\PO)$-structure is actually a $(\RR^{1,2},\Iso(\RR^{1,2}))$-structure.
     We conclude that there exists $g'$ in the conformal class of $g$ which is flat, and which is preserved by 
     $\Iso(M,g)$.  We can thus apply the results of Section \ref{sec.minkowski}. Theorem \ref{thm.bieberbach.minkowski}
      and Proposition \ref{prop.topologie.minkowski} say that $(M,g')$ is the quotient of $\RR^3$, $\Heis$ or $\Sol$
       by a lattice.  But  $\rho(\pi_1(M)) \subset G$ by Lemma \ref{lem.normalisateur}, and $G$ does not contain any subgroup
        isomorphic to $\Sol$.  We thus get that $M$ is homeomorphic to $\TT^3$ or to a torus bundle $\TT_A^3$ with 
        $A \subset \SL(2,\ZZ)$ parabolic.  This proves points $1)$ and $2)$ of Theorem \ref{th.theoplat}.

 Finally, Carri\`ere's completeness result \cite{carriere} says that $\delta: \tm \to \RR^{1,2}$ is a conformal 
 diffeomorphism.  It follows that if the coordinates associated to $(e_1,e_2,e_3)$ in $\RR^{1,2}$ are $(u,t,v)$,
  the metric $\tilde{g}$ is of the form
 $$ a(u,t,v)(dt^2+2dudv).$$

 It remains to check that the function $a$ depends only on $v$.  First, the foliation by planes with direction
  $\operatorname{Span}(e_1,e_2)$ is totally geodesic. If $\tilde{\nabla}$ denotes the Levi-Civita
   connection of $\tilde{g}$, we thus have: 
   $$0=\tilde{g}(\tilde{\nabla}_{\partial_t}\partial_t,\partial_u)=
   -\frac{1}{2}\partial_u.\tilde{g}(\partial_t,\partial_t)=
   -\frac{1}{2}\frac{\partial a}{\partial u}.$$
   
   Identifying $\tth$ and $\rho(\tth)$, we saw that $$ \tth= \left(  \begin{array}{ccc} 
  1&- y&-\frac{ y^2}{2 }\\ 0&1 & y\\0&0& 1\\ 
  \end{array}\right) + 
  \left(  \begin{array}{c} \mu_0 \\ \nu_0 \\ 0 \\ \end{array}\right)$$
  where $y \not =0$. 

  It follows that $\rho(\tth)$ acts on each hyperplane $v=v_0$ by the affine transformation:
  $$ \left( \begin{array}{c} u\\ t  \end{array} \right) \mapsto \left( \begin{array}{c} u-yt + \mu(v_0)\\ t+ \nu(v_0)
  \end{array}  \right),$$
where $\mu(v_0)=\mu_0 - \frac{y^2}{2}v_0$ and $\nu(v_0)=\nu_0 + y v_0$.

The group $\Gamma_0$ is generated by two translations $\tau_1,\tau_2$ of (linearly independant) vectors $ \left( \begin{array}{c} a\\ b  \end{array} \right)$
 and $ \left( \begin{array}{c} c\\ d  \end{array} \right)$  respectively.
 The $w$-coordinate of $\tth^k\circ \tau_1^m \circ \tau_2 ^n \left( \begin{array}{c} u\\ t  \end{array} \right)$
 is $t+k \nu(v_0)+mb+nd$.  Because $\tth^k\circ \tau_1^m \circ \tau_2 ^n$ acts isometrically for $\tilde{g}$
  this leads to $a(t,v)=a(t+k \nu(v)+mb+nd,v)$ for every $(k,m,n) \in \ZZ^3$.  Since $b$ and $d$ can not be both zero 
  (let say $b \not =0$), and because $\nu(v)$ and $b$ are rationally independant for almost every value of $w$ (because $y \not =0$), we get that
   for almost every $w$, $t \mapsto a(t,v)$ is constant.  As a consequence, $a=a(v)$, and the fact that it is a periodic function follows easily from the
    compactness of $M$.  
    
    Finally the group $N$ which comprises transformations of the form 
    $$  \left(  \begin{array}{ccc} 
  1&- y&-\frac{ y^2}{2 }\\ 0&1 & y\\0&0& 1\\ 
  \end{array}\right) + 
  \left(  \begin{array}{c} z\\ x\\ 0 \\ \end{array}\right), \ \ x,y,z \in \RR$$
  is isomorphic to $\Heis$ and acts isometrically on $(\tm,\tilde{g})$.
    This concludes the proof of Theorem \ref{th.theoplat}.

\section{Conclusions}
\label{sec.fin}

The study made in Section \ref{sec.homogeneous}, as well as Theorems \ref{theo.hyperbolic} and \ref{th.theoplat} provide
 all possible topologies for a closed $3$-dimensional, orientable and time-orientable, Lorentz manifold with 
  a noncompact isometry group.  Those are the $3$-dimensional torus, hyperbolic or parabolic torus bundles, and compact quotients
   $\Gamma \backslash \widetilde{\PSL}(2,\RR)$. Together with the examples provided in Section \ref{sec.panorama}, this yields Theorem 
   \ref{thm.topological}.
   
   Let us now look at the  geometries which can occur on those manifolds, and prove Theorem \ref{thm.geometrical}.  
    The manifolds $\Gamma \backslash \widetilde{\PSL}(2,\RR)$ occur only in Proposition \ref{prop.topologieads}.  Hence
     the only metrics on such manifolds which admit a noncompact isometry group are covered by
   $\widetilde{\PSL}(2,\RR)$, endowed with  a Lorentzian, 
     non-Riemannian, left-invariant metric.  In particular those manifolds $(M,g)$ are locally homogeneous and 
     $(\tm,\tilde{g})$ admits an isometric action of  $\widetilde{\PSL}(2,\RR)$. 
     
   Parabolic torus bundles appear in Proposition \ref{prop.topologie.minkowski} and Theorem \ref{th.theoplat}.
     We saw there that the universal cover is isometric to $\RR^3$ endowed with a metric
     $$ a(v)(dt^2+2dudv),$$
      with $a$ smooth and periodic. This universal cover admits an isometric action of $\Heis$.
      If the manifold $(M,g)$ is locally homogeneous, Proposition  
      \ref{prop.topologie.minkowski}  ensures that $g$ is flat or locally isometric to the Lorentz-Heisenberg metric.
      
  Hyperbolic torus bundles appear only in Proposition \ref{prop.topologie.minkowski} 
   Theorem \ref{theo.hyperbolic}.  We saw that the universal cover is isometric with $\RR^3$ endowed
    with a metric $dt^2+2a(t)dudv$, with $a$ smooth and periodic.  There is an isometric action of $\Sol$ on this universal cover.
     The manifold $(M,g)$ is locally homogeneous if and only if it is flat.
     
   Finally, $3$-tori appear in  Proposition \ref{prop.topologie.minkowski} 
   Theorem \ref{theo.hyperbolic} and Theorem \ref{th.theoplat}.  The metric on the universal cover $\tm$ is provided
    by those two last theorems, and there is always an isometric action of $\Heis$ or $\Sol$ on $(\tm,\tilde{g})$.
     Finally, $(M,g)$ is locally homogeneous if and only if it is flat.

     Those results alltogether prove Theorem \ref{thm.geometrical} and Corollary \ref{coro.action.universel}.

\section{Annex A: Some computations}

We present here the necessary computations leading to Proposition \ref{prop.heis}.  

\subsubsection{The curvature module}

We consider on  $\RR^3$ the Lorentzian form, with matrix in a   basis $e,h,f$ given by 
 $ J=\left(  \begin{array}{ccc}    0&0&1\\ 0&1&0 \\ 1&0&0\\ \end{array}  \right).$
 
We call ${\operatorname O}(1,2)$ the subgroup of  $\GL(3,\RR)$ preserving the bilinear form determined by $J$. 
Its Lie algebra is denoted by   $\oo(1,2)$, and admits the following basis :
$$ E=\left(  \begin{array}{ccc}   0&1&0\\ 0&0& -1\\ 0&0&0\\ \end{array}  \right),  
H= \left(  \begin{array}{ccc} 1&0&0\\ 0&0&0 \\ 0& 0& -1   \end{array}  \right), 
 F= \left(  \begin{array}{ccc}   0&0&0\\ 1&0&0 \\ 0&-1& 0 \end{array}  \right).$$
We thus have the commutation relations $[H,E]=E, [H,F]=-F$ and $[E,F]=H$.

Let $(M,g)$ be  $3$-dimensional Lorentz manifold, and denote by  $\hm$ its bundle of orthonormal frames.  At each $\hx \in \hm$, the 
 curvature $\kappa(\hx)$ is an element of $\Hom(\wedge^2(\RR^3),\oo(1,2))$. Because of Bianchi's identities,
 the curvature module is actually a $6$-dimensional
  submodule of $\Hom(\wedge^2(\RR^3),\oo(1,2))$.
Choosing $e \wedge h$, $e \wedge f$, $h \wedge f$ as a basis for $\wedge^2(\RR^3)$, and $E,H,F$ as a basis for $\oo(1,2)$, 
an element of $\Hom(\wedge^2(\RR^3),\oo(1,2))$ is merely given by a $3 \times 3$ matrix.   The action of ${\operatorname O}(1,2)$ on 
$\Hom(\wedge^2(\RR^3),\oo(1,2))$ corresponds to  the conjugation on matrices.

Scalar matrices are ${\operatorname O}(1,2)$-invariant, and form a $1$-dimensional irreducible submodule 
(corresponding to constant sectional curvature).

The other irreducible submodule of the curvature module is $5$-dimensional, spanned by the matrices:

$$ \left(  \begin{array}{ccc}   0&0&1\\ 0&0&0\\ 0&0&0\\  \end{array}  \right), 
 \left(  \begin{array}{ccc}   0&1&0\\ 0&0&1\\ 0&0&0\\ \end{array}  \right), 
 \left(  \begin{array}{ccc}  1&0&0\\ 0&-2 & 0\\0& 0& 1 \\  \end{array}  \right),
 \left(  \begin{array}{ccc}  0&0&0\\ 1&0 &0 \\ 0 & 1 & 0\\  \end{array}  \right), \left(  \begin{array}{ccc}    
0&0&0 \\ 0&0&0 \\1&0&0 \\ \end{array}  \right).$$

We call $\kappa_0$ the element of $\Hom(\wedge^2(\RR^3),\oo(1,2))$  corresponding to the identity matrix, namely
$\kappa_0$ maps $e \wedge h$ to $E$, $e \wedge f$ to $H$ and  $h \wedge f$ to $F$.  We also call $\kappa_1$ 
the element of  $\Hom(\wedge^2(\RR^3),\oo(1,2))$ 
corresponding to the matrix $\left(  \begin{array}{ccc}   0&0&1\\ 0&0&0\\ 0&0&0\\  \end{array}  \right)$.  

The two dimensional vector space spanned by $\kappa_0$ and $\kappa_1$ is the set of fixed points of the action of $\{ e^{tE} \}_{t \in \RR}$ on 
the curvature module.

\subsubsection{Identification of the $\kiloc$-algebra} 
We consider a parabolic component $\calm$ which
 is not locally homogeneous.  In such a component, the points are either parabolic, or 
 points where the isotropy algebra is 
 $3$-dimensional and  the sectional curvature is constant.  The set of parabolic points  is thus 
 a dense open set $\Omega \subset \calm$.  Observe that at a parabolic point $x \in \Omega$, if $X$ 
a local Killing field around $x$, generating the isotropy ${\mathfrak Is}(x)$, the
 $1$-parameter group $D_x\varphi_X^t$ is unipotent in $\operatorname{O}(T_xM)$. 
  In a suitable basis $(u_1,u_2,u_3)$ of $T_xM$ satisfying
  $g(u_1,u_3)=1=g(u_2,u_2)$ and all the other products are $0$, the matrix of $D_x\varphi_X^t$ reads
  $$\left( \begin{array}{ccc}
    1&t&-t^2/2\\
  0&1&-t\\
   0&0&1\\

      \end{array}
      \right)
 $$
 
 We quickly check that the only $2$-plane stable by $D_x\varphi_X^t$ is spanned by $u_1$ and $u_2$, so that on $\Omega$, the
  $\kiloc$-orbits must be lightlike surfaces.

Let us now fix  a point $x \in \Omega$. 
%
%
 We 
work in the fiber bundle $\hm$  (and lift all local Killing fields there). After multiplying 
$X$ by a suitable constant, we can
find $\hx \in \hm$ in the fiber of $x$ such that 
 $\omega(X(\hx))=E$.  We now choose $Z$ and $Y$  two local Killing fields around $x$
   such that $Z(x)=u_1$ and $Y(x)=u_2$.
 After adding to $Z$ and $Y$ a suitable multiple of $X$, we can write, at $\hx$:

$$ \omega(Z)=e+\beta H + \gamma F \ \text{   and     }\   \omega(Y)=h+\alpha H +\nu F.$$
The curvature  $\kappa(\hx)$ is $\Ad(e^{tE})$-invariant, hence is of the form $\kappa=\sigma \kappa_0 +b \kappa_1$. 
In particular, the following identities hold at $\hx$:
\begin{equation}
  \label{eq.courbure}
  \kappa(e \wedge h)= \sigma E, \ \kappa(e \wedge f)=\sigma H, \ \kappa(h \wedge f)=bE+\sigma F.
 \end{equation}
Notice that  $\sigma$, $b$,
$\alpha$,$\beta$,$\gamma$,$\nu$ depend on $x$ and  $\hx$, but since those points are fixed, 
there will be considered as constant in the sequel.  

 Cartan's formula $L_X \omega= \iota_X  d \omega + d(\iota_X \omega)$ shows that whenever $U,V$ are 
 two Killing fields on $\hm$,  the following relation holds:

\begin{equation}
\label{courbure.killing}
\omega([U,V])=K(U,V)-[\omega(U),\omega(V)]
\end{equation}

In the sequel, we will call  ${\mathcal H}$ the span of $\omega(Z),\omega(X),\omega(Y)$ at $\hx$, and we are going
 to write Equation (\ref{courbure.killing}) at $\hx$, using identities (\ref{eq.courbure}), when $U$ and $V$ range 
 over $Z$, $X$, $Y$.  
  For instance, the first equation is:
\begin{eqnarray}
\omega([Z,X])& = & -[\omega(Z),\omega(X)] =  -[e+ \beta H + \gamma F,E] \nonumber\\
  &=& - \beta E + \gamma H =  \omega(- \beta X) + \gamma H. \nonumber
  \end{eqnarray}

The fact that
 $\kiloc(x)$ is a Lie algebra, together with the property 
$H \not \in {\mathcal H}$ forces  $\gamma$ to vanish. Next, two Killing
fields which coincide at $\hx$ must be equal (by freeness of the action of isometries on the orthonormal frames),
 which implies $[Z,X]= - \beta X$. To summarize: 
\begin{equation}
\label{eq.zx}
\gamma =0\  \text{ and }\  [Z,X]= - \beta X.
\end{equation}

We proceed exactly in the same way for the two other equations:

 \begin{eqnarray}
\omega([X,Y])& = &-[\omega(X),\omega(Y)]=-[E, h + \alpha H + \nu F] \nonumber\\
 &=& - e + \alpha E - \nu H= \omega(-Z +  \alpha X) +(\beta - \nu) H. \nonumber
 \end{eqnarray}
leads to: 
\begin{equation}
\label{eq.xy}
\beta=\nu \ \text{ and }\  [X,Y]=-Z+ \alpha X.
\end{equation}
Finally
\begin{eqnarray}
\omega([Z,Y])& = & \ka(e,h)-[e+ \nu H, h + \alpha H + \nu F]=\sigma E + \alpha e + \nu h + \nu^2 F \nonumber\\
  &=&  \omega(\alpha Z + \sigma X + \nu Y) -2 \alpha \nu H. \nonumber
\end{eqnarray}
 implies: 
\begin{equation}
\label{eq.zy}
\alpha \nu =0 \text{ and } [Z,Y]= \alpha Z + \sigma X + \nu Y.
\end{equation}

Notice that establishing (\ref{eq.zx}) and (\ref{eq.xy}), we have actually shown that $\ad(X)$ is a nilpotent endomorphism of 
$\kiloc(x)$.  This property did not use anything special on $x$, so that we actually have:
\begin{fact}
    \label{fait.nilpotent}
    At each $z \in \Omega$, if $U$ is a local Killing field around $z$ generating the isotropy at $z$, then
    $\operatorname{ad}(U)$ is a nilpotent endomorphism of $\kiloc(z)$.
   \end{fact}

%
%
%
%
%
%


At $x$, $Z(x)$ is lightlike and nonzero and $Y(x)$ is spacelike, orthogonal to $Z(x)$.  
 The orthogonal to  $Y(x)$ 
at $x$  is a Lorentzian plane spanned by $Z(x)$ and another vector $w \in T_xM$.  
 Let us call $t \mapsto \gamma(t)$    the geodesic through $x$ satisfying 
 $\dot{\gamma}(0)=w$. Clairault's equation ensures that the quantities 
 $g(\dot{\gamma}(t),Z(\gamma(t)))$, $g(\dot{\gamma}(t),X(\gamma(t)))$ and 
 $g(\dot{\gamma}(t),Y(\gamma(t)))$ do not depend on $t$. In particular, for
  $t>0$, both $Y(\gamma(t))$ and $X(\gamma(t))$ are orthogonal to $\dot{\gamma}(t)$ 
  while  $Z(\gamma(t))$ is not. 
 For $t>0$ small enough, $Y(\gamma(t))$ is still spacelike, hence nonzero, and $\gamma(t)$ belongs to $\Omega$.
  In particular, the $\kiloc$-orbit at $\gamma(t)$ is $2$-dimensional, so that 
 $Y$ and $X$  must be colinear at $\gamma(t)$.   One then has 
 $X(\gamma(t))=\lambda_t Y(\gamma(t))$, for some real $\lambda_t$. Observe finally that $w$ is not fixed by $D_x\phi_X^t$, hence is
  transverse to the set where  $X$ vanishes.   In particular, for $t \geq 0$ small, $X(\gamma(t))=0$ only for $t = 0$,
  and thus
   $\lambda_t \not =0$ if $t \not =0$.

We claim that those considerations lead necessarily to  $\alpha=0$.  Indeed, using the bracket relations 
(\ref{eq.xy}),(\ref{eq.zy}) and (\ref{eq.zx}), we compute
$$ \operatorname{Trace}(\ad(\lambda_t Y - X))=-2 \lambda_t \alpha.$$
For $t \geq 0$ small, $X,Y,\ Z$ generate $\kiloc(\gamma(t))$, hence $-2 \lambda_t \alpha=0$ because of 
Fact \ref{fait.nilpotent}. Since 
 $\lambda_t \not = 0$ if 
$t \not =0$, we get $\alpha=0$.  Injecting this data in equation (\ref{eq.zy}) and (\ref{eq.zx}), we find that the matrix
of $\ad(\lambda_t Y - X)$ in the basis $Y,Z,X$ is:
  
  $$ \left(\begin{array}{ccc}
  
  0 & -\lambda_t \nu &0\\
  1 & 0 & \lambda_t\\
  0 & -\lambda_t \sigma -\nu & 0\\
  
  \end{array} \right).$$
The characteristic polynomial of  $\ad(\lambda_t Y - X)$ is
$$ Q(x)=-x^3-\lambda_t x(\lambda_t \sigma + 2 \nu). $$
Hence, the nilpotency of $\ad(\lambda_t Y - X)$ (Fact \ref{fait.nilpotent})  implies 
\begin{equation}
\label{eq.relation}
\lambda_t \sigma + 2 \nu=0, 
\end{equation}
If $\sigma \not =0$, we get that $t \mapsto \lambda_t$ is constant, which is not the case since we observed that 
$\lambda_0=0$ but $\lambda_t \not =0$ for $t>0$ small.  
  We end up with the equality $\sigma=\nu=0$.   The vector fields
  $Z,X,Y$ then satisfy the bracket relations:
  $$ [Z,X]=0=[Z,Y], \text{ and } [X,Y]=-Z,$$
   showing that  Lie algebra $\kiloc(x)$ is  isomorphic to $\heis$. We also proved that $\sigma$, the scalar curvature at
   $x$, vanishes, but since $x$ was arbitrary in the open set $\Omega$, 
   we finally get the vanishing of the scalar curvature 
   on $\Omega$, and then on $\calm$ by density.

\subsubsection{Description of the  $\kiloc$-orbits}  
 The fact that the local Killing algebra is isomorphic to $\heis$ shows that no point in $\calm$ has a $3$-dimensional isotropy algebra.
  Indeed, the isotropy representation at those points would yield an embedding $\heis \to \oo(1,2)$, what is impossible.
   We thus get $\Omega=\calm$, and all the $\kiloc$-orbits on $\calm$ are $2$-dimenional and lightlike.
   
   On the other hand, since the isotropy algebra ${\mathfrak{Is}}(x)$ generates a parabolic $1$-parameter subgroup of 
    $\OO(1,2)$ at each $x$, there is a totally geodesic lightlike hypersurface $F(x)$, whose tangent space is left invariant by 
    the isotropy (see \cite[Lemma 3.5]{sorin.ghani} and its proof).  We already observed that at $x \in \calm$, 
     the local isotropy preserves only one $2$-plane of $T_xM$.  This implies that the $\kiloc$-orbits are everywhere
      tangent to a leaf of a totally geodesic foliation of $\calm$, hence the $\kiloc$-orbits are themselves totally geodesic.
       This concludes the proof of Proposition 
  \ref{prop.heisenberg.flat}.

%
%

\end{document}